\documentclass{article}

\usepackage{amsfonts}

\usepackage{amsmath,theorem}

\usepackage[all]{xy}
\usepackage{graphicx}

\usepackage{amssymb}

\usepackage{latexsym}

\DeclareMathAlphabet\EuFrak{U}{euf}{m}{n}       
\SetMathAlphabet\EuFrak{bold}{U}{euf}{b}{n}     
\newcommand{\halmos}{\ensuremath{\raisebox{-.3ex}{\rule{1ex}{1em}}}}

\newcommand{\ra}{\rightarrow}

\newcommand{\ovl}{\overline}

\newcommand{\sC}{{\it C*}-}
\newcommand{\bC}{{\mathbb C}}

\newcommand{\mA}{\mathcal A}

\newcommand{\mC}{\mathcal C}
\newcommand{\mE}{\mathcal E}

\newcommand{\mU}{\mathcal U}

\newcommand{\Hom}{\mathrm{Hom}}
\newcommand{\End}{\mathrm{End}}

\newcommand{\om}{\omega}

\theoremheaderfont{\scshape}
\newtheorem{defn}{Definition}[section]
\newtheorem{lem}[defn]{Lemma}
\newtheorem{prop}[defn]{Proposition}
\newtheorem{thm}[defn]{Theorem}
\newtheorem{cor}[defn]{Corollary}

\newtheorem{ass}[defn]{Assumption}
\theorembodyfont{\rmfamily}
\newtheorem{rem}[defn]{Remark}
\newtheorem{ex}[defn]{Example}

\numberwithin{equation}{section}

\begin{document}     
\author{ Pasquale A. Zito \\Dipartimento di Matematica\\  Universit\`a di Roma ``Tor
  Vergata''\\ via della Ricerca scientifica 1, I-00133, Roma, Italy\\ e-mail: zito@mat.uniroma2.it}
\title{\bf{2-$C^*$-Categories with non-simple units }\\}
\maketitle

\begin{abstract}
We study the general structure of 2-$C^*$-categories closed under conjugation,
projections and direct sums. We do not assume units to be simple,
i.e.\  for $\iota_A$ the 1-unit corresponding to an object $A$, the space
$\Hom(\iota_A,\iota_A)$ is a commutative unital $C^*$-algebra. We show that
2-arrows can be viewed as continuous sections in Hilbert  bundles and describe
the behaviour of the fibres with respect to the categorical structure. 
We give an example of a 2-$C^*$-Category giving rise to bundles of finite
Hopf-algebras in duality. We make some remarks concerning Frobenius
algebras and $Q$-systems in the general context of tensor $C^*$-categories with non-simple units.    
\end{abstract}

\section*{Introduction}
\label{Introduction}

Tensor categories arise in many mathematical contexts, such as quantum field
theory, topological field theory, the theory of topological invariants for
knots and links, quantum groups. The various examples may be characterised by
more specific properties such as the notions of duality or of spherical
(resp. braided, symmetric, modular) category. These have in common a
``tensor'' structure, i.e.\   a $\otimes$ product defined on  objects
and on arrows and a particular object $\iota$ in the category, the unit with
respect to this tensor product. Tensor $C^*$-categories are naturally related
to the theory of operator algebras, and are at the heart of the structure of
super-selection theory in algebraic quantum field theory and  duality for
compact groups (see  \cite{DR89}, \cite{DR89A}).

A 2-category may be viewed as a generalisation of the notion of  tensor category. 
We recover the notion of a tensor category as the particular case of a 2-category with
only one object, where we can regard  1-arrows as the objects of the tensor category and
2-arrows as the arrows in the tensor category. The horizontal composition gives the tensor
structure.  
2-$C^*$-categories describe much of the structure in the theory of subfactors 
of von Neumann algebras. Many tools used in the past years which have proved to be
fundamental in order to characterise subfactors can be described in the context of
2-$C^*$-categories. Furthermore the language of categories makes easier the link with
other fields of mathematics. 

 In a  2-category for each object $A$ we have a particular 1-arrow $ A
 \xleftarrow{\iota_{A}} A$, a unit with respect to the $\otimes$ composition. In the
 tensor category case this corresponds to the unit object. 
In most of the literature units are supposed to be simple,
i.e.\  $\End(\iota_{A})\equiv\Hom(\iota_{A},\iota_{A}) = \mathbb{K} 1_{\iota_A}$,
the spaces of 2-arrows connecting the unit to itself is a field. 
In this work we drop this hypothesis. More precisely, our object of study will be
2-$C^*$-categories closed under conjugation, projections and direct sums. 
The axioms for the $\otimes$, $\circ$ products imply that the spaces
$\End(\iota_{A})$ are commutative $C^*$-algebras,
i.e.\  $\End(\iota_{A}) \cong C(\Omega_{A})$, the continuous functions over
a compact Hausdorff space $\Omega_{A}$. 
We show that for each couple $B \xleftarrow{\rho} A$, $B \xleftarrow{\sigma}
A$ the space of 2-arrows $\Hom(\rho,\sigma)$ connecting $\rho$ to $\sigma$ may be
furnished with the
structure of a Hilbert (or more generally a Banach) bundle over the
topological space $\Omega_{A}$, derived by the conjugation relations and
depending on the choice of solutions up to continuous linear isomorphisms
which respect the fibre structure (Proposition \ref{bundlea}).
This generalises the well known situation of simple units,
where
for each pair of 1-arrows $\rho,\sigma$ the space of 2-arrows $\Hom(\rho,\sigma)$
is finite dimensional. We show that each element of $\Hom(\rho,\sigma)$ can be
envisaged as a continuous section in this Banach bundle. We study the relation
of this bundle structure with respect to the categorical structure, in
particular  how the fibres behave under $\otimes$ and $\circ$
composition, conjugation and $*$ involution. 
Furthermore, $\End(\rho)$ has the structure of a
$C^*$-algebra bundle, and this implies that even for the case $\rho \neq
\sigma$ it is possible to give $\Hom(\rho,\sigma)$ a Banach bundle structure such
that the norm of each fibre satisfies the $C^*$-property. This fact
implies, for example, that in the particular case of a
$\End(\iota)$-linear tensor $C^*$-category one can describe the
situation as that of a ``bundle of tensor categories with simple units''
(Proposition \ref{fibre of tensor categories}).
We discuss the existence of a class of ``standard'' solutions to the
conjugation equation, a
generalisation of a notion appearing naturally in the simple units case, stable under
composition of arrows, direct sums and projections.



The hypothesis of a non-simple unit has been considered by some authors in related  contexts:

\cite{BL01} and \cite{BL03} deal with the categorical structure arising from the extension
of a $C^*$-algebra with non trivial centre by a Hilbert $C^*$-system; 

\cite{FI} deals with the existence of minimal conditional expectations for the inclusion
of two von Neumann algebras with discrete centres; 

\cite{SZ} deals with the structure of depth two inclusions of finite dimensional $C^*$-algebras.

\cite{KPW03} studies a notion of Jones index and its relation to conjugation in the
context of the 2-$C^*$-category of Hilbert $C^*$-bimodules; 

\cite{V2} deals with crossed products of $C^*$-algebras with nontrivial centres by
endomorphism, $C^*$-algebra bundles, group bundles.

Our work  lives in a more abstract setting, and it  may be viewed as a development of some
of the results of \cite{LR97} for the case $\End(\iota) \ncong \bC $.

The work is organised as follows:

\begin{itemize}
\item In the first part we recall the basic definitions and state some basic
  results concerning the structure of 2-$C^*$-categories. 

\item In the second part we develop further the structure of 2-$C^*$-categories. In
particular the fibre structure appears.

\item In the third  part we  define and study the class of ``standard solutions
  `` to the conjugate equations, giving some conditions for their existence.
\item   In the fourth part we give an example of a 2-category generated by two
  conjugate 1-arrows $\rho, \bar \rho$ giving rise to two continuous bundles of
  finite dimensional Hopf-algebras in duality, a result analogous to the case
  of an irreducible depth two subfactor. 
\item In the fifth part we show, among other things, that  in the tensor $C^*$-category context
  every Frobenius algebra is  equivalent to a $Q$-system (Corollary \ref{frbstrspb}). This fact 
  strengthens a result contained in \cite{Mu} dealing with the equivalence of
  Frobenius algebras and pairs of conjugate 1-arrows in a 2-category, as it
  implies that in the $C^*$-  
  case   some of the hypothesis assumed to hold are automatically verified (Proposition \ref{thmmichaelplus}).

\end{itemize}

\section{Preliminaries}
\label{Preliminaries}

We recall here briefly some basic concepts from category theory (the classical
reference is \cite{ML98}, see also the  introduction of
\cite{Le04} for a quick review).
We will be interested only in small categories and  we will not need more general
definitions which avoid the notion of ``set'' from the beginning.

\begin{defn}
A \emph{category} $\mathcal{C}$ consists of  a set of objects $\mathcal{C}_{0}$ and
a set of arrows $\mathcal{C}_{1}$ together with the following structure:
\begin{itemize}
\item a  \emph{source} map: $\mathcal{C}_1 \rightarrow \mathcal{C}_0$ assigning an
  object $\emph{s}(f)$ to each arrow $f \in \mathcal{C}_1$,

\item a \emph{target} map: $\mathcal{C}_1\rightarrow  \mathcal{C}_0$ assigning an
  object $\emph{t}(f)$ to each arrow $f \in \mathcal{C}_1$,

\item an identity $1$ map $: \mathcal{C}_0 \rightarrow \mathcal{C}_1$
  assigning to each object $a$ an arrow $1_a$ with
  $\emph{s}(1_a) = \emph{t}(1_a) = a$,

\item a composition map $\circ: \mathcal{C}_1 \times \mathcal{C}_1 \rightarrow
  \mathcal{C}_1$ assigning to each pair of arrows $f,g$ s.t. $\emph{s}(g) =
  \emph{t}(f)$ a third arrow $g \circ f $ with $\emph{s}(g \circ f) =
  \emph{s}(f)$ and $\emph{t}(g \circ f) = \emph{t}(g)$.

\item The composition $\circ$ is associative, i.e.\  $h \circ (g \circ f) = (h \circ
g) \circ f$ whenever these compositions make
sense.

\item the identity map satisfies $ f \circ 1_a =  f, \ \ \forall f
\ \ s.t. \ \ \emph{s}(f) = a$ and \\ $ 1_a \circ g = g, \ \ \forall g
\ \ s.t. \ \ \emph{t}(g) = a $.
\end{itemize}
\end{defn}

\begin{rem}
An alternative way of defining a category is that of a set $\mathcal{C}$
together with a set of pairs $(f,g)$ of elements in $\mathcal{C}$, which are
said to be composable, and a composition map $\circ: (f,g ) \mapsto f
\circ g$ assigning to each composable couple an element in $\mathcal{C}$. This
map is supposed to be strictly associative. Elements $u$ such that $f \circ u
= f$ and $u \circ g = g$ for any $f,g$ such that the compositions are defined
are called units. For each $f \in \mathcal{C}$ there exist right and left
units. Thus in this approach one   deals  only with arrows and the objects are identified
with their units .  
\end{rem}

\begin{defn}
  A functor $F: \mathcal{C} \rightarrow \mathcal{C'}$ from the category
  $\mathcal{C}$ to the category $\mathcal{C'}$ consists of a map sending
  objects to objects and arrows to arrows such that:
  \begin{itemize}
   \item for any object $A \in \mathcal{C}_{0}$ we have $F(1_A) = 1_{F(A)}$
   \item for any arrow $ f \in \mathcal{C}_{1}$ we have $\emph{s}(F(f)) =
   F(\emph{s}(f))$ and $\emph{t}(F(f)) = F(\emph{t}(f))$.
   \item $F(g \circ f) = F(g) \circ F(f)$, if $f$ and $g$ are composable.
  \end{itemize}

\end{defn}

\begin{defn}
  Let $F,G: \mathcal{C} \rightarrow \mathcal{C'}$ be functors from the
  category $\mathcal{C}$ to the category $\mathcal{C'}$. A natural
  transformation $\eta: F \rightarrow G$ is a family of arrows 
$ \eta(A):   F(A) \rightarrow G(A),\ A\in\mathcal{C}$ in $\mathcal{C'}$ such that for any
morphism $f: A  \rightarrow B$ in $\mathcal{C}$ the following holds:
 $$ G(f) \circ \eta(A) =  \eta(B) \circ F(f) . $$
\end{defn}
\begin{defn}
  A strict 2-category $\mathcal{B}$ consists of

  \begin{itemize}
  \item a set $\mathcal{B}_{0}$, whose elements are called objects.

  \item a set $\mathcal{B}_{1}$, whose elements are called 1-arrows. A source
  map $\emph{s}: \mathcal{B}_1 \rightarrow \mathcal{B}_0$ and a target map
  $\emph{t}: \mathcal{B}_{1} \rightarrow \mathcal{B}_{0}$ assigning to
  each 1-arrow a source and a target object. 
 
  \item an identity map $\iota: \mathcal{B}_0 \rightarrow \mathcal{B}_1$
  assigning to each object $A$ a 1-arrow $\iota_A$ such that
  $\emph{s}(\iota_A) = \emph{t}(\iota_A) = A$.

  \item a composition map $\otimes: \mathcal{B}_1 \times \mathcal{B}_{1}
  \rightarrow \mathcal{B}_{1} $ making $(\mathcal{B}_{0}, \mathcal{B}_{1})$
  into a category, i.e.\  for any $\rho , \sigma \in \mathcal{B}_{1}$ such that
  $\emph{s}(\sigma) = \emph{t}(\rho)$, \ $\emph{s}(\sigma \otimes
  \rho) = \emph{s}(\rho)$ and $\emph{t}(\sigma \otimes \rho) = \emph{t}(\sigma)$hold.

  \item a set $\mathcal{B}_{2}$, whose elements are called 2-arrows. A source
  map (we use the same symbols as above) $\emph{s}: \mathcal{B}_{2}
  \rightarrow \mathcal{B}_{1}$ and a target map $\emph{t}: \mathcal{B}_{2}
  \rightarrow \mathcal{B}_{1}$ assigning to each 2-arrow $S \in
  \mathcal{B}_{2}$ source and target 1-arrows. The maps $\emph{s}$ and
  $\emph{t}$ satisfy $\emph{s}(\emph{s}(S)) = \emph{s}(\emph{t}(S))$ and
  $\emph{t}(\emph{t}(S)) = \emph{t}(\emph{s}(S))$ for any $S \in \mathcal{B}_{2}$. 

  \item a composition map  (the ``vertical'' composition) $\circ:
  \mathcal{B}_{2} \times \mathcal{B}_{2} \rightarrow \mathcal{B}_{2}$ defined
  for each pair of $S,T $ for $\mathcal{B}_{2}$  such that $\emph{s}(T) = \emph{t}(S)$.

  \item a unit map $1: \mathcal{B}_{1} \rightarrow \mathcal{B}_{2}$ assigning
  to each 1-arrow $\rho$ a 2-arrow $1_{\rho}$ such that $\emph{s}(1_{\rho}) =
  \rho = \emph{t}(1_{\rho})$. 
 
  \item a composition map (the ``horizontal'' composition, for which we use
  the same symbol as above) $ \otimes: \mathcal{B}_{2} \times
  \mathcal{B}_{2} \rightarrow \mathcal{B}_{2}$ defined for each pair of
  2-arrows $S,T$ such that $\emph{s}(\emph{s}(T)) = \emph{t}(\emph{s}(S))$ and
  $\emph{s}(\emph{t}(T)) = \emph{t}(\emph{t}(S))$.

  \item the following equation for elements of $\mathcal{B}_{2}$ holds, whenever the
    compositions make sense 
$$(S   \otimes T) \circ (S' \otimes T') = (S \circ S') \otimes (T \circ T').$$

  \item for each $A \in \mathcal{B}_{0}$  the following equations hold,
  whenever the composition with elements $S,T \in \mathcal{B}_{2}$ make sense
  $$T \otimes 1_{\iota_{A}} = T, \ \ 1_{\iota_{A}} \otimes S = S. $$
\end{itemize}
\end{defn}


Note that for any pair of objects $A,B$, we have a category $HOM(A,B)$, with
the set of  1-arrows with source $A$ and target $B$  as objects of $HOM(A,B)$
and the set of 2-arrows with the latter as source and target as arrows and
$\circ$ as composition and $1$ as unit map.   
Note that each object $A$ determines a unit 1-arrow $\iota_A$, which in turn
determines a unit 2-arrow $1_{\iota_A}$. 

\begin{rem}
Concerning notation, we will give precedence, when not specified otherwise,  to the
$\otimes$ product over
the $\circ$ product, i.e.\  $S \otimes T \circ U$ has to be read $(S \otimes T)
\circ U$ and not $S \otimes (T \circ U)$.
\end{rem}

\begin{rem}
We will often write $\End(\rho)$ instead of $\Hom(\rho,\rho)$. 
\end{rem}

\begin{lem}
The set of 2-arrows $\End(\iota_A)$ is a commutative monoid.   
\end{lem}
{\noindent{\it Proof.}}  In fact  for any $w,z \in \End(\iota_A)$ we have $ w \otimes z =
(1_{\iota_A} \circ w) \otimes (z \circ 1_{\iota_A}) = (1_{\iota_A} \otimes z)
\circ (w \otimes 1_{\iota_A}) = z \circ w = (z \otimes 1_{\iota_A}) \circ
(1_{\iota_A} \otimes w) = (z \circ 1_{\iota_A}) \otimes (1_{\iota_A} \circ w)
= z \otimes w$, thus $\otimes$ and $\circ$ in this case agree and define a
commutative monoid with unit $1_{\iota_A}$.  
\hfill \halmos

\begin{rem}
As for categories, we can describe $2$-categories in terms of
$2$-arrows. One considers a set with two operations $\otimes , \circ$ and
units such that each operation gives the set the structure of a
category. Furthermore one supposes that all $\otimes$-units are also
$\circ$-units and that an associativity relation as above 
$$(S \otimes T) \circ (S' \otimes T') = (S \circ S') \otimes (T \circ T')$$ holds for the two products. 
\end{rem}

A weak 2-category (the term  bicategory is also used) is a 2-category as above, where the
associativity and unit 
identities are replaced by natural isomorphisms satisfying pentagon and triangle
axioms. As we will deal (almost) only with the strict case, will not state
these relations explicitly.

 A 2-$C^*$-category is a $2$-category for which the following hold:

 \begin{itemize}
 \item for each pair of $1$-arrows $\rho,\sigma$ the space $\Hom(\rho,\sigma)$ is a complex Banach space.
 \item there is an anti-linear involution $*$ acting on $2$-arrows, i.e.\  $*: \Hom(\rho,\sigma) \rightarrow
\Hom(\sigma,\rho),\ S\mapsto S^*$.
 \item the Banach norm is sub-multiplicative (i.e.\  $\| T \circ S \| \leq \| T
 \| \| S \|,$  when the composition is defined) and satisfies the $C^*$- condition $\| S^* \circ S \| =
 \|S\|^2$.
 \item for any $2$-arrow $S \in \Hom(\rho,\sigma)$, $S^* \circ S$ is a positive
 element in $\End(\rho)$.  
 
 \end{itemize}

\begin{rem}
  The above axioms imply that for each $1$-arrow $\rho$ the space
  $\End(\rho)$ is a unital $C^*$-algebra. For each unit $1$-arrow $\iota_A$ the
  space $\End(\iota_A)$ is a commutative unital $C^*$-algebra.
\end{rem}

We assume our category to be closed under projections (or retractions). By this
we mean the following: take a 1-arrow $B \xleftarrow {\rho} A$ and consider the
space $\End(\rho)$, which has also the structure of an algebra, as we have
seen. Then for each projection $P \in \End(\rho)$ there exists a
corresponding sub-1-arrow $ A \xleftarrow {\rho_{_{P}}} B$ and an isometry $W \in
\Hom(\rho_{_{P}} , \rho)$ such that $W^* \circ W = 1_{\rho_{_{P}}}$ and $W  \circ W^* = P$. 

We assume our category is closed under direct sums. By this we mean that for
each pair of 1-arrows $ B \xleftarrow {\rho_{1},\rho_{2}} A$ there exists a
1-arrow $B \xleftarrow {\rho} A$ and
isometries $W_1 \in \Hom(\rho_{1},\rho)$ , $W_2 \in \Hom(\rho_{2},\rho)$ such that
$W_1 \circ W_1^* + W_2 \circ W^*_2 = 1_{\rho}$ and $W_{i}^* \circ W_{j} =
1_{\rho_{i}} \delta_{i,j}$. 
Consistently with the previous definition, $\rho_1,\rho_2$ are sub 1-arrows
of
$\rho$. 
We will sometimes simply write $\rho_1 \oplus \rho_2$ for a direct
sum. Analogously we may identify a projection $P \in \End(\rho)$ with the
unit $1_{\rho_{_{P}}}$ of its
corresponding sub-1-arrow.
Thus if  $W_1, W_2, \rho_1, \rho_2 , \rho$ are as above and $T_1 \in
\End(\rho_1)$, $ T_2 \in\End(\rho_2)$,
we will simply indicate by $T_1 \oplus T_2 \in\End(\rho_1 \oplus \rho_2)$ the
2-arrow 
$W_1 \circ T_1 \circ W_1^* + W_2 \circ T_2 \circ  W_2^*$.

We assume that the category is closed under conjugation, that is, for each 1-arrow $\rho$
going from $A$ to $B$ there exists another 1-arrow$\bar \rho $ from $B$ to $A$ and two
2-arrows $R_\rho \in \Hom(\iota_{A},\bar \rho \rho)$ and $\bar R_\rho \in \Hom(\iota_{B}, \rho
\bar \rho)$  satisfying the following relations: 
$$\bar R_\rho^*\otimes 1_\rho \circ 1_\rho\otimes R_\rho = 1_\rho ; \quad
R_\rho^*\otimes 1_{\bar \rho} \,\circ \, 1_{\bar\rho}\,\otimes \bar
R_\rho = 1_{\bar \rho}.$$ 
This property is symmetric, i.e.\  if $\bar \rho$ is a conjugate for $\rho$,
then $\rho$ is a
conjugate for $\bar \rho$, as is easily seen by taking   $R_{\bar \rho}: = \bar
R_{\rho}, \ \bar R_{\bar \rho}:= R_{\rho}$ as solutions.
 $R_\rho$ and $\bar{R}_\rho$ are  fixed up to a choice of an invertible element in
$\End(\rho)$, i.e.\  if $R'_\rho$ and $\bar{R}'_\rho$ is another solution, then there exists an
invertible $A \in \End(\rho)$ such that $R'_\rho = (1_{\bar{\rho}} \otimes A) \circ R_\rho$
and $\bar{R}'_\rho = (A^{-1 *} \otimes 1_{\rho}) \circ \bar{R}_\rho$. In fact, simply take
$A = (\bar R^*_\rho \otimes 1_{\rho}) \circ (1_{\rho} \otimes R'_\rho)$. The
same holds for $\End(\bar \rho).$

Conjugacy is determined up to isomorphism, i.e.\  given
conjugate a 1-arrows $\rho$ and $\bar  \rho$ with solution $R_\rho,\bar R_\rho$, any
other  1-arrow $\bar \rho'$ conjugate to $\rho$ is isomorphic to $\bar \rho$.
In fact, let $R'_\rho,\bar R'_\rho$ be solutions for $\rho$ and $\bar \rho'$, then
$(1_{\bar \rho} \otimes \bar R^{' *}_\rho ) \circ ( R_\rho \otimes 1_{\bar \rho'}) \ \in
\Hom(\bar \rho',\bar \rho)$  is invertible.

Given two pairs of conjugate $1$-arrows $\rho_1, \bar \rho_1$ and $\rho_2,
\bar \rho_2$ with solutions $R_1, \bar R_1$ and $R_2, \bar R_2$ respectively,
one can check  that  their  sum $R_1 \oplus R_2, \bar R_1 \oplus \bar R_2$  is 
a solution for  the couple of conjugate 1-arrows $\rho_1 \oplus \rho_2$ and $\bar \rho_2 \oplus \bar \rho_2$.
Analogously given two pairs $\rho, \bar \rho$ and $\sigma, \bar \sigma$
with solutions $R_{\rho}, \bar R_{\rho}$ and $R_{\sigma},\bar R_{\sigma}$
respectively, such that the composition $\sigma \otimes \rho$ is defined, one
can consider the product solution for $\sigma \otimes \rho, \bar \rho \otimes
\bar \sigma$ defined as $R_{\sigma \otimes \rho}:= (1_{\bar \rho} \otimes
R_{\sigma} \otimes 1_{\rho}) \circ R_{\rho}, \ \bar R_{\sigma \otimes \rho}:=
(1_{\sigma} \otimes \bar R_{\rho} \otimes 1_{\bar \sigma}) \circ \bar
R_{\sigma}$.

The conjugate  relations imply, among other things, Frobenius duality, i.e.\ the following isomorphisms:
$\Hom(\rho, \sigma \otimes \eta) \cong\Hom (\rho \otimes \bar{\eta}, \sigma) \cong
\Hom(\bar{\sigma} \otimes \rho, \eta) \cong \Hom(\bar{\eta} \otimes \bar{\sigma} \otimes \rho,
\iota_{A})$.

We recall the definition of the $\bullet$ map introduced in \cite{LR97}.

\begin{defn}
Given two  1-arrows $\rho, \sigma$ , their conjugates $\bar{\rho}, \bar{\sigma}$, and  a
choice of solutions to the conjugation equations 
$ R_{\rho}$,$\bar{R}_{\rho}$, $R_{\sigma}$,$\bar{R}_{\sigma}$, 
 we define  the map $\bullet: \Hom(\rho,\sigma) \rightarrow \Hom(\bar{\rho},\bar{\sigma})$  by 
$$  \  S^{\bullet}:=  (1_{\bar{\sigma}} \otimes \bar{R}_{\rho}^*) \circ
(1_{\bar{\sigma}} \otimes S^* \otimes 1_{\bar{\rho}}) \circ (R_{\sigma}
\otimes 1_{\bar{\rho}}), \ \forall S \in \Hom(\rho,\sigma),
$$ 

and $\bullet: \Hom(\bar \rho, \bar \sigma) \rightarrow \Hom(\rho,\sigma)$ by 
$$ 
T^{\bullet}:= (1_{\sigma} \otimes \bar R^*_{\bar
  \rho}) \circ (1_{\sigma} \otimes  T^* \otimes 1_{\rho}) \circ ( R_{\bar
    \sigma} \otimes 1_{\rho}), \ \forall
T \in\Hom(\bar \rho, \bar \sigma).  
$$
\end{defn}

It is an  anti-linear isomorphism, and its square is the identity.
When $\rho = \sigma$ it is an algebraic anti-linear isomorphism and it satisfies $1_{\rho}^{\bullet} = 1_{\bar \rho}$.
Notice that $\bullet$ depends on the choice of the solution $R, \bar R$, and
that in general it does not commute with the $*$ operation.

$R^* \circ R$ and $\bar R^* \circ \bar R$ are positive elements of the commutative $C^*$-
algebras $\End(\iota_A)$ and $\End(\iota_B)$, respectively,  
so they can be thought of as positive functions in $C(\Omega_A)$ and $C(\Omega_B)$,
$\Omega_A$ and $\Omega_B$ the spectra of the two commutative algebras. 

But we can say more, the functions $R^* \circ R$ and $\bar R^* \circ \bar R$
are strictly positive on their supports.
Thus each 1-arrow $\rho$ defines a projection in $C(\Omega_A)$, namely the
projection on the (clopen) support of  $R^* \circ R$, and
analogously for the support of  $\bar R^* \circ \bar R$ in $C(\Omega_B)$. These
projections do not depend on the choice of the solutions of the conjugate equations $R$
and $\bar R$. 
In particular, if $\Omega_A$ and $\Omega_B$ are  connected, then $R^* \circ R$ and
$\bar{R}^* \circ \bar{R}$ are  positive invertible functions. 

These assertions are consequences of  the following lemmas and propositions,
most of which have been taken from \cite{Rob} and \cite{LR97}.


\begin{lem}
\label{circleunit}
Let $w \in\End(\iota_{A})$,  then the following  conditions are  equivalent:

$ a) \ 1_{\rho} \otimes w = 0,$

$b) \ R^{*} \circ R \circ w = 0.$

\noindent Similarly if $z \in\End(\iota_{B})$ the following conditions are equivalent:

$ a') \ z \otimes 1_{\rho}  = 0,$

$b') \ z \circ \bar{R}^{*} \circ \bar{R} = 0. $

\end{lem} 

{\noindent{\it Proof.}} 
The implication $a) \Rightarrow b)$ is obvious. Suppose, without loss of generality, that
$w$  is positive. Then $R^{*} \circ R \circ w = R^{*} \otimes w^{\frac{1}{2}} \circ R
\otimes w^{\frac{1}{2}} = 0$, which implies $R \otimes w^{\frac{1}{2}} = 0$ by the
$C^*$-property of the norm. 
But then $1_{\rho} \otimes w^{\frac{1}{2}} = \bar{R}^{*} \otimes 1_{\rho} \circ 1_{\rho}
\otimes R \otimes w^{\frac{1}{2}} = 0$, thus $1_{\rho} \otimes w = 1_{\rho} \otimes
w^{\frac{1}{2}}\otimes w^{\frac{1}{2}} = 0$. 
The proof of $a') \Leftrightarrow b')$ is analogous.
\hfill \halmos

The maps $ \End(\iota_{A}) \ni w  \mapsto 1_{\rho} \otimes w \in
Z(\End(\rho))$ (the centre of $\End(\rho)$) and $
\End(\iota_{B}) \ni z \mapsto z \otimes 1_{\rho} \in Z(\End(\rho))$
are $C^{*}$homomorphisms into $Z(\End(\rho)).$ Denote by $S_{l}(\rho)$ and
$S_{r}(\rho)$ the closed subspaces of $\Omega_{A}$ and $\Omega_{B}$ corresponding to the
kernels of these maps. 

\begin{lem}
\label{support}
$S_{l}(\rho)$ and $S_{r}(\rho)$ are the supports of $R^{*} \circ R$ and $\bar
{R}^{*} \circ \bar {R}$ respectively.
\end{lem}
{\noindent{\it Proof.}} 
Suppose $u \subset \Omega_{A}$ is an open subset such that  $(R^{*} \circ R)_{|_{u}} = 0$,
then for any $w \in (\iota_{A},\iota_{A})$ with support in $u$ we have $R^{*} \circ R
\circ w = 0 $ . But this implies $1_{\rho} \otimes w = 0$ by Lemma \ref{circleunit}, thus
$u \cap S_{l}(\rho) = \emptyset$ and $S_{l}(\rho) \subset \mathrm{supp}(R^{*} \circ R) $. 

If $\omega \notin S_{l}(\rho)$ then since $S_{l}(\rho)$ is closed we can find a $w \in
\End(\iota_A)$ such that $\omega(w) \neq 0$ and $\omega'(w) = 0 \ \forall \omega'
\in S_{l}(\rho)$. Thus $1_{\rho} \otimes w = 0 $, so by Lemma \ref{circleunit} $0 =
\omega(R^{*} \circ R \circ w) = \omega (R^{*} \circ R ) \omega(w)$ which implies $\omega
(R^{*} \circ R) = 0$ and $\mathrm{supp}(R^{*} \circ R) \subset S_{l}(\rho)$. 
The proof for $S_{r}(\rho)$ is analogous.
\hfill \halmos

\begin{cor}
The supports of $R^{*} \circ R$ and $\bar {R}^{*} \circ \bar {R}$ do not depend on the
choice of $R$ and $\bar{R}$. 
\end{cor} 

\begin{cor}
\label{supportcutting}  Let $E_{S_l(\rho)}$ and $E_{S_r(\rho)}$ denote the projections
onto the supports 
of $R^* \circ R$ and $\bar R^* \circ \bar R$ respectively. Then the following
equalities hold: $1_\rho = 1_\rho \otimes E_{S_l(\rho)} = E_{S_r(\rho)}
\otimes 1_\rho.$
\end{cor}

\begin{lem}
\label{cupunit}
The following inequalities hold:

$R \circ R^{*} \leq (R^{*} \circ R) \otimes 1_{\bar{\rho} \rho}$ 

$R \circ R^{*} \leq 1_{\bar{\rho} \rho}  \otimes (R^{*} \circ R)$

$\bar{R} \circ \bar{R}^{*} \leq (\bar{R}^{*} \circ \bar{R}) \otimes 1_{\rho \bar{\rho}}  $ 

$\bar{R} \circ \bar{R}^{*} \leq 1_{\rho \bar{\rho}}   \otimes (\bar{R}^{*} \circ \bar{R})$.
\end{lem}

{\noindent{\it Proof.}} 
Notice that $(R \circ R^{*}) \circ (R \circ R^{*}) = (R \circ R^{*}) \circ 1_{\bar{\rho}
  \rho} \otimes (R^{*} \circ R) = (R \circ R^{*}) \circ (R^{*} \circ R) \otimes
1_{\bar{\rho} \rho},$ 
where we are regarding  $(R \circ R^{*})$, $1_{\bar{\rho} \rho} \otimes (R^{*} \circ R)$,
$(R^{*} \circ R) \otimes 1_{\bar{\rho} \rho}$ as positive elements of the algebra
$\End(\rho)$. 
In particular $1_{\bar{\rho} \rho} \otimes (R^{*} \circ R)$ and $(R^{*} \circ
R ) \otimes 1_{\bar{\rho} \rho}$ are elements of the centre $Z(\End(\bar \rho \rho))$.
Analogous relations hold for $ \bar{R} \circ \bar{R}^{*}$.
Now, in general, if we have a positive element $X$ in a $C^*$-algebra $A$ such that $X^2 = X Z$, where Z is a positive element of the centre of $A$, we have $Z \geq X$.
In fact, take a faithful representation $(\pi,H)$ of the algebra $A$, take two
generic vectors $\alpha  \in  H, \beta \in (X^{\frac{1}{2}}H)^{\perp}.$ Then 
$(X^{\frac{1}{2}} \alpha + \beta, X (X^{\frac{1}{2}} \alpha + \beta)) \  = \  (\alpha, X^{2} \alpha)$   
 and $(X^{\frac{1}{2}} \alpha + \beta, Z (X^{\frac{1}{2}} \alpha + \beta)) \ = \ (\alpha, X^{2} \alpha) + (\beta, Z \beta).$
Thus $Z \geq X$.
\hfill \halmos

\begin{prop}
\label{pimsner-popa}
For each positive $X \in \End(\rho \otimes \sigma)$ the following inequality holds:

$X \leq ( \bar{R}^* \circ \bar{R}) \otimes 1_{\rho} \otimes (R^* \otimes 1_{\sigma}\circ  (1_{\bar{\rho}} \otimes X)  \circ R \otimes 1_{\sigma})$.

\end{prop}

{\noindent{\it Proof.}}  $$ X = (1_{\rho} \otimes R^{*} \otimes 1_{\sigma}) \circ (\bar{R} \circ
\bar{R}^{*}) \otimes X \circ (1_{\rho} \otimes R \otimes 1_{\sigma}) \leq $$

$$ 
(1_{\rho} \otimes R^{*} \otimes 1_{\sigma} ) \circ (\bar{R}^{*} \circ
\bar{R})\otimes 1_{\rho \bar \rho} \otimes X \circ (1_{\rho} \otimes  R
\otimes 1_{\sigma})$$

$$ = (\bar{R}^{*} \circ \bar{R}) \otimes 1_{\rho}  \otimes (R^{*} \otimes 1_{\sigma} \circ
(1_{\bar{\rho}} \otimes X) \circ   R \otimes 1_{\sigma}),$$
where  in the first line we have used the conjugation equations and in the
 second we have used the third inequality of the preceding lemma.
\hfill \halmos

\begin{cor}
\label{serpentwegg}
The following inequality holds:

$(\bar{R}^{*} \circ \bar{R}) \otimes 1_{\rho} \circ 1_{\rho} \otimes (R^{*} \circ R) \geq 1_{\rho}$.
\end{cor}

\begin{cor}
\label{closets}
The following hold:

$i)$ $1_{\rho} \otimes R^{*} \circ R \geq \frac{1}{||\bar{R}||^{2}} 1_{\rho}$
; $\bar{R}^{*} \circ \bar{R} \otimes 1_{\rho} \geq \frac{1}{||R||^{2}} 1_{\rho}$

$ii)$ $(R^* \circ R)_{|_{S_{l}(\rho)}} \geq \frac{1}{||\bar{R}||^{2}} $ ; $(\bar{R}^* \circ \bar{R})_{|_{S_{r}(\rho)}} \geq \frac{1}{||R||^{2}}$

$iii)$ $S_{l}(\rho)$ and $S_{r}(\rho)$ are open and closed.
\end{cor}

\begin{lem}
\label{swish}
Let $\rho$,$\sigma$ be 1-arrows from  $A$ to $B$ and 
$E_{S_{l}(\rho)}$ and $E_{S_{l}(\sigma)}$ the associated   projections on
$S_{l}(\rho)$ , $S_{l}(\sigma)$ in $\Hom(\iota_A,\iota_A)$. Suppose $E_{S_{l}(\rho)} E_{S_{l}(\sigma)} = 0$, then $\Hom(\rho,\sigma)=0$.  
An analogous assertion holds for the right supports $S_{r}(\rho)$,$S_{r}(\sigma)$.
\end{lem}
{\noindent{\it Proof.}} 
If  $T \in \Hom(\rho,\sigma)$, then $T = 1_{\sigma} \circ T \circ 1_{\rho} = (1_{\sigma} \otimes E_{S_{l}(\sigma)}) \circ T \circ (1_{\rho} \otimes E_{S_{l}(\rho)}) = 1_{\sigma} \circ (T \otimes (E_{S_{l}(\sigma)} E_{S_{l}(\rho)}) \circ 1_{\rho} = 0$. \hfill \halmos

We will meet a refinement of this lemma in section \ref{bundles}.

We introduce now a construction which will be useful in the sequel. 
Suppose we have chosen for each object $A$ of our 2-$C^*$-category $\mathcal{A}$ a set  of projections $ \{ 1_{\iota_{A_i}} 
\} $  in the
associated commutative $C^*$-algebra 
$\End(\iota_A)$  
 such that  $\sum_i 1_{\iota_{A_i}} = 1_{\iota_A}$ 
(i.e.\  the set is complete) and $1_{\iota_{A_i}} \otimes 1_{\iota_{A_j}} = \delta_{i,j} 1_{\iota_{A_i}}$ (i.e.\  the
projections are orthogonal).
To each projection $1_{\iota_{A_i}}$ there will correspond a 1-arrow, which we will call
$\iota_{A_i}$. 
We would like to  think of these 1-arrows as units corresponding to
objects.
In other words, we would like to ``decompose'', in some sense,  each object into sub-objects
corresponding  to our original choice of sets of projections.
1-arrows and 2-arrows should be decomposed accordingly. We must show that
this can be done in a consistent manner.
We define a new 2-$C^*$-category $\mathcal{B}$  the following way (we will use the ``only $2$-arrows approach'' mentioned above):
\begin{itemize}

\item  define as  $2$-arrows of the new category $\mathcal{B}$ all  the elements of the  form $$\{ 1_{\iota_{B_j}} \otimes S
  \otimes 1_{\iota_{A_i}}, \ \ \forall \  1_{\iota_{B_j}}, 1_{\iota_{A_i}}, \ B \xleftarrow{\rho} A, \ B
  \xleftarrow{\sigma} A, \ S \in \Hom(\rho,\sigma), \ A,B \in \mathcal{A} \}$$ with the same $\otimes$ and $\circ$ operations of the
  original category.
\item we set each projection $1_{\iota_{A_i}}$ to be a $\otimes$-unit (thus, also a
  $\circ$-unit),
\item we define the  $\circ$-units to be the set $ \{ 1_{\iota_{B_j}} \otimes
  1_{\rho} \otimes 1_{\iota_{A_i}} , \ \ \forall \ \ 1_{\iota_{A_i}}, \ 1_{\iota_{B_j}}, \ A \xleftarrow{\rho} B\}$

\end{itemize}

The units satisfy the necessary properties by definition and   compatibility between the
$\otimes$ and the $\circ$ products descends from 
the original one. The new $2$-category is still closed under conjugation. In
fact, let $B \xleftarrow{\rho} A$ and $ A \xleftarrow{\bar \rho} B$ be two
conjugate $1$-arrows in the original category. Then  each $ 1_{\iota_{B_j}} \otimes
1_{\rho} \otimes 1_{\iota_{A_i}}$ has as conjugate $1_{\iota_{A_i}} \otimes 1_{\bar \rho}
\otimes 1_{\iota_{B_j}}$ with   conjugate solutions
$ R_{i,j}:= (1_{\bar \rho} \otimes 1_{\iota_{B_j}} \otimes
1_{\rho}) \circ R\circ 1_{\iota_{A_i}}$
and  $ \bar R_{i,j}:= (1_{ \rho} \otimes 1_{\iota_{A_i}} \otimes 1_{\bar \rho}) \circ
\bar R \circ 1_{\iota_{B_j}}.$ 


\begin{rem}
The above construction is not an inclusion of the 2-$C^*$-category
$\mathcal{A}$ in $\mathcal{B}$ by in the proper
sense. In fact, each unit $1_{\iota_A}$ is not sent to one  but to a  set of units  $1_{\iota_{A_i}}$. 
The original  2-$C^*$-category $\mathcal{A}$ is easily recovered from
the new one by considering
linear combinations of the elements of $\mathcal{B}$ (the fact that linear
combinations of 2-arrows in $\mathcal{B}$ span all of the 2-arrows in
$\mathcal{A}$ is a consequence of Lemma \ref{swish}).       
\end{rem}

\begin{rem}
  Obviously this construction depends on the choice of the sets of projections
  $\{ 1_{\iota_{A_i}} \}$. If we suppose the topological spaces
  $\Omega_A$ corresponding to each object $A$ to be a finite union of
  connected components $ \{ \Omega_{A_i} \}$ then the  central projections
  associated to each component would be a
  natural choice. Notice that in this case in the new category $\mathcal{B}$
  the objects $A_i$ would have connected spectrum, thus for any $ B_j
  \xleftarrow{\rho} A_i$ the elements $R^{*}{_\rho} \circ R_\rho$ and ${\bar
  R}^*_\rho \circ {\bar R}_\rho$ would be invertible as a consequence of Corollary \ref{closets}.
\end{rem}

\section{Bimodules and bundles.}
\label{bundles}

In this section we will pursue a parallel path with respect to the case of simple
 units. We will show that Banach, Hilbert and $C^*$-algebra bundles appear as a
 generalisation of finite dimensional spaces. Several analogous results are
 obtained, such as a finite upper bound on the dimension of the fibres, given
 by the non scalar analogue of the dimension function.
We investigate the behaviour of these fibres with respect
 to the categorical structure: the $\circ$ composition, the $*$ involution and
 conjugation naturally preserve the fibres and are straightforward to
 describe. The behaviour under the $\otimes$ composition is less friendly and
 for this purpose we introduce a further hypothesis. Our assumption is not the
 most general a priori, but general enough to comprehend most (if not all) of the
 known examples (e.g.\  $\End(\iota)$-linear tensor categories, such as braided
 categories, are included). Furthermore it does not seem easy to produce a
 counterexample, a question we leave open for
 the future. 

 We begin by noticing that  the spaces $\Hom(\rho,\sigma)$ have  a structure of $\End(\iota_B)-\End(\iota_A)$ Hilbert bimodule
 given by the conjugation relations and the $\otimes$ product. In fact, given
 $z \in \End(\iota_B)$ and $w \in \End(\iota_A)$ and $T \in
 \Hom(\rho,\sigma),$ we can consider the tensor products $z \otimes T$ and $T
 \otimes w$ both still in $\Hom(\rho,\sigma)$. Given $T,S \in \Hom(\rho,\sigma)$ we
 define the right $\End(\iota_A)$-valued product $\langle T,S
 \rangle^{(\rho,\sigma)}_{{{\End}(\iota_A)}}$ by $R^* \circ (1_{\bar \rho} \otimes (S^* \circ T)) \circ R$. 
The same way we define the left $\End(\iota_B)$-valued product
 $_{{\End}(\iota_B)}\langle T,S \rangle ^{(\rho,\sigma)}$ by $\bar R^* \circ ((S^* \circ T) \otimes 1_{\bar \rho}) \circ \bar R$. 
Both of these products are non-degenerate.

We recall the definition of a Banach bundle (the terminology ``continuous
field of Banach spaces'' is also used in the literature, cf.\ e.g., \cite{DD63}).

\begin{defn}
Let $\Omega$ be a compact Hausdorff topological space. A Banach bundle $E$ over $\Omega$ is a family  of Banach spaces $\{ E_{\omega}  , \  \| . \|^{\omega} \ , \  \omega \in \Omega \}$ with a  set  $\Gamma \subset \prod_{\omega \in \Omega} E_{\omega}$ such that:

i ) $\Gamma$ is a linear subspace of $\prod_{\omega \in \Omega} E_{\omega}$

ii )  $\forall \omega \in \Omega$  $ \{ S_{|_{\omega}}$ , $ S \in \Gamma \}$    is dense in $E_{\omega}$

iii ) $\forall S \in \Gamma $ the norm function $\omega \rightarrow \|S_{|_{\omega}}\|^{\omega}$ is a continuous function on $\Omega$

iv)  Let $X \in \prod_{\omega \in \Omega} E_{\omega} $ . If $\forall \omega
 \in \Omega $ and $ \forall \epsilon > 0  \ \exists S \in \Gamma $ such
 that $ \|X_{|_{\omega}} - S_{|_{\omega}}\|^{\omega} < \epsilon$ in a neighbourhood of $\omega$, then $X \in \Gamma$.

\end{defn}

The elements of $\prod_{\omega \in \Omega_{A}} E_{\omega}$ are called sections and those of $\Gamma$ continuous sections.
Analogously if the spaces $E_{\omega}$ have the structure of Hilbert spaces and the norm is given by the inner product, we will talk about Hilbert bundles. If the $E_{\omega}$ have the structure of $C^*$-algebras, and the space of continuous sections $\Gamma$ is closed under multiplication and  the $*$ operation, we will talk about a $C^*$-algebra bundle.
In a Banach bundle the fibre space might vary according to the base point. So
it is a more general situation than that of a locally trivial bundle. The choice of a set $\Gamma \subset \prod_{\omega} E_{\omega}$ as the space of continuous sections is part of the initial data, as in general  we have no local charts, with the implicit notion of continuity given by them.

\begin{prop}
\label{bundlea}
Given two 1-arrows $\rho,\sigma$ from objects $A$ to $B$ and a choice of the
conjugation equations $R_{\rho}, \bar R_{\rho}$ for $\rho$, $\Hom(\rho,\sigma)$ has
the structure of a Hilbert bundle. 
\end{prop}

{\noindent{\it Proof.}}  We evaluate the product $\langle S,T \rangle_{{\End}(\iota_A)}^{(\rho,\sigma)}$ on each $\omega \in \Omega_A$ for any $S,T \in \Hom(\rho,\sigma)$. 
The  procedure of the GNS construction  gives us for each point $\omega$ a Hilbert space,
which we shall denote by $\Hom(\rho,\sigma)_{\omega}$. We take $\prod_{\omega \in
\Omega_{A}} \Hom(\rho,\sigma)_{\omega}$ as fibre bundle and the image of
$\Hom(\rho,\sigma)$ (which we will still denote  by $\Hom(\rho,\sigma)$) as the module of
continuous sections. 

Note that for $S \in \Hom(\rho,\sigma)$ the topology given by $sup_{|_{\omega \in \Omega}} \| S_{\omega} \|^{\omega}$ is equivalent to the original one.
In fact, on the one hand we have 
$$sup_{|_{\omega \in \Omega_{A}}} \|S_{|_{\omega}}\|^{\omega} = \| \langle S,S
\rangle^{\frac{1}{2}}_{\End(\iota_{A})}\| = \| \langle S,S \rangle_{\End(\iota_{A})}\|^{\frac{1}{2}}$$
  $$\leq \| S^{*} S \|^{\frac{1}{2}} \| \langle 1_{\rho},1_{\rho}\rangle
_{\End(\iota_{A})} \|^{\frac{1}{2}} = \| S^{*} S \|^{\frac{1}{2}} \| R^*_\rho
\circ R_\rho \|^{\frac{1}{2}} = \|S \| \| R_\rho \|$$ 
and on the other hand we have 
 $$\| S \| = \| S^* \circ S \|^{\frac{1}{2}} \leq \| (\bar R_\rho^* \circ \bar R_\rho )
\otimes 1_\rho \otimes  \langle S,S
\rangle_{\End(\iota_{A})}\|^{\frac{1}{2}}$$ 
$$ \leq \| (\bar R_\rho^{*} \circ \bar
R_\rho) \|^{\frac{1}{2}} 
\| \langle S,S \rangle_{\End(\iota_{A})}\|^{\frac{1}{2}}   = \|
\bar R^*_\rho \circ \bar R_\rho \|^{\frac{1}{2}} \| \langle S,S
\rangle_{{\End}(\iota_A)}S \|^{\frac{1}{2}}$$ $$ = \| \bar R_\rho \|
(sup_{|_{\omega   \in \Omega_A}} \| S_\omega \|^{\omega}) = \| \bar R_\rho \|
\| S_\omega \|,$$
where we have used the expression of the definition of the inner product $\langle \cdot,\cdot
\rangle_{\End(\iota_{A})}$, the monotonicity  of the square root function and the
inequality for $S^{*} S$ given by Proposition \ref{pimsner-popa}. 
So $\Hom(\rho,\sigma)$ is closed even as a subspace of the Banach bundle. 
The first three conditions either very easy to prove.
We prove only the last one.
Suppose we have $X \in \prod_{\omega \in \Omega_{A}} \Hom(\rho,\sigma)_{\omega}$
satisfying condition $iv)$. Then, as $\Omega_{A}$ is compact,  $\forall \epsilon > 0 $ we
can choose a finite family of elements $S^{\alpha} \in \Hom(\rho,\sigma)$ and a
corresponding finite open covering  $\{ U_{\alpha} \}$ of $\Omega_{A}$ such that $\|X -
S^{\alpha}\|^{\omega}_{|_{{U_{\alpha}}}} \leq \epsilon$. Take a partition of unity
$f_{\alpha}$ subordinate to the open covering. Then $\|X - \sum_{\alpha} f_{\alpha}
S^{\alpha}\|^{\omega} \leq \epsilon$ because of   convexity of the norm. 
Thus $X \in \ovl{\Hom(\rho,\sigma)} = \Hom(\rho,\sigma)$. \hfill \halmos

This gives the $\End(\iota_B) - \End(\iota_A)$ bimodule
 $\Hom(\rho,\sigma)$ the structure of a Hilbert bundle over the compact
 topological space $\Omega_A$. So we can think of each element $T$ as a
 continuous section $T_{|_{\omega}}$ in this Hilbert bundle. The right action of
 $\End(\iota_A)$ is given by multiplying  functions in $C(\Omega_A)$. 
 

 Having shown that each $\Hom(\rho,\sigma)$ has a Hilbert bundle structure, we
 would like to show, as already mentioned, its behaviour in the context of the
 whole category.   
In fact $\Hom(\rho,\sigma)$ is not only a bimodule: its elements can be regarded as operators between other spaces
by $\circ$ composition of 2-arrows.
For example  the elements of $\Hom(\rho,\sigma)$ can be regarded as operators from  $\Hom(\eta,\rho)$ to $\Hom(\eta,\sigma)$.
Suppose we have chosen solutions to the conjugation equations for $\rho$ and for $\eta$,
with the respective induced Hilbert bundle structures.  
We have the following
\begin{prop}
\label{bundleb}
If  $T, T' \in \Hom(\rho,\sigma) $ such that $T_{|_{\omega}} = T'_{|_{\omega}}$ ,$ \  P,
P'  \in \Hom(\eta,\rho)$ such that  $P_{|_{\omega}} = P'_{|_{\omega}}$ then $(T\circ
P)_{|_{\omega}} = (T' \circ P)_{|_{\omega}} = (T\circ P')_{|_{\omega}} = (T' \circ
P')_{|_{\omega}}$. 
\end{prop}

{\noindent{\it Proof.}}  $T_{\omega} = T'_{\omega}$ means $\langle T,S
\rangle^{(\rho,\sigma)}_{\End(\iota_{A}) \ |_{ \omega}} = \ \langle T',S
\rangle^{(\rho,\sigma)}_{\End(\iota_{A}) \ |_{\omega}} \forall S \in \Hom(\rho,\sigma)$. Analogous relations hold for $P,P' \in \Hom(\eta,\rho)$.
We have $$\langle T \circ P, Q \rangle^{(\eta,\sigma)}_{\End(\iota_{A}) \
  |_{\omega}} = \ \langle P,T^{*} \circ Q
\rangle^{(\eta,\sigma)}_{\End(\iota_{A}) \ |_{\omega}}  = \ \langle P',T^{*}
\circ Q \rangle^{(\eta,\sigma)}_{\End(\iota_{A}) \ |_{\omega}} $$ $$= \ \langle
T \circ P', Q \rangle^{(\eta,\sigma)}_{\End(\iota{_A}) \ |_{\omega}}   , \ \forall Q \in \Hom(\eta,\rho).$$
Thus $(T \circ P)_{|_{\omega}} = \ (T \circ P')_{|_{\omega}}$, and in the same way we have $(T' \circ P)_{|_{\omega}} = (T' \circ P')_{|_{\omega}}$.

Also we have $$\langle T \circ P, Q \rangle^{(\eta,\sigma)}_{\End(\iota_{A}) \
  |_{\omega}} =  \langle T , Q \circ P^{* \bullet * \bullet *}
\rangle^{(\rho,\sigma)}_{\End(\iota_{A}) \ |_{\omega}} =  \langle T' , Q \circ P^{* \bullet * \bullet *} \rangle^{(\rho,\sigma)}_{\End(\iota_{A}) \ |_{\omega}}$$
$$ = \langle T' \circ P, Q \rangle^{(\eta,\sigma)}_{\End(\iota_{A}) \
  |_{\omega}} ,  \ \forall Q \in \Hom(\eta,\rho),$$ thus $ (T\circ P)_{|_{\omega}} = (T' \circ P)_{|_{\omega}}$.
\hfill \halmos

The preceding proposition implies  the following
\begin{cor}
  For each $\omega \in \Omega_A$ the set $$\mathit{I}_{\omega}:= \{ S \in \End(\rho) \ s.t.  \langle S,S\rangle_{{\End}(\iota_A) |_{\omega}} =
  0 \}$$ is a closed  two-sided ideal of $\End(\rho)$.
\end{cor}

\begin{cor}
  Choose solutions $R_\rho,\bar R_\rho$ and $R_\sigma, \bar R_\sigma$  for 
  $\rho$ and $\sigma$, respectively. For  $T, T' \in \Hom(\rho,\sigma) $ if
  $T_{|_{\omega}} = T'_{|_{\omega}}$ then $T_{|_{\omega}}^* = {T'}_{|_{\omega}}^*$.  
\end{cor}

{\noindent{\it Proof.}}  By the previous proposition if $T_{|_\omega} = T'_{|_\omega}$ (i.e.\  $(T -
T')_{|_{\omega}} = 0)$) then $((T - T') \circ (T - T')^*)_{|_\omega} = 0 $,
which implies $ \langle  ((T - T') \circ (T - T')^*, 1_{\sigma}
\rangle^{(\sigma,\rho)}_{{\End}(\iota_A) \ |_\omega} = \langle (T - T')^*,(T - T')^*
\rangle^{(\sigma,\rho)}_{{\End}(\iota_A) \ |_\omega} = 0$,
i.e.\  $T^*_{|_\omega} = {T'}^*_{|_{\omega}}$. \hfill \halmos

Thus the $\circ$ and  $*$ operations preserve the fibre structure.

\begin{rem}
  $\End(\rho) / \mathit{I}_{\omega}$ is the pre-Hilbert space that gives rise to the
  fibre-Hilbert space $\End(\rho)_{\omega}$ when completed with respect to the
  pre-scalar product norm.  For each $\omega$ we can pursue the whole GNS
  construction and obtain a $C^*$-algebra $\pi_{\omega}(\End(\rho))$ acting on
  this Hilbert space. The preceding corollary shows that
  $\pi_{\omega}(\End(\rho))$ is the completion of the same
  pre-Hilbert space  $\End(\rho) / \mathit{I}_{\omega}$ with
  respect to the $C^*$-norm given by the GNS construction. 
\end{rem}

We have the following

\begin{prop}
\label{cbundle}
$\End(\rho)$ has the structure of a $C^*$-algebra bundle
\end{prop}

{\noindent{\it Proof.}} 
Proceed as in the beginning of Proposition \ref{bundlea} and for each $\omega$
consider the GNS construction. 
We must show the continuity of  this $C^*$-norm with respect to
the base point $\omega$.
For $A \in \End(\rho) / \mathit{I}_{\omega}$ we define
$$\|A\|^{C^*_1}_{\omega}:= \sup _{\tilde y \in \End(\rho) /
  \mathit{I}_{\omega}  , \   \| \tilde y \|^{Hilbert} \ \leq 1} \ \ \|A \tilde y\|^{Hilbert}_{\omega},$$ where by
$\|A \tilde y \|_{\omega}^{Hilbert}$ we mean $ \langle Ay,Ay\rangle^{\frac{1}{2}}_{\omega}$,
for any $y \in \End(\rho)$ such that $y_{|_{\omega}} = \tilde y$.
As this norm is defined as a $\sup $ over continuous functions, it is a
priori only lower semicontinuous.
As $\mathit{I}_{\omega}$ is a two-sided $C^*$-ideal, there is another candidate
$C^*$-norm, namely $\|A\|^{C^*_2}_{\omega}:= \inf_{y \in
  \mathit{I}_{\omega} }\|A - y \|$, the $C^*$-norm of the quotient
$C^*$-algebra $\End(\rho) / \mathit{I}_{\omega}$.
We show that these two norms are the same.

\begin{lem}
$\| . \|_{\omega}^{C_2^*} = \| . \|_{\omega}^{C_1^*}.$  
\end{lem}

{\noindent{\it Proof.}} 
Take an approximate unit $u_{\lambda}$ for $\mathit{I}_{\omega}$.
Then $$ \inf_{y \in \mathit{I}_{\omega}} \| A - y \| =
\lim _{\lambda \rightarrow \infty} \|A(1-u_{\lambda})\| =
$$ 
$$\lim _{\lambda \rightarrow \infty} \ \sup _{\phi \in
  \mathcal{S}\End(\rho)} \ \phi ( (A(1-u_{\lambda}))^*
A(1-u_{\lambda}))^{\frac{1}{2}} = $$
$$\lim _{\lambda \rightarrow \infty} \ \sup _{\phi \in
  \mathcal{P}\End(\rho)} \ \phi ( (A(1-u_{\lambda}))^*
A(1-u_{\lambda}))^{\frac{1}{2}},$$
where  $\mathcal{S}\End(\rho)$ and $\mathcal{P}\End(\rho)$ are the states
and the pure states, respectively,  of $\End(\rho)$.

Notice that as $ \langle \cdot ,\cdot \rangle$ is a non-degenerate $C(\Omega_A)$-valued inner
product, we can restrict ourselves evaluating  the supremum on pure states of the algebra $\End(\rho)$  dominated by 
states of the form $ \langle y,y\rangle_{|_{\omega'}}$ for some $\omega' \in
\Omega_A$ and some $y \in \End(\rho)$ such that $ \langle y,y\rangle_{|_{\omega'}} = 1$.

For each $u_{\lambda}$ choose a sequence $\phi^{\lambda}_n \in \mathcal{P}\End(\rho)$ such that
$\lim _{n \rightarrow \infty}
\ \phi^{\lambda}_n((A(1-u_{\lambda}))^*A(1-u_{\lambda}))^{\frac{1}{2}} =
\|A(1-u_{\lambda}\|$. Then choose a diagonal sequence $\phi_{\lambda}$ such
that$\lim _{\lambda \rightarrow \infty}
\ \phi_{\lambda}((A(1-u_{\lambda}))^*A(1-u_{\lambda}))^{\frac{1}{2}} =
\lim _{\lambda \rightarrow \infty} \|A(1-u_{\lambda}\| .$
Let $\phi_{0}$ be an accumulation point of this last sequence. Then
$\lim _{\lambda \rightarrow \infty}
\ \phi_{0}((A(1-u_{\lambda}))^*A(1-u_{\lambda}))^{\frac{1}{2}} =
\lim _{\lambda \rightarrow \infty} \ \|A(1-u_{\lambda})\|.$ Suppose that $\phi_{0}$ is dominated by a state of the
kind $ \langle y,.y\rangle_{\omega'}$ for some $\omega' \neq \omega$.
Then choosing a continuous function $g \in C(\Omega_A)$ such that $g(\omega')
= 1$ and $g(\omega) = 0$ we have $\lim _{\lambda \rightarrow \infty}
\phi_{0}((A(1-u_{\lambda}))^*A(1-u_{\lambda})) = \lim _{\lambda
  \rightarrow \infty }
\phi_{0}(g(A(1-u_{\lambda}))^*A(1-u_{\lambda}))^{\frac{1}{2}}$. But the last
term is zero, as $g \in \mathit{I}_{\omega}$.
So any accumulation point must be dominated by a state of the form
$ \langle y,.y\rangle_{|_{\omega}}$. But in $\omega$ each $u_{\lambda} = 0$, so
$\phi_{0}((A(1-u_{\lambda}))^*A(1-u_{\lambda})) =\phi_{0}(A^*A)$.

So we conclude $$\| A \|_{\omega}^{C_2^*} = \inf _{y \in \mathit{I}_{\omega}} \| A - y\| =
\phi_{0}(A^*A)^{\frac{1}{2}} = $$
$$\sup _{\tilde y \in \End(\rho) / 
  {\mathit{I}}_{\omega} , \ {\| \tilde y \|}^{Hilbert}_{\omega} = 1}   \langle A
\tilde y, A \tilde y\rangle^{\frac{1}{2}}_{|_{\omega}} = \| A \|_{\omega}^{C_1^*}.$$

Now we must show that $\|.\|^{C_2^*}$ is upper semicontinuous. Suppose $y^{\omega}
\in \mathit{I}_{\omega}$ such that $\| A - y^{\omega} \| = \|A
\|_{\omega}^{C_2^*} + \epsilon$. Then choose a neighbourhood $U_{\omega}
\subset \Omega_A$ such that $\| y^{\omega} \|_{\omega' }^{Hilbert} \leq
\epsilon \ \forall \omega' \in U_{\omega}.$

 Then $\| A \|_{\omega'}^{C_2^*} = \inf _{y^{\omega'} \in
 \mathit{I}_{\omega'} } \| A - y^{\omega'} \| \leq \inf _{y^{\omega'}
 \in \mathit{I}_{\omega'}} (\|A - y^{\omega} \| + \| y^{\omega} - y^{\omega'} \|).$

But $\| y^{\omega} -y^{\omega'}\| \leq \mathit{const} \|  \langle (y^{\omega}
- y^{\omega'}), (y^{\omega} -y^{\omega'})\rangle \|^{\frac{1}{2}}$, where the constant
is independent of $y^{\omega}, y^{\omega'}$ and  is given by Lemma
\ref{pimsner-popa}, and $\inf _{y^{\omega'} \in \mathit{I}_{\omega'}}
\|  \langle (y^{\omega} - y^{\omega'}),(y^{\omega} - y^{\omega'})\rangle \|^{\frac{1}{2}}
\leq \epsilon.$
Thus $\| A \|_{\omega'}^{C_2^*} \leq \|A \|_{\omega}^{C_2^*} + \epsilon +
\mathit{const} \epsilon$ for any $\omega' \in U_{\omega}$, which means that
$\|.\|_{\omega}^{C_2^*}$ is upper semicontinuous.

Having proven continuity of the $C^*$-norm, the rest of the proof follows as in
Proposition \ref{bundlea}. $\hfill \halmos$

\begin{rem}

We have introduced two kinds of bundle structures for $\End(\rho)$, the first one giving each
fibre a scalar product norm, the second a $C^*$-norm. 
As the inner $\End(\iota_A)$ and $\End(\iota_B)$ products depend on
the choice of solutions to the conjugation equation, the Hilbert bundle structure is defined
only up to an isomorphism. In fact, consider $\Hom(\rho,\sigma)$ with the bundle structure given by a
solution $R,\bar R$ for $\rho$ and a second bundle structure given by a second solution
$(1_{\bar \rho} \otimes A) \circ R , \ (A^{-1 *} \otimes 1_{\bar \rho}) \circ
\bar R$, where $A \in \End(\rho)$ is an invertible.
Then the map $ S \mapsto S \circ A^{-1} \ \forall S \in \Hom(\rho,\sigma)$ is a unitary
map between the two Hilbert bundles structures.

Even for $\rho \neq \sigma$ we can define for $\Hom(\rho,\sigma)$ a Banach bundle structure where the fibre norm satisfies a $C^*$-condition by defining $\| S_{|_{\omega}} \|_{\omega}:= (\| (S^* \circ S)_{|_{\omega}}\|_{\omega}^{C^*})^{\frac{1}{2}}$. 
We will consider this structure by default, if not specified otherwise.
As we shall see in the sequel, the fibres are finite dimensional, so the two
types of norms are equivalent. The latter bundle structure endows each Banach space fibre with its unique
$C^*$-norm.
\end{rem}
We also have the following
\begin{prop}
\label{finitefibres}
For each $\End(\rho)$ and for each point $\omega$ of the base space $\Omega$ the associated $C^*$-algebra fibre space $\End(\rho)_{\omega}$ is finite dimensional.  
\end{prop}

{\noindent{\it Proof.}}  This is essentially the same proof as in \cite{LR97}. 
Consider a set $\{ X_{i |_{\omega}} \in \End(\rho)_{\omega} \}$ of positive elements of norm one with $\sum_{i} X_{i |_{\omega}} \leq 1_{\rho |_{\omega}} $. 
We can think, without loss of generality, of each $ X_{i |_\omega}$ as the
value in $\omega$ of a positive $ X_i \in \End(\rho)$ (if not, just take
$(X^*_i \circ X_i)^{\frac{1}{2}}$ instead).  
By the inequality in Corollary \ref{pi-po.om} we have 
$$    X_i   \leq ((\bar {R}^{*} \circ \bar {R})  \otimes 1_{\rho }
\otimes (R^{*}  \circ 1_{\bar {\rho} } \otimes  X_{i } \circ R)) 
= (\bar R^* \circ \bar R) \otimes 1_\rho \otimes \langle X_i , 1_\rho \rangle
_{{\End}(\iota_A)}.$$
Notice that $ ((\bar R^* \circ \bar R) \otimes 1_\rho \otimes \langle X_i , 1_\rho \rangle
_{{\End}(\iota_A)})_{|_\omega}$ is simply $((\bar R^* \otimes \bar R)
\otimes 1_\rho)_{|_{\omega}}$ times the positive constant  $(\langle X_i , 1_\rho \rangle
_{{\End}(\iota_A)}) _{ |_{\omega}}.$
As the norm of  each single $X_{i |_{\omega}}$ is one, we have
  $$ 1 = \| X_{i |_{\omega}} \|^{\omega} \leq \| ((\bar {R}^{*} \circ \bar {R})  \otimes
  1_\rho )_{|_\omega} \|^{\omega} \ ( \langle X_i , 1_\rho \rangle
  _{{\End}(\iota_A)} ) _{|_{\omega}}.$$
Summing over $i$  we have 
 $$ n = \sum_i \| X_{i |_{\omega}} \|^{\omega} \leq \sum_i \| ((\bar {R}^{*} \circ \bar {R})  \otimes
  1_\rho )_{|_\omega} \|^{\omega} \ ( \langle X_i , 1_\rho \rangle
  _{{\End}(\iota_A)} ) \ _{|_{\omega}} $$
$$\leq \| ((\bar R^* \circ R) \otimes
  1_\rho)_{|_\omega} \|^{\omega}  (\langle 1_\rho , 1_\rho \rangle
  _{{\End}(\iota_A)})_{|_\omega} =  \| ((\bar {R}^{*} \circ \bar {R}) \otimes
  1_\rho)_{|_\omega} \|^{\omega}   (R^{*} \circ R)_{|_{\omega}},$$ which is 
  finite.\hfill \halmos


\begin{rem}
A posteriori we see that the fibre spaces $\End(\rho) / \mathit{I}_{\omega}$
are finite dimensional, thus the completion  in the Hilbert and
the $C^*$-norm was superfluous. They are finite dimensional  $C^*$-algebras.

\end{rem}

\begin{cor}
For each $\Hom(\rho,\sigma)$ and for each point of the base space $\Omega$ the
associated Hilbert fibre space $\Hom(\rho,\sigma)_{\omega}$ is finite dimensional.
\end{cor}

{\noindent{\it Proof.}}   That 
$\Hom(\sigma,\rho)_{\omega} \circ \Hom(\rho,\sigma)_{\omega}$ is finite dimensional  follows  from the preceding proposition and remark 
by
embedding  it into $\End(\rho)_{\omega}$. But the map $\Hom(\rho,\sigma)_\omega
\rightarrow \Hom(\rho,\sigma)_\omega \circ \Hom(\sigma,\rho)_\omega: \ S_\omega \mapsto
S^*_\omega \circ S_\omega, \ S_\omega \in \Hom(\rho,\sigma)_\omega$ is injective.
Thus $\Hom(\rho,\sigma)_\omega$ is finite dimensional as well.
$\hfill \halmos$



So far we were able to obtain many results under fairly general assumptions
and in a straightforward way. In particular, a picture analogous to the case of
simple units appears, where finite dimensional spaces are replaced by Banach
bundles with finite dimensional fibres and $\circ$ composition of 2-arrow
takes place fibrewise. As already mentioned, the behaviour under the $\otimes$
composition is harder to describe, and in order to do so we introduce some
additional hypothesis.

We begin with the following definition:

\begin{defn}
  We call a 1-arrow $ B \xleftarrow{\rho} A $ ``centrally balanced'' if the
  following holds:
\label{balance} 
$$ \End(\iota_B) \otimes 1_\rho = 1_\rho \otimes \End(\iota_A).$$
\end{defn}
We don't claim that 2-arrows of this form exhaust all of $Z(\End(\rho))$, but simply that
for each $z \in \End(\iota_B) $ there exists $w \in \End(\iota_A)$ s.t. $z \otimes 1_\rho
= 1_\rho \otimes w$ (and vice-versa). 
Notice that this property is closed for $\otimes$ products and sub-objects.

Now we make the following  
\begin{ass} [Balanced Decomposition]
\label{balancedass}
We assume every 1-arrow $\rho$ to be a direct sum $\oplus_i \rho_i$ of
centrally balanced 1-arrows. 
\end{ass}
We have the following
\begin{prop}
Let $ B \xleftarrow{\rho} A$ be a centrally balanced 1-arrow with $S_l(\rho) = \Omega_A$ and
$S_r(\rho) = \Omega_B$. Then $\rho$ establishes an isomorphism of the two algebras
$\End(\iota_A)$ and $\End(\iota_B)$. The isomorphism, which we shall denote
$\theta_{\rho}: \End(\iota_A) \rightarrow \End(\iota_B)$, is independent
of the choice of solutions of the conjugation equation and is given by the expression
$\theta_\rho(w)=  \frac{\ _{{\End}(\iota_B)} 
\langle 1_{\rho} \otimes w, 1_\rho \rangle^{(\rho,\rho)}}  { _{{\End}(\iota_B)}
\langle 1_\rho, 1_\rho \rangle^{(\rho,\rho)}}$.
\end{prop}

{\noindent{\it Proof.}}  The first sentence is just a restatement of Definition \ref{balance}:
the maps 
$s\mapsto s\otimes 1_{\rho}$ and $s\mapsto 1_{\rho}\otimes s$, $s\in\End(\iota_{B})$,
are injective (this follows from Lemma \ref{circleunit} 
and the fact that $R^* \circ R$ and $\bar R^* \circ \bar R$ are invertible, as
$S_l(\rho) = \Omega_A$ and $S_r(\rho) = \Omega_B$). Definition  \ref{balance} tells us that they have the same images in $Z(\End(\rho))$.

In order to give the rest of the proof  we first introduce a simple lemma:
\begin{lem}
For $\rho$ as above one has $w \otimes 1_{\bar \rho \rho} = 1_{\bar \rho \rho} \otimes w$
for all $w \in\End(\iota_A)$.   
\end{lem}
$Proof.$ As there exists $w' \in \End(\iota_A)$ s.t. $w \otimes 1_{\bar \rho \rho} = 1_{\bar \rho \rho} \otimes w'$, we only have to prove that they are the same.
But $w \otimes 1_{\bar \rho \rho} = 1_{\bar \rho \rho} \otimes w'$ implies $ w
\otimes R = R \otimes w'$. And $w \otimes R = R \circ w = R \otimes w$ which
implies $w' = w$, as we have  shown that tensoring elements of
$\End(\iota_A)$ with $R$ is an injective map when $R^* \circ R$ is invertible.\hfill \halmos

 We show that $\theta_\rho(w) \otimes 1_\rho = 1_\rho \otimes w$.
To do so we take the difference of the two elements   and evaluate the product
$\langle(1_{\rho} \otimes w - \theta_\rho(w) \otimes 1_{\rho}), (1_{\rho} \otimes w
- \theta_\rho(w) \otimes 1_{\rho}) \rangle^{(\rho,\rho)}_{{\End}(\iota_A)}$. Making use of the previous lemma shows that the product is zero,
so the two objects must be equal as the $\End(\iota_A)$ inner product is
non-degenerate. Also notice that the right hand side    of  $\theta_\rho(w) \otimes 1_\rho =
1_\rho \otimes w$ does not depend on the choice of solutions $R, \bar R$, thus
the isomorphism $\theta_{\rho}$ must be independent as well. \hfill \halmos

\begin{rem}
In the same way we have an expression for  $\theta^{-1}_{\rho}: \End(\iota_B) 
\rightarrow \End(\iota_A)$ with $\theta^{-1}_{\rho}(z):= \frac{\langle z
\otimes 1_{\rho}, 1_\rho \rangle ^{(\rho,\rho)}_{{\End}(\iota_A)}} {\langle
1_\rho, 1_{\rho} \rangle^{(\rho,\rho)}_{{\End}(\iota_A)}}$.
Also note that $\theta_{\rho} = \theta^{-1}_{\bar \rho}$ and $ \theta_{\rho'} \circ \theta_{\rho} = \theta_{\rho' \otimes \rho}$ (when the composition of arrows  is defined).
\end{rem}

\begin{rem}
  In the preceding proposition we have supposed $S_l(\rho) = \Omega_A$ and
  $S_r(\rho) = \Omega_B$. In the general case an analogous isomorphism holds for the
  sub-algebras $ E_{S_l(\rho)} \otimes \End(\iota_A)$ and $ E_{S_r(\rho)} \otimes
  \End(\iota_B)$, which will be denoted by the same symbol $\theta_\rho:
  E_{S_l(\rho)} \otimes \End(\iota_A) \rightarrow E_{S_r(\rho)} \otimes
  \End(\iota_B)$. The same way, $\theta^{-1 *}$ will indicate the
  homeomorphism between the two subspaces $S_l(\rho) \subset \Omega_A$ and
$S_r(\rho) \subset \Omega_B$.
 \end{rem}
We now describe the behaviour of the $\otimes$ product.

We begin with a remark about the supports of continuous sections.

\begin{lem}
\label{wheresupport}
  Let $ T \in \Hom(\rho,\sigma)$, for $\rho,\sigma$ two generic 1-arrows. Then $\mbox{support T } $ $\subset S_l(\rho) \cap S_l(\sigma).$
\end{lem}

{\noindent{\it Proof.}}  We have $T = 1_{\sigma} \circ T \circ 1_{\rho} = (1_\sigma
\otimes E_{S_l(\sigma)}) \circ T \circ (1_\rho \otimes E_{S_l(\rho)}) = T
\otimes E_{S_l(\sigma)} \otimes E_{S_l(\rho)},$ where in the second equality we
have used Corollary \ref{supportcutting}.  $ \mbox{ } $ \hfill \halmos

\begin{rem}
  The same way we have $ T = E_{S_r(\sigma)} \otimes E_{S_r(\rho)} \otimes T.$
\end{rem}


Let $ T,T' \in (C \xleftarrow {\rho'} B, C \xleftarrow {\sigma'}
B)$ and $ B \xleftarrow {\rho} A$,$ D \xleftarrow{\rho''} C$. Then we have the
following

\begin{cor}
If $S_l(\sigma') \cap S_l(\rho') \cap S_r(\rho) = \emptyset$, then $T \otimes
   1_\rho = 0.$

Analogously, if $S_r(\sigma') \cap S_r(\rho') \cap S_l(\rho'') = \emptyset,$ then
   $1_{\rho''} \otimes T = 0.$
\end{cor}

For the rest of this
section, we will consider only centrally balanced 1-arrows. This is not a real limitation, as we
have supposed Assumption \ref{balancedass} to hold.
\begin{prop}
\label{otimes}
For $\theta^{-1 *}_{\rho}(\omega) \in S_l(\rho'),$  $(T \otimes 1_{\rho})_{ |_{\omega}}= \ (T' \otimes 1_{\rho})_{|_{\omega}}$
  iff $T_{|_{\theta^{-1 *}_{ \rho}(\omega)}} = \ T'_{|_{\theta^{-1 *}_{\rho}(\omega)}}.$

The same way, for $\theta^{-1 *}_{\rho'}(\alpha) \in S_l(\rho''),$  $(1_{\rho''} \otimes T)_{|_{\alpha}}= (1_{\rho''} \otimes T')_{|_{\alpha}}$
iff $T_{|_{\alpha}} = \ T'_{|_{\alpha}}.$
\end{prop}

{\noindent{\it Proof.}}  We fix solutions of the conjugation equations for $\rho, \rho', \rho''$. For the product of the $1-$arrows we take the product of the solutions. For example: $R_{\rho' \otimes \rho}:= (1_{\bar{\rho}} \otimes R_{\rho'} \otimes 1_{\rho}) \circ R_{\rho}$,  $\bar{R}_{\rho' \otimes \rho}:= (1_{\rho'} \otimes \bar{R}_{\rho} \otimes 1_{\bar{\rho'}}) \circ \bar{R}_{\rho'}$.

$(T \otimes 1_{\rho})_{|_{ \omega}}= (T' \otimes 1_{\rho})_{|_{\omega}} $ if and only if 

$  0  \ = \  \langle (T-T') \otimes 1_{\rho},(T-T') \otimes 1_{\rho}\rangle^{(\rho' \otimes \rho \ , \sigma' \otimes \rho)}_{\End(\iota_{A}) \ |_{\omega}}  $

$ = \
\theta^{-1}_{\rho}( \langle
(T-T'),(T-T')\rangle^{(\rho',\sigma')}_{{\End}(\iota_B)})_{|_{\omega}}
\times 
\langle 1_\rho, 1_\rho \rangle^{(\rho,\rho)}_{{\End}(\iota_A) \ |_{\omega}}$.

Thus $T_{|_{\theta^{-1 *}_{ \rho }(\omega)}} = \ T'_{|_{\theta^{-1
      *}_{\rho}(\omega)}},$ as $\langle 1_\rho, 1_\rho
\rangle^{(\rho,\rho)}_{{\End}(\iota_A) \ |_{\omega}} \neq 0$.  
 
The proof of the second statement is analogous. \hfill \halmos

We summarise the situation as follows. Associated to any non zero centrally
balanced 1-arrow (say, $B \xleftarrow{\rho} A)$) we have
homeomorphic  subspaces $S_l(\rho) \subset \Omega_A$ and $S_r(\rho)
\subset \Omega_B$. The conjugation relations give explicit expressions of the isomorphism, depending on the choice of the 1-arrow $\rho$, but not on the choice of the solutions to the conjugation equations.

We are now in the position to give the following
\begin{defn}
\label{local.products}  
We  define a $\circ$ and $\otimes$ product on the fibres:
Let $S \in \Hom(\rho,\sigma),$ $ T \in \Hom(\eta, \rho), Q \in (\rho',\sigma');$ then

$S_{\omega} \circ T_{\omega}:= (S \circ T)_{\omega};$  

$ Q_{\theta^{-1 *}_{\rho}(\omega)} \otimes S_{\omega}:= (Q \otimes
S)_{\omega}, \ \ for \ \omega \in S_l(\rho) \cap S_l(\sigma). $

\end{defn}
The consistency of these definitions is ensured by Propositions \ref{bundleb} and \ref{otimes}.

Notice that the $\circ$ composition is defined between fibres with the same base
point, while the $\otimes$ composition is defined for fibres with base points
in distinct topological spaces, the correspondence given by the homeomorphism
between the supports.   
Thus, we  can think of the $\circ$ product as composing 2-arrows fibrewise and
the $\otimes$ product acting on the fibre structure by ``gluing'' the two
supports of the bundles by means of the homeomorphism
$\theta^{-1 *}$.

As a consequence we have

\begin{cor} Let $\omega \in S_l(\rho),$ then
the map $$ T_{|_{\theta^{-1 *}_{\rho}(\omega)}} \mapsto T_{|_{\theta^{-1
      *}_{\rho}(\omega)}} \otimes 1_{\rho |_{\omega} } $$  is injective; 

the same way, let $ \theta^{-1 *}_{\rho'}(\alpha) \in S_l(\rho''),$ then the
map $$
T_{|_{\alpha}} \mapsto 1_{\rho'' |_{\theta^{-1
      *}_{\rho'} (\alpha)}} \otimes T_{|_\alpha} $$ is injective.
\end{cor}

\begin{cor}
\label{local.conjugation} Let $\omega \in S_l(\rho).$ Then 
$R_{|_{\omega}}$ and $\bar{R}_{|_{\theta^{-1 *}_{\rho}(\omega)}}$ satisfy the conjugation relations: $$(\bar{R}^*_{|_{\theta^{-1 *}_{\rho}(\omega)}} \otimes 1_{\rho |_{\omega}}) \circ (1_{\rho |_{\omega}} \otimes R_{|_{\omega}}) = 1_{\rho |_{\omega}} ; $$ 
$$(R^*_{|_{\omega}} \otimes 1_{\bar{\rho} |_{\theta^{-1 *}_{\rho}(\omega)}}) \circ (1_{\bar{\rho} |_{\theta^{-1 *}_{\rho}(\omega)}} \otimes \bar{R}_{|_{\theta^{-1 *}_{\rho}(\omega)}}) = 1_{\rho |_{\theta^{-1 *}_{\rho}(\omega)}}.$$
\end{cor}

\begin{cor}
The map $\bullet$ induces a conjugate linear isomorphism, which we will indicate with the
same symbol, between  the fibres 
$ \bullet: \Hom(\rho,\sigma)_{|_{\omega}} \rightarrow \Hom(\bar{\rho},\bar{\sigma})_{|_{\theta^{-1 *}_{\rho}(\omega)}}$.
\end{cor}


\begin{cor}
\label{pi-po.om}
Let $X_{|_{\omega}} \in \End(\rho)_{|_{\omega}}$ be  positive . Then the following inequality holds in $\End(\rho)_{|_{\omega}}$: $ X_{|_{\omega}} \leq   (\bar{R}^* \circ \bar{R})_{|_{\theta^{-1 *}_{\rho}(\omega)}} \otimes  1_{\rho |_{\omega}} \otimes (R^*_{|_{\omega}} \circ (1_{\bar{\rho} |_{\theta^{-1 *}_{\rho}(\omega)}} \otimes X_{|_{\omega}}) \circ R_{|_{\omega}})$.
\end{cor}

These assertions are the ``local'' version of the ones already encountered in the
introduction. For example Corollary \ref{pi-po.om} is valid as long as $X_{|_{\omega}}$ is
positive, independently of the value $X$ in other points (thus, $X$ need not be positive
in $\End(\rho)$).

\begin{lem}
  Let $T \in (B \xleftarrow{\rho} A, B \xleftarrow{\sigma} A)$, with $\rho$ and
  $\sigma$ centrally balanced $1$-arrows. Then $\forall
  \omega \in \Omega_A \ | \ T_{\omega} \neq 0 \ \ $ the homeomorphisms between
  $\Omega_A$ and $\Omega_B$ induced by $\rho $and $\sigma$  coincide, i.e.\ $ \ \ \theta^{* -1}_{\rho}(\omega) =
  \theta^{* -1}_{\sigma}(\omega)$.
\end{lem}
{\noindent{\it Proof.}}  For any $w \in \End(\iota_A)$ we have $$T \otimes w = T \circ
(1_{\rho} \otimes w) = T \circ (\theta_{\rho}(w) \otimes 1_{\rho})
= \theta_{\rho}(w) \otimes T$$ 
and in the same way we have 
$$T
\otimes w = (1_{\sigma} \otimes w) \circ T = \theta_{\sigma}(w)
\otimes 1_{\sigma}) \circ T = \theta_{\sigma}(w) \otimes T.$$
In particular 
$$ (\theta_{\rho}(w) \otimes T)_{|_{\omega}} =
((\theta_{\sigma}(w) \otimes T)_{|_{\omega}}
$$ 
which by Proposition \ref{otimes} is equivalent to 
 $$\theta_{\rho} (w)_{|_{\theta^{* -1}_{\rho}(\omega)}} \otimes
 T_{|_{\omega}}= \theta_{\sigma} (w)_{|_{\theta^{* -1}_{\rho}(\omega)}}
 \otimes T_{|_{\omega}},$$ which in turn gives $\theta_{\rho} (w)_{|_{\theta^{* -1}_{\rho}(\omega)}} = \theta_{\sigma} (w)_{|_{\theta^{* -1}_{\rho}(\omega)}}$.
As this must hold for any $w \in \End(\iota_A)$ we have $ \theta^{*
 -1}_{\rho}(\omega) = \theta^{* -1}_{\sigma} (\omega)$. \hfill \halmos

We show now that starting  from a 2-$C^*$-category $ \mathcal{C}$, for  each
$\omega_0 \in \Omega_A$, where $A$ is an arbitrary object, it is possible to
construct a  
2-$C^*$-category, which we will indicate by $\mathcal{C}^{\omega_0 , A}$, with simple units.
Consider the full sub 2-$C^*$-category of $\mathcal{T}$ of $C$ generated by
centrally balanced $1$-arrows, their products and their sub$1$-arrows.
 This way to each $B \xleftarrow{\rho} A \  \in \mathcal{T}$ will be associated  a homeomorphism
 $\theta^{* -1}_{\rho}:  S_l(\rho) \rightarrow S_r(\rho).$

Now define:
\begin{itemize}
\item as objects the set $\{ \omega_B  \in \Omega_B, \ \forall B \in \mathcal{C} \}$ such that there exists  a $ B \xleftarrow{\rho} A
  \ \in \mathcal{T}$ with $ \omega_0 \in S_l(\rho), \omega_B \in S_r(\rho)$
  and $\omega_B = \theta^{* -1}_{\rho}(\omega_0)$, 
\item as $1$-arrows  for any $ \omega_B,  \omega_C $ as above  a $
  \omega_C \xleftarrow{\sigma} \omega_B$ in correspondence to any $C
  \xleftarrow{\sigma} B $ with $ \omega_B \in S_l(\sigma), \omega_C \in
  S_r(\sigma),$ verifying $\theta^{* -1}_{\sigma}(\omega_B) = \omega_C,$
\item as $2$-arrows, for any $\omega_C \xleftarrow{\sigma} \omega_B, \ \omega_C
  \xleftarrow {\eta} \omega_B$ as above,  $\Hom(\omega_C \xleftarrow{\sigma} \omega_B, \omega_C \xleftarrow{\eta} \omega_B)$ to be the set $ \{ T_{|_{\omega_B}}, \ \forall T \in
  \Hom(\sigma,\eta) \}$.
   
\end{itemize}

\begin{rem}
  By the preceding lemma we see that for $\omega_C \xleftarrow{\sigma}
  \omega_B \ , \ \omega'_C \xleftarrow{\eta} \omega_B \ , \ \omega_C \neq
  \omega'_C \,\ \forall T \in  \Hom(\sigma,\eta) \ T_{|_{\omega_B}} = 0$, i.e.\  we
  don't loose any information by  considering $\omega_C$ and $\omega'_C$  as distinct objects in the new category.
\end{rem}

We define the $\circ$ and $\otimes$ products of Definition \ref{local.products}
as the operations of our new category and we endow the spaces of $2$-arrows
with the fibre $C^*$-norm introduced above, thus obtaining a 2-$C^*$-category. For each object $\omega_B$ we have the
$1$-unit $\omega_B \xleftarrow{\iota_{\omega_B}} \omega_B$ corresponding to $\iota_B$. For
each $1$-arrow $\omega_C \xleftarrow{\sigma} \omega_B$ we have the $2$-unit
$1_{\sigma |_{\omega_B}}$. As $$\End(\omega_B \xleftarrow{\iota_{\omega_B}}
\omega_B) \cong
\End(\iota_B)_{\omega_B} \cong \mathbb{C},$$ we see that each $1$-unit is simple. Corollary
\ref{local.conjugation} 
ensures that this category is closed under conjugation. 
Closure with respect to projections is not automatically ensured. For example, if
$P^{\omega_B} \in (\omega_C \xleftarrow{\sigma} \omega_B,\omega_C \xleftarrow{\sigma}
\omega_B)$ is a projection, there does not  necessarily exist a projection (thus a
corresponding $1$-arrow) $P$ in $\Hom(\sigma,\sigma)$ such that $P_{|_{\omega_B}}
= P^{\omega_B}$. 
We will consider the completion under projections of the above category, and
denote it by $\mathcal{C}^{\omega_0,A}$.

Note that in $\mathcal{C}^{\omega_0,A}$ several points $\omega_B, \omega'_B
\dots$ of $\Omega_B$ may appear as distinct objects.
In fact, as the various maps $\theta^{* -1}_{\rho}: S_l(\rho) \rightarrow
S_r(\rho) $ are invertible, we see that each point $\omega_0$ determines an
orbit in the spaces $\Omega_A, \Omega_B, \dots$, and that the construction
leading to $\mathcal{C}^{\omega_0,A}$ depends only on the choice of one of these (disjoint) orbits.

\begin{ex}

Let  $B \xleftarrow{\rho} A, A
\xleftarrow{\bar{\rho}} B$, be a pair
of centrally balanced , conjugate $1$-arrows with $E_{S_l(\rho)} = \Omega_A$
and $ E_{S_r(\rho)} = \Omega_B$. We can consider the full sub
2-$C^*$-category generated by these two elements, i.e.\  their compositions $
\rho, \ \rho \otimes \bar \rho, \ 
\rho \otimes \bar \rho \otimes \rho \dots $ and
their sub-$1$-arrows. This is  the categorical analogue of Jones' basic
construction. 
We can take the sequence of algebras 
$$\End(\rho), \ \End(\rho \otimes \bar \rho), \ \End( \rho \otimes \bar \rho \otimes \rho), \dots.$$
and realise a sequence of injective inclusions as follows: $$\End(\rho) \ni X \hookrightarrow X \otimes 1_{\bar \rho} \in \End(\rho \otimes \bar \rho),$$
$$ \End(\rho \otimes \bar
\rho) \ni Z \hookrightarrow Z \otimes 1_\rho \in \End(\rho \otimes \bar \rho \otimes
\rho), \dots $$
We choose $\omega_0 \in \Omega_A$ and construct $\mathcal{C}^{\omega_0,A}$ as
above. In this case we have only two objects, namely $\omega_0$ and $\theta^{*
  -1}_{\rho}(\omega_0)$, as $\theta^*_{\bar \rho} = \theta^{* -1}_{\rho}$  implies
$\theta^{* -1}_{\rho \otimes \bar \rho} = \theta^{* -1}_{\rho} \circ \theta^{*
  -1}_{\bar \rho} = id$. Thus, for
example, $$\theta^{* -1}_{\rho \otimes \bar \rho \otimes \rho}(\omega_0)
\xleftarrow{\rho \otimes \bar \rho \otimes \rho} \omega_0 =  \ \theta^{* -1}_{\rho}(\omega_0) \xleftarrow{\rho
  \otimes \bar \rho \otimes \rho} \omega_0 ,$$  
and analogously for the other $1$-arrows of $\mathcal{A}^{\omega_0,A}$.
As in the case of subfactors, we have a sequence of inclusions of finite
dimensional $C^*$-algebras:
$$\End(\rho)_{\omega_0} \ni X_{| _{\omega_0}} \hookrightarrow X_{| _{\omega_0}} \otimes 1_{\bar \rho |_{\theta^{* -1}_{\rho}(\omega_0)}}  \in \End(\rho \otimes \bar \rho)_{\theta^{* -1}_{\rho}(\omega_0)}$$
$$\End(\rho \otimes \bar \rho)_{\theta^{*
    -1}_{\rho}(\omega_0)} \ni Z_{| _{\theta^{* -1}_\rho (\omega_0)}} 
\hookrightarrow Z_{| _{\theta^{* -1}_\rho (\omega_0)}} \otimes 1_{\rho |_{\omega_0}}
\in \End(\rho \otimes \bar \rho \otimes \rho)_{\omega_0} \dots $$
\end{ex}

\begin{rem}
  A particular case is that of an $\End(\iota)$-linear tensor $C^*$-category
  $\mathcal{T}$, where $\iota$ is the unit element. $\End(\iota)$-linear means
  that for any $k \in End(\iota) \cong C(\Omega)$ and for any $\rho \in
  \mathcal{T}$, one has $\rho \otimes k = k \otimes \rho.$ Clearly our Assumption
  \ref{balancedass} holds, namely every object of the tensor category is 
  centrally balanced: for each object $\rho$ the homeomorphism $\theta^{-1
  *}_\rho: S_l(\rho) \rightarrow S_r(\rho) = S_l(\rho)$ is the identity map. The above
  construction assumes a clear form in this case, which we are tempted to
  describe as that of a ``bundle of tensor $C^*$-categories with simple units''
  over the topological space $\Omega$: 
\end{rem}

\begin{prop}
\label{fibre of tensor categories}
  Let $\mathcal{T}$ be an $\End(\iota)$-linear tensor $C^*$-category closed for
  conjugation, sub-objects and direct sums. Let
  $\Omega$ be the compact Hausdorff topological space associated to
  $\End(\iota) \cong C(\Omega).$ Then for each  $\omega \in \Omega$ there is
  an associated tensor $C^*$-category $\mathcal{T}_\omega$ with simple unit
  object, closed for conjugation and direct sums, the fibre category at the
  point $\omega$. The arrows of the original
  category $\mathcal{T}$ can be viewed as continuous sections taking values in
  the $Hom$ spaces of the fibre categories. The $\circ, \otimes, * $
  operations, as well as the conjugation operation of $\mathcal{T}$, preserve the fibre structure.
\end{prop}



\section{Standard solutions}
\label{standardsolutions}
In this section we will introduce  a particular class of solutions to the
conjugation equations.

We begin by recalling the analogous definition and some basic facts
concerning the case of simple units, contained
in \cite{LR97}. 
So, for the moment, we suppose that $\rho$ is a 1-arrow going from $A$ to $B$
and that $\End(\iota_A) \cong \mathbb{C}$, $\End(\iota_B) \cong \mathbb{C}$.

Let $\bar \rho$ be a conjugate 1-arrow going from $B$ to $A$.
The main difference  (and simplification) from the general case is that the algebra $\End(\rho)$
and its isomorphic $\End(\bar \rho)$ are finite dimensional $C^*$-algebras, i.e.\  a finite direct sum of matrix
algebras. It is thus possible to decompose the unit arrow $1_{\rho}$ as a sum
of minimal projections $e_i$ in the algebra $\End(\rho)$. To each of these minimal projections will
correspond an irreducible $\rho_i$, i.e.\  $\End(\rho_i) \cong \mathbb{C}$,  and
we think of   $\rho$ as a direct sum  $ \oplus_i \rho_i$, i.e.\  there exists a
complete family of  isometric 2-arrows $W_i \in \Hom(\rho_i, \rho)$
s.t. $W_i \circ W^*_i = e_i$ , $W^*_i \circ W_i = 1_{\rho_i} $. Two projections
$e_i,e_j$ dominated by the same minimal central projection in $\End(\rho)$ will lead to equivalent
irreducible 1-arrows, i.e.\  there will exist a unitary $V_{i,j} \in\Hom (\rho_j)$.

For  irreducible $\rho$, i.e.\  $\End(\rho) \cong  \mathbb{C}$, we give the
following definition (notice that by duality, $\bar \rho$ is 
irreducible too):

\begin{defn}
  Let $\rho, \bar \rho$ be irreducible. $R,\bar R$ are said to be a standard
  solution if $R^*\circ R = \bar R^* \circ \bar R$.
\end{defn}

Let now $\rho$ be not necessarily irreducible and $\oplus_i \rho_i$ its
decomposition into irreducibles. Let \{ $\bar \rho_i$ \} be irreducibles, each
conjugate to its corresponding $\rho_i$ and $R_i, \ \bar R_i$ standard solutions
for each of these couples . 

Then $\oplus_i \bar \rho_i$ and  $ \rho $ are conjugate and $ \oplus_i R_i ,
\oplus_i \bar R_i$ is a solution. 

\begin{defn}
\label{standsol}
  Let $\rho = \oplus \rho_i$ and $\bar \rho = \oplus_i \bar \rho_i$ be
  conjugate. A solution of the form $\oplus_i R_i$, $\oplus_i \bar R_i$, where the $R_i$ and
  $\bar R_i$ are  irreducible and standard as defined above, is called a standard solution. 
\end{defn}
  
As we will see, standard solutions are uniquely defined up to a unitary in $\End(\rho)$.

In this context standard solutions always exist:

\begin{prop}
\label{discreteconiug}
 Let $\rho$ and $\bar \rho$ be conjugate 1-arrows between objects $A$ and $B$,
 with $\End(\iota_A) \cong \mathbb{C}$ and $\End(\iota_B) \cong \mathbb{C}$.
Then  standard solutions always exist.  
\end{prop}

{\noindent{\it Proof.}} 
Choose an arbitrary 
solution $R',\bar R'$. Decompose $\rho$ into a sum of irreducibles $\oplus_i
\rho_i$ by means of a complete family of orthogonal projections $\{ e_i \}$
in $\End(\rho)$.  The corresponding elements $e^{\bullet}_i \in\End (\bar\rho)$ will be a
complete family of disjoint 
idempotents(i.e.\  $e_i^{\bullet} e^{\bullet}_j = \delta_{i,j} e^{\bullet}_i$), as the
$\bullet$ operation is an anti-isomorphism between the algebras $\End(\rho)$
and $\End(\bar \rho)$. But they will fail, in general,  to be self
adjoint. Nevertheless there will exist an invertible $A \in \End(\bar\rho)$ such that 
$\{\bar e_i:= A e^{\bullet}_i A^{-1} \}$ is a complete family of orthogonal projections. 
As before, we can decompose $\bar \rho = \oplus_i \bar \rho_i$, where each
$\bar \rho_i$ corresponds to the projection $\bar e_i$. Taking now $R'':= (A \otimes
1_{\rho}) \circ R'$ and $\bar R'':= (1_{\rho} \otimes A^{-1 *}) \circ \bar R'$ as
solutions, we see that each pair $\rho_i$ and $\bar \rho_i$ is a couple of
conjugates, with $R''_i:= (\bar e_i \otimes e_i) \circ R''$, $\bar R''_i:= ( e_i \otimes \bar e_i)
\circ \bar R''$ as solutions.

Furthermore  each pair $R''_i, \bar R''_i$ is determined up to a constant as $\rho_i$ and $\bar
\rho_i$ are irreducible. We rescale each couple  by a  constant so to have the equality
$R^{'' *}_i \circ R''_i = {\bar R}_i^{''*} \circ \bar R''_i$. This amounts to multiplying
the invertible element $A \in\End ( \bar \rho)$ by a diagonal element 
$D$ in the same algebra. 
We set $R:= (D \otimes 1_{\rho}) \circ R''$ and $\bar R:= (1_{\rho} \otimes
D^{-1}) \circ \bar R''$ and $R_i:= (\bar e_i \otimes e_i) \circ
R$ , $ \bar R_i:= (e_i \otimes \bar e_i) \circ \bar R$. It is easily verified
that $R = \oplus_i R_i$ and $\bar R = \oplus_i \bar R_i $.\hfill \halmos

The number $R^*_i \circ R_i = \bar R^*_i \circ \bar R_i$ is called
the  dimension of the irreducible $\rho_i$ and depends only on the equivalence class of the
$\rho_i$. In fact, suppose $\rho_i$ and $\rho_j$ are equivalent, i.e.\  there
exists a partial isometry  $V_{i,j} \in \End(\rho)$ such that $V^*_{i,j} \circ V_{i,j} =
e_i$ and $V_{i,j} \circ V^*_{i,j} = e_j$ (in other words, $e_i$ and $e_j$ have same central
support). Then, by duality, the same will hold for $\bar \rho_i$
and $\bar \rho_j$, with a partial isometry $\bar V_{i,j} \in\End (\bar \rho)$. 
Take $(\bar V_{i,j} \otimes V_{i,j}) \circ R_i $ and $(V_{i,j} \otimes
\bar V_{i,j}) \circ \bar R_i $ as a solution for $\rho_j$ and $\bar \rho_j$.
This solution can differ from $R_j$,$\bar R_j$  only by a 
invertible element in $\End(\rho_j) = \mathbb{C}$, i.e.\  $(\bar V_{i,j} \otimes V_{i,j})
\circ R_i = \lambda R_j$ and 
$(V_{i,j} \otimes \bar V_{i,j}) \circ \bar R_i = \lambda^{-1} \bar R_j$.
But this implies $(R_i^* \circ R_i) ({\bar R}_i^ * \circ \bar R_i) =
(R_j^* \circ R_j) ({\bar R}_j^* \circ \bar R_j)$ and , because of the
normalization we have chosen above, $R^*_i \circ R_i = R^*_j \circ
R_j.$ 

The number $R^* \circ R = \sum_i R^*_i \circ R_i $ is the  
dimension of the object $\rho$. It is additive with respect to direct sums.

The class of  standard solution defines a trace on the algebra $\End(\rho)$
(and in an analogous way on the algebra $\End(\bar \rho)$). 
Consider an element $ S \in \End(\rho) $. Thinking of $\End(\rho)$ as a direct
sum of matrix algebras, we indicate by $S_{i,j}$ the matrix
elements of $S$ corresponding to a representation given by the projectors ${
\{ e_i \} }$.
Thus $(\bar R^* \circ ( S \otimes 1_{\bar \rho} ) \circ \bar R)   = \sum_i S_{i,i} \bar R^*_i \circ \bar R_i   = 
  \sum_i S_{i,i} R^*_i
\circ R_i =   R^* \circ (1_{\bar \rho} \otimes S)\circ R$ which is a trace, as the
 dimensions $R^*_i \circ R_i$ depend only on the central supports of the
 $e_i$.

Notice that we have shown that  for a standard solution the following holds:

$$ R^* \circ (1_{\bar \rho} \otimes S) \circ R = \bar R^* \circ (S \otimes
1_{\bar \rho}) \circ \bar R, \ \ \forall S \in \End(\rho).$$
This can be used as an equivalent definition for standardness (see \cite{LR97}).

\begin{rem}
We can summarize the situation as follows: fix a faithful normalised trace
$tr$ (i.e.\  $tr(1) = 1$) on the algebra $\End(\rho)$. A generic solution $R',
\bar R'$ to the conjugation relations  will induce faithful functionals
$\Phi'(S):= R'
\circ (1_{\rho} \otimes S) \circ R' =  tr(D K S)$ ,  $ \Psi^{'}( S ):=  \bar
R' \circ (S \otimes 1_{\rho}) \circ R' = tr(D K^{-1} S)$, where $D$ is a positive
invertible central element in $\End(\rho)$, whose trace is the dimension of
$\rho$, and $K$ is a positive invertible element of $\End(\rho)$.  
Taking $R:= (1_{\bar \rho} \otimes K^{-1}) \circ R'$ and $\bar R:= (K
\otimes 1_{\bar \rho}) \circ \bar R'$ will give a standard solution.
Notice that the element $K$ is uniquely defined by the original choice $R',
\bar R'$. The arbitrariness of the choice of the trace is expressed by the 
central element $D$.   
\end{rem}

We  now drop the hypothesis of simple units. We begin by  giving   the
definition of standardness for centrally balanced 1-arrows.

\begin{defn}
Let $\rho, \bar \rho$ be centrally balanced. Let $R, \bar R$ be a solution to the
conjugation equations. We say $R, \bar R$ to be standard if $\forall X \in \End(\rho)$ the
following holds:  $$ 1_{\rho} \otimes (R^* \circ (1_{\bar \rho}
\otimes X) \circ R) =   ({\bar R}^* \circ (X \otimes 1_{\bar \rho})
\circ \bar R) \otimes 1_\rho .$$
    
\end{defn}

The following lemma shows that another appropriate name could have been ``minimal'':

\begin{lem}
Let $R,\bar R$ be  standard solutions  for centrally balanced $\rho , \bar \rho$. Then for any other solution $R',\bar R'$
we have: $  (\bar {R'}^* \circ \bar {R'}) \otimes
1_{\rho} \otimes ({R'}^* \circ {R'})  \geq   (\bar {R}^* \circ \bar R) \otimes 1_{\rho}
\otimes (R^* \circ R) .$ 
\end{lem}

{\noindent{\it Proof.}} 
We have 
$R' = (1_{\bar \rho} \otimes X) \circ R$ and  $\bar R' = (X^{-1 *}\otimes
1_{\bar \rho}) \circ \bar R$ for some invertible $X \in \End(\rho)$.
Thus  $$(\bar R'^* \circ \bar R') \otimes 1_{\rho} \otimes (R'^* \circ R') 
 = $$
$$ (\bar R^* \circ (X^{-1 } \circ X^{-1 *}) \otimes 1_{\bar \rho} \circ \bar R) \otimes 1_{\rho} \otimes (R^* \circ 1_{\bar \rho} \otimes (X^* \circ
X)\circ R)  = $$
$$  1_{\rho} \otimes ((R^* \circ 1_{\bar \rho}
\otimes (X^* \circ X) \circ R)  \circ  (R^* \circ 1_{\bar \rho}
\otimes (X^{-1} \circ X^{-1 *}) \circ R)).$$

The claim is implied by the inequality  $$(R^* \circ 1_{\bar \rho}
\otimes (X^* \circ X) \circ R)  \circ (R^* \circ 1_{\bar \rho}
\otimes (X^{} \circ X^{-1 *}) \circ R) \geq (R^* \circ R)^{2},$$ which can be written
equivalently $$ \langle X,X\rangle_{\End(\iota_A)}   \langle X^{-1 *} ,X^{-1 *}\rangle_{(\iota_A ,
  \iota_A)}  \geq \langle 1_\rho, 1_\rho \rangle^2_{{\End}(\iota_A)}$$ \ (where we have used the $ \langle  \ , \ \rangle_{\End(\iota_{A})}$ inner product defined by $R$
and $\bar R$). 

It is sufficient to prove this inequality for each $\omega \in
\Omega_{A}$, which is easily done by means of a  Cauchy-Schwarz argument:
first rewrite $ \langle X,X\rangle_{{\End}(\iota_A)}$ and $ \langle X^{-1 *},X^{-1
  *}\rangle_{{\End}(\iota_A)}$ as $ \langle (X^* \circ X)^{\frac{1}{2}},(X^* \circ
X)^{\frac{1}{2}}\rangle_{(\iota_A.\iota_A)}$ and $ \langle (X^* \circ X)^{- \frac{1}{2}},
(X^* \circ X)^{- \frac{1}{2}}\rangle_{{\End}(\iota_A)}$ respectively.
Then $$ \langle (X \circ X^*)^{\frac{1}{2}},(X \circ
X^*)^{\frac{1}{2}}\rangle_{{\End}(\iota_A)}  \langle (X \circ X^*)^{- \frac{1}{2}},(X
\circ X^*)^{- \frac{1}{2}}\rangle_{{\End}(\iota_A)}$$ $$ \geq \ \ \ (  \langle 1_\rho,
1_\rho\rangle_{{\End}(\iota_A)})^{2}. \ \ \  \hfill \halmos$$

\begin{lem}
\label{stun}
Suppose that  $R,\bar R$ and $R',\bar R'$ are two pairs of  standard
solutions. Then there exists a unitary $U \in \End(\rho)$ such that $R' =
(1_{\bar \rho}\otimes U) \circ R$ and $\bar R' = (U \otimes 1_{\bar \rho}) \circ \bar R$. 
\end{lem}
{\noindent{\it Proof.}} 
As we have seen, there exists an invertible $U$ satisfying $R' =
(1_{\bar \rho}\otimes U) \circ R$ and $\bar R' = (U^{-1 *} \otimes 1_{\bar
  \rho}) \circ \bar R$.   We must prove that $U$ is unitary.

By the definition of standardness we have $$(\bar R \circ A \otimes 1_{\bar
  \rho} \circ \bar R) \otimes 1_\rho = 1_\rho \otimes (R^* \circ 1_\rho
  \otimes A \circ R), \ \forall A \in \End(\rho).$$
This implies $(\bar R \circ (U^{-1} A U^{-1}) \otimes 1_{\bar \rho} \circ \bar
  R) \otimes 1_\rho = 1_\rho \otimes (R^* \circ 1_{\bar \rho} \otimes (U^* A
  U^{-1 *}) \circ R).$ As $R'$ and $\bar R'$ are standard solutions as well,
  we have $$ (\bar R \circ(U^{-1} A U^{-1 *}) \otimes 1_{\bar \rho} \circ \bar
  R) \otimes 1_\rho = 1_\rho \otimes (R^* \circ 1_{\bar \rho} \otimes (U^* A
  U) \circ R), \ , \ \forall \in \End(\rho).$$
This implies $$ \langle U,AU\rangle_{{\End}(\iota_A)} \  = \  \langle U^{-1 *},A
  U^{-1 *}\rangle_{{\End}(\iota_A)} \ \forall A \in \End(\rho).$$
In particular $ \langle U,U\rangle_{{\End}(\iota_A)} \ = \  \langle U^{-1 *},U^{-1
  *}\rangle_{{\End}(\iota_A)}$. $R'$ and $\bar R'$ also satisfy the minimality
  condition of the preceding lemma, i.e.\  $ (\bar R^{' *} \circ R^{' *})
  \otimes 1_\rho \otimes (R^{' *} \circ R') = (\bar R^* \circ \bar R) \otimes 1_\rho
  \otimes (R^* \circ R)$. In other words we have $ \langle U,U\rangle_{{\End}(\iota_A)} \ =
  \  \langle U^{-1 *}, U^{-1 *}\rangle_{{\End}(\iota_A)} = \langle 1_\rho,1_\rho \rangle_{{\End}(\iota_A)}$ (where we are using as before
  the
  $ \langle \ ,\ \rangle_{{\End}(\iota_A)}$ inner product relative to the $R, \bar R$ solution). 
Noticing that $$ \langle U^* U, U^* U\rangle_{{End}(\iota_A)} \ = \  \langle U U^* U,
  U\rangle_{{\End}(\iota_A)}$$ $$ = \  \langle U^{-1 *} U^* U, U^{-1
  *}\rangle_{{\End}(\iota_A)} \ =
  \  \langle 1_\rho,1_\rho\rangle_{{\End}(\iota_A)}$$ we see that the product $ \langle (1 - U^*
  U),(1-U^* U)\rangle_{{\End}(\iota_A)}$ is zero.

Thus $U = U^{-1 *}$, as the $\End(\iota_A)$-valued inner product is non-degenerate.
$\hfill \halmos$

\begin{rem}
  In an analogous way one can show that there exists a unitary $\bar U \in
  (\bar \rho, \bar \rho)$ such that $R' = \bar U \otimes 1_\rho \circ R$ and
  $\bar R ' = 1_\rho \otimes \bar U \circ \bar R$.
\end{rem}

\begin{prop}
\label{sttr}
  Let $R,\bar R$ be a standard solution for $\rho,\bar \rho$. Then the
  associated inner product $ \langle \cdot,\cdot\rangle_{{\End}(\iota_A)}$ is tracial,
  i.e.\  $$ \langle ST,1_\rho\rangle_{{\End}(\iota_A)} \ = \  \langle TS,1_\rho\rangle_{{\End}(\iota_A)}
  \ \forall S,T \in \End(\rho).$$ 
\end{prop}
{\noindent{\it Proof.}}  It suffices to prove $$ \langle U^*SU,1_\rho\rangle_{{\End}(\iota_A)} \ =
\    \langle S,1_\rho\rangle_{{\End}(\iota_A)}$$ for any unitary $U \in
\End(\rho), \ \forall S \in \End(\rho).$ 
We have $$ \langle U^*SU,1_\rho\rangle_{{\End}(\iota_A)} \ =  \  R^*
  (1_{\bar \rho} \otimes U^*SU) \circ  R $$ $$ = \  R^* \circ (1_{\bar \rho}
  \otimes U^*) \circ (1_{\bar \rho} \otimes S) \circ (1_{\bar \rho} \otimes U)
  \circ R.$$ But $R':= (1_{\bar \rho} \otimes U) \circ R$ ,$ \bar R':= \ (U \otimes 1_{\bar
  \rho}) \circ \bar R $ is still a standard solution, so there exists a 
unitary $ \bar U \ \in (\bar \rho,\bar \rho)$ such that $ R' = (\bar U \otimes
  1_\rho) \circ R$ and $\bar R' = (1_\rho \otimes \bar U) \circ \bar R$.
 Thus $$ \langle USU^*,1_\rho\rangle_{{\End}(\iota_A)} \ = \  R^* \circ (\bar U^* \otimes
 1_\rho) \circ (1_{\bar \rho} \otimes S) \circ (\bar U \otimes 1_\rho) \circ R $$
$$ = R^* \circ (1_{\bar \rho} \otimes S) \circ R \ = \
   \langle S,1_\rho\rangle_{{\End}(\iota_A)}.$$ \hfill \halmos 

  \begin{rem}
    Analogously one proves that the $_{{\End}(\iota_B)} \langle \ , \ \rangle$ inner product
    is tracial too. 
  \end{rem}

\begin{lem}
\label{stpr}
  Let $\rho$ and $\sigma$ be two centrally balanced 1-arrows with standard solutions  $R_{\rho},
  \bar R_{\rho}$ and $R_{\sigma}, \bar R_{\sigma}$ respectively. Then
  the product solution, defined by
$$R_{\sigma \otimes \rho}:= 1_{\bar \rho} \otimes R_{\sigma} \otimes
1_{\rho} \circ R_{\rho}, \quad \bar R_{\sigma \otimes \rho}:= 1_{\sigma} \otimes
\bar R_{\rho} \otimes 1_{\bar \sigma} \circ \bar R_{\sigma},$$
is standard.
\end{lem}

{\noindent{\it Proof.}} 
For every $A \in\End (\sigma \otimes \rho)$ we have: 
$$ 1_{\sigma \otimes \rho} \otimes (R^*_{\sigma \otimes \rho} \circ (1_{\bar \rho \otimes \bar \sigma} \otimes A)
\circ R_{\sigma \otimes \rho}) 
= $$
$$  1_{\sigma} \otimes 1_{\rho} \otimes (R_{\rho}^* \circ(1_{\bar \rho} \otimes R_{\sigma}^* \otimes 1_{\rho}) \circ (1_{\bar \rho}
\otimes 1_{\bar \sigma} \otimes A) \circ  (1_{\bar \rho}
\otimes R_{\sigma} \otimes 1_{\rho}) \circ R_{\rho}).$$

Using the standard solution property of $R_{\rho}, \bar R_{\rho}$ on the
 element $R_{\sigma}^* \otimes 1_{\rho} \circ (1_{\bar \sigma} \otimes A) \circ
 R_{\sigma} \otimes 1_{\rho} \in \End(\rho)$ we get:

 $$1_{\sigma} \otimes  1_{\rho} \otimes (R_{\rho}^* \circ(1_{\bar \rho} \otimes R_{\sigma}^* \otimes 1_{\rho}) \circ (1_{\bar \rho}
\otimes 1_{\bar \sigma} \otimes A) \circ  (1_{\bar \rho}
\otimes R_{\sigma} \otimes 1_{\rho}) \circ R_{\rho})
 = $$
$$1_{\sigma}  \otimes ( \bar R_{\rho}^* \circ ((R_{\sigma}^* \otimes 1_{\rho} \circ 1_{\bar \sigma} \otimes A \circ
 R_{\sigma} \otimes 1_{\rho}) \otimes 1_{\bar \rho}) \circ \bar R_{\rho}) \otimes
 1_{\rho} = $$
$$1_{\sigma} \otimes (R_{\sigma}^* \otimes \bar R_{\rho}^* \circ (1_{\bar \sigma } \otimes A
 \otimes 1_{\rho}) \circ R_{\sigma} \otimes \bar R_{\rho}) \otimes 1_{\rho}. $$

The same reasoning applied to the solution $R_{\sigma},\bar R_{\sigma}$ and
the element $1_\sigma \otimes  \bar R_{\rho}^* 
\circ (A
\otimes 1_{\bar \rho} )\circ 1_\sigma \otimes \bar R_{\rho} \in
\End(\sigma)$ gives:
$$ 1_\sigma \otimes (R_{\sigma}^* \otimes \bar R_{\rho}^* \circ (1_{\bar \sigma } \otimes A
 \otimes 1_{\rho}) \circ R_{\sigma} \otimes \bar R_{\rho}) \otimes 1_{\rho} = $$
$$(\bar R_{\sigma}^* \circ (1_{\sigma} \otimes \bar R_{\rho}^* \otimes 1_{\bar
  \sigma}) \circ (A \otimes 1_{\bar \rho} \otimes 1_{\bar \sigma}) \circ (1_{\sigma} \otimes \bar R_{\rho} \otimes 1_{\bar
   \sigma}) \circ \bar R_{\sigma}) \otimes 1_{\sigma} \otimes 1_{\rho}. $$
Thus $$1_{\sigma \otimes \rho} \otimes (R_{\sigma \otimes \rho}^* \circ (1_{\bar
  \rho \otimes \bar \sigma} \otimes A) \circ R_{\sigma \otimes \rho}) = (\bar
R_{\sigma \otimes \rho}^* \circ (A \otimes 1_{\sigma \otimes \rho}) \circ
\bar R_{\sigma \otimes \rho}) \otimes 1_{\sigma \otimes \rho}.$$ \hfill \halmos

We now give the natural definition of standardness for a generic 1-arrow.

\begin{defn}
  Let $\rho$ and $\bar \rho$ be a direct sum of centrally balanced 1-arrows
  $\oplus_i \rho_i$ and $\oplus_i \bar \rho_i$ respectively,
  the decomposition given by complete sets of orthogonal partial isometries $\{ W_i \in
\Hom (\rho_i,\rho) \}$ and $\{ \bar W_i \in \Hom(\bar \rho_i,\bar \rho) \}.$ A
  solution  $R,
  \bar R$ to the conjugation equations for $ \rho, \bar \rho $ is said to be
  standard if it is of the form $ \sum_i (\bar W_i \otimes W_i) \circ R_i ,
   \ \sum_i (W_i \otimes \bar W_i) \circ \bar R_i$, where each couple $R_i,\bar R_i$ is
  standard for  $\rho_i, \bar \rho_i$, respectively .   
\end{defn}

\begin{rem}
  With this definition, Proposition \ref{sttr} and Lemmas \ref{stun}, \ref{stpr}  are easily seen to
  hold in the general case too. 
\end{rem}

Thus we have the following
\begin{prop}
  The class of standard solutions is stable under the operations of direct
  sum, tensor product, projections and conjugation.
\end{prop}





The natural question is whether a choice of standard solution is available for
all 1-arrows. In the rest of this chapter we will give some partial answers.

In order to do so, we first recollect some results concerning Banach and
$C^*$-algebra bundles. As already mentioned, the notion of Banach bundle is more general than the more familiar notion of locally  trivial bundle, even in the case of bundles
with finite dimensional fibres.
The following example shows the necessity for considering  such a notion.

\begin{ex}

Consider the following tensor $C^*$-category (i.e.\  2-$C^*$-category with one
object): take $SU_n \times [0,1]$, the trivial group bundle.  $\Gamma (SU_{n}
\times [0,1])$ the group of continuous sections and $\mU:= \{ \xi \ \in \
\Gamma(SU_{n} \times [0,1]) \ s.t. \ \xi_{\omega} = 1_{n}$ (the identity of
$SU_{n}$)$, \forall \omega \in [\frac{1}{2},1] \}$, a closed subgroup. Consider the trivial bundle
$H \times [0,1]$,  where $H = \bC^n$, together with    the natural action 
of $\mU$ on it and  denote it  by $\rho$. Denote by $\rho^{\otimes n}$ the nth tensor product (fibre tensor product over the base space $[0,1]$) of the same bundle
with itself, with the  natural action of $\mU$.  For $n = 0$ let 
$\rho^0 = \iota =: \bC \times [0,1]$, the trivial line bundle  with the trivial action of $\mU$. 

The powers of $\rho$ induce  a tensor $C^*$-category, where $\Hom(\rho^n,\rho^m)$
are sections of intertwining operators, i.e.\  continuous sections 
$S \in (H^n,H^m) \times [0,1] \ $ such that $ \ S\rho^n(g) \xi = \rho^m(g) S
\xi$ ,$ \forall \xi \in \Gamma(H^n \times [0,1])$ , $ \forall g \in \mU$.   
It is easy to see that $\iota$ is the unit object in this tensor category.
The conjugate object $\bar{\rho}$ is the conjugate fibre $ \bar{H} \times
[0,1]$ with the conjugate  action of  $\mU$. A standard solution for the
conjugation equations is given by $ R \in \Hom(\iota, \bar{\rho} \rho):=
\sum_{i} \bar{e_i} \otimes e_{i}$ and  $ \bar{R} \in \Hom(\iota, \rho \bar \rho):=
\sum_{i} e_{i} \otimes \bar{e_i} $, where $e_i$ and $\bar e_i$ are the
constant sections given by the canonical basis for $H$ and $\bar H$ respectively.

Then $\End(\rho)    \ = \ \{S \in \Gamma(M_n \times [0,1]) $ such that$
S_{\omega} \in \bC 1_n \forall \omega \in [0,\frac{1}{2}] \}.$
Thus the fibres of the Banach bundle $\End(\rho)$ are of two types, $\bC$ for $\omega \in [0,\frac{1}{2}]$ and $M_n$ for $\omega \in (\frac{1}{2},1]$.
\end{ex}

We already know from the definition that $\forall \omega \in \Omega, \ \forall K^{\omega} \in \End(\rho)_{\omega},
\ \exists A
\in \End(\rho) $ such that $A_{|_{\omega}} = K^{\omega}$. In other words, each element
of a single fibre can be extended to a continuous section defined on the whole
base space. 
More can be said. The following lemma (whose proof can be found for example in \cite{Du}) will be useful in
the sequel. 

\begin{lem}
\label{est.imm}
  Let $F$ be a finite dimensional $C^*$-algebra, $\Xi$ a $C^*$-algebra bundle
  over a normal topological space $\Omega$ and $\Omega \times F$ be the
  product (trivial) $C^*$-algebra bundle with fibre $F$. Let
  $\Phi: A \times F \rightarrow \Xi_{|_{A}}$ be a $C^*$-algebra bundle embedding of
  the reduced bundles over a closed subset $A \subset \Omega$. Then there exists an
  open subset $U    \supset A$ and an embedding $\Psi: U \times F \rightarrow
  \Xi_{|_{U}}$ which extends $\Phi$.
\end{lem}

This enables us to prove the following

\begin{lem}
 Let $K^{\omega}$ be a positive invertible element in the finite dimensional algebra
 $\End(\rho)_{\omega}$. Then there exists an invertible $A \in \End(\rho)$
 and 
 an open neighbourhood 
$U_{\omega} $ of $\omega$ 
 such that $A_{|{\omega}} = K^{\omega}$ and $A_{|_{\omega'}} = 1_{\rho |_{\omega'}} \ \forall
 \omega' \notin U_{\omega}$. 
\end{lem}
{\noindent{\it Proof.}}  Take an open set $U_{\omega} \ni \omega$ and a $C^*$algebra bundle
embedding $\Psi: U_{\omega} \times \End(\rho)_{\omega} \rightarrow
\End(\rho)_{|_{U_{\omega}}}$ extending the identity bundle embedding.
Take a positive invertible section $H$ in the product bundle $U_{\omega} \times
\End(\rho)_{\omega}$ which extends $K^{\omega}$
(for example the constant section). Now take a second open set $W  \ni \omega$ such that $
U_{\omega} \supset \ovl{W}$ (we can do so, as the base space is normal) and a continuous complex-valued function $f$
defined on $U_{\omega}$ such $f = 0$ out of $W$ and $f(\omega) = 1$.
Then $\Psi(\exp ( f \ln H ))$ will be a continuous section with value $K^{\omega}$ in
$\omega$ and value $1_{\rho |_{\omega'}}$ for $\omega'$ out of $W$. Extending
it with the identity section on the rest of the base space $\Omega$ we get a
globally defined continuous section with the desired property. \hfill \halmos

We can now state a first result concerning standardness:
\begin{prop}
\label{weakstandardsolution}
  Let $\rho, \bar \rho$ be centrally balanced. For each $ \omega \in \Omega_A \ $ there exists a solution to the conjugation
  equations $R^{^{\omega}}, \bar R^{^{\omega}} \in \End(\rho)$ such that   $
  (1_{\rho} \otimes (R^{^{\omega} *} \circ
  (1_{\bar \rho} \otimes X) \circ R^{^{\omega}}) )_{|_{\omega}} =   
  (({\bar R}^{^{\omega} *} \circ (X
  \otimes 1_{ \rho}) \circ \bar R^{^{\omega}}) \otimes 1_{\rho})_{|_{\omega}}.$
\end{prop}

{\noindent{\it Proof.}}  Take a general solution $R, \bar R$, and suppose $\omega \in
S_l(\rho),$ otherwise the case being trivial. As Corollary
\ref{local.conjugation} points out, we can, mutatis mutandis, repeat the argument of the
simple unit case for the algebras $\End(\rho)_{\omega}$  and 
$\End(\bar \rho)_{\theta^{-1 *}_{\rho}(\omega)}$ and the local solutions
$R'_{|_{\omega}}, \bar R'_{|_{\theta^{-1 *}_{\rho}(\omega)}}$. By the remarks after 
Proposition \ref{discreteconiug} we see that there exists a positive
invertible $K^{\omega} \in \End(\rho)_{\omega}$ such that   
$$R^*_{|_{\omega}} \circ (1_{\bar \rho |_{\theta^{* -1}_{\rho}(\omega)}}
\otimes (K^{\omega -1 } \circ S_{|_{\omega}} \circ K^{\omega -1}) \circ
R_{|_{\omega}} =  $$ 
$$  \bar R^*_{|_{\theta^{-1 *}_{\rho}(\omega)}} \circ
((K^{\omega } \circ S_{|_{\omega}} \circ K^{\omega } \otimes 1_{\bar \rho
  |_{\theta^{-1 *}_{\rho} (\omega)}} )\circ \bar R_{|_{\theta^{-1 *}_{\rho}(\omega)}} $$     where  we have used the $\circ$ and $\otimes $
products of fibre elements introduced in the preceding section.
By the above lemma we can choose a positive invertible element $K \in
\End(\rho)$ such that $K_{|_{\omega}} = K^{\omega}$. 
Let $R^{\omega}:= (1_{\bar \rho} \otimes K^{-1}) \circ R$ and $\bar
R^{\omega}:= (K \otimes
1_{\bar \rho}) \circ \bar R$. \hfill \halmos

\begin{rem}
The above proposition may be viewed as a local version of standardness.  
It implies that for each $\omega \in \Omega,$  $R^{\omega *}_{|_{\omega}} \circ (1_{\bar \rho |_{\omega}} \otimes X_{|_{\omega}}) \circ
R_{|_{\omega}}$ is a uniquely defined (up to normalization) trace on $\End(\rho)_{\omega}.$
It would be tempting to use this trace as a definition for a standard
$\End(\iota_A)$-valued trace. 
Unfortunately in the general case the section $K^{\omega}$is not a priori
continuous, thus the above formula does not give a continuous trace, but
only an upper  semicontinuous one, as  it is the inferior limit of a family of
continuous functionals.
\end{rem}

Nevertheless, we have the following

\begin{prop}
\label{locallytrivialstandard}
  Suppose $\End(\rho)$ is a locally trivial bundle. Then a standard solution exists.
\end{prop}

{\noindent{\it Proof.}}  Suppose, for simplicity, $\rho, \bar \rho$ to be centrally balanced (if not,
decompose and consider each component separately, the associated
$\Hom(\rho_i,\rho_i)$ will still be locally trivial). If $\End(\rho)$ is locally trivial it has constant fibre, i.e.\  for
each point of the base space $\omega$, $\End(\rho)_{\omega}$ is isomorphic to
a finite dimensional algebra $F$. 
We can choose a finite atlas of local charts, i.e.\  maps $\Theta_{\alpha}:=
U_{\alpha} \times F \rightarrow \End(\rho)_{|_{U_{\alpha}}}$ which are local
isomorphisms of the trivial bundle $U_{\alpha} \times F $ onto the restriction
of $\End(\rho)$ over the open space $U_{\alpha } \subset \Omega$.
The sets  $U_{\alpha}$ form an open covering of $\Omega$. Where two maps overlap we have transition functions, i.e.\  when for example $U_{\alpha} \cap
U_{\beta} \neq \emptyset,$ we have unitary sections in $W_{\alpha,\beta} \in (U_{\alpha} \cap
U_{\beta}) \times F$  such that $\Theta_{\alpha |_{U_{\alpha} \cap U_{\beta}}}
= W_{\alpha,\beta} \circ \Theta_{\beta |_{U_{\alpha} \cap U_{\beta}}} \circ W^*_{\alpha,\beta}.$

Fix a faithful trace $tr$ on the algebra $F$. This defines a trace on each
local chart $\Theta_{\alpha}$. As the transition functions are unitary
sections, these local traces paste together into a continuous trace defined on
the whole bundle. Take a generic solution $R', \bar R'$ of the conjugation
equations. As shown in the remark following Proposition \ref{discreteconiug}  for each point
$\omega$ of the base space we have $ (R' \circ (1_{\bar \rho} \otimes S) \circ
R')_{|_{\omega}} = tr(D^{\omega}K^{\omega} S_{|_{\omega}})$ and $(\bar R' \circ (S
\otimes 1_{\bar \rho}) \circ \bar R')_{|_{\theta^{-1 *}_{\rho}(\omega)}} =
tr(D^{\omega} K^{\omega^{ -1}} S_{|_{\omega}})$ 
for a positive invertible $K^{\omega} $ and a positive invertible central
$D^{\omega}$ in
the fibre $\End(\rho)_{\omega}$ . Essentially we only have to prove that the  section  realized
by these $K^{\omega}$ is a continuous section. Then there will be a
corresponding element in $\End(\rho)$ fulfilling our requirements.

But as the left hand sides of the above equations are continuous, and as the
trace is continuous, this implies that both the sections $D^{\omega}K^{\omega}$
and $D^{\omega}K^{\omega^{-1}}$ are continuous sections in the locally trivial
bundle $\End(\rho)$. Thus $K^{\omega}$ is a continuous section, and a corresponding element
$K \in \End(\rho)$ exists.
$R:= (1_{\bar \rho} \otimes K^{-1})\circ R', \  \bar R:= (K \otimes 1_{\bar
  \rho}) \circ \bar R'$ will be a standard solution. \hfill \halmos

\section{Bundles of Hopf algebras}
\label{hopfbundles}

Finite irreducible subfactors of depth two are characterised by the action of finite
dimensional Hopf algebras (see, for example, \cite{Lo94}, \cite{Sz}). 
This situation corresponds, in the context of 2-$C^*$-categories with simple units, to an irreducible
$1$-arrow $\rho$
(i.e.\  $\End(\rho) = \mathbb{C}1_{\rho}$) generating a 2-$C^*$-category of
depth two (see below). 
The context can be generalised to the case of units with discrete and finite spectra (see \cite{SZ}), leading to the appearance of Weak-Hopf algebras. 

In the following we pursue the ``orthogonal'' direction, i.e. that of units with connected spectra,  obtaining a continuous bundle of finite dimensional Hopf algebras
in duality. We will follow the exposition given in \cite{Mu}, as it is closer
to our context.

We begin  by considering a centrally balanced $1$-arrow $B \xleftarrow{\rho} A$ such that  $\End(\rho) =
1_{\rho} \otimes \End(\iota_A) = \End(\iota_B) \otimes 1_\rho,$ with
$\Omega_A$ and  $\Omega_B$ connected homeomorphic spaces together with its conjugate $A
\xleftarrow{\bar \rho} B$. We will call such $\rho$ and $\bar \rho$ ``irreducible''.
We denote by $\mathcal{C}$ the 2-$C^*$-category generated by compositions of
$\rho$ and $\bar \rho$ (i.e.\  $ \rho, \ \rho \otimes \bar \rho, \ \rho \otimes \bar
\rho \otimes \rho \dots$) and their projections.

In the case of categories with simple units,  there is a well established
notion of finite depth: 
the 2-$C^*$-category generated  by $\rho$ and $\bar
\rho$  has finite depth $n$ if the
number of isomorphism classes of $1$-arrows is  finite (i.e.\  it's rational) and all of
them  appear as sub-$1$-arrows of the first $n$ products $ \rho, \rho
\otimes \bar \rho \dots$. 
For the sequence  of inclusions
$\End(\rho) \otimes 1_{\bar \rho} \subset \End(\rho \otimes \bar \rho) 
\dots$ this means that the corresponding principal part of the Bratteli
diagram is finite, with depth $n$.

In the general case of non-simple units we say  that $\rho$ has finite depth
$n$ if for each $\omega \in \Omega_A$ the associated  $\mathcal{C}^{\omega, A}$ has finite depth  and the maximum depth
among them is $n$.   
In particular we will be interested in the case of an irreducible $\rho$ of
depth two.

As $\rho$ is irreducible, it is possible to find standard solutions $R_\rho,
\bar R_\rho$. The same will be true for the
products $\rho \otimes \bar \rho, \ \rho \otimes \bar \rho \otimes \rho \dots$,
as the product solutions of standard solutions are still standard.
In order to simplify the notation, we set $$\mathcal{B}:= (\rho \otimes \bar \rho, \rho
\otimes \bar \rho), \ \mathcal{A}:= (\bar \rho \otimes \rho, \bar \rho \otimes \rho), \ \mathcal{D}
:=(\rho \otimes \bar \rho \otimes \rho, \rho \otimes \bar \rho \otimes
\rho).$$
$\mathcal{A}$ and $\mathcal{B}$ are in a natural way $C(\Omega_A)$
(respectively, $C(\Omega_B)$) bimodules, where left and right actions coincide
(as $\theta_{\rho \otimes \bar \rho}$ and $\theta_{\bar \rho \otimes \rho}$ are
    the identity isomorphisms).
As we  have standard solutions $R_{\bar \rho \otimes  \rho}, \bar R_{\bar \rho
  \otimes  \rho}$ for $\bar \rho \otimes  \rho$, the right and left $\End(\iota_A)$-valued inner products coincide as well
and give a faithful (non normalised) trace $Tr_{A}$ on $\mathcal{A}$:
$$ a \in \mathcal{A}, \ \ Tr_{A}(a):=   \langle a,1_{\bar \rho \otimes 
  \rho}\rangle_{{\End}(\iota_A)} = \ R^*_{\bar \rho \otimes  \rho} \circ (1_{\bar \rho \otimes 
  \rho} \otimes a ) \circ R_{\bar \rho \otimes  \rho}.$$
The same way we have a $\End(\iota_B) \cong C(\Omega_B)$-valued trace on $\mathcal{B}$:
$$ b \in \mathcal{B}, \ \ Tr_{B}(b):=   \langle b,1_{ \rho \otimes
  \bar \rho}\rangle_{{\End}(\iota_B)} = \ R^*_{ \rho \otimes \bar \rho} \circ
  (1_{ \rho \otimes
 \bar  \rho} \otimes b ) \circ R_{ \rho \otimes \bar  \rho}.$$

We will indicate for convenience $R^*_{\bar \rho \otimes \rho}
\circ R_{\bar \rho \otimes \rho}$ and $R^*_{\rho \otimes \bar \rho} \circ
R_{\rho \otimes \bar \rho}$ by $d_A$ and $d_B$, respectively.

  One defines a Fourier transform $ \mathcal{F}: \mathcal{A} \rightarrow \mathcal{B}$ as the linear
map defined by $$ \mathcal{F}(a):=  1_{ \rho \otimes \bar \rho} \otimes 
\bar R_{\rho}^* \circ (1_{ \rho} \otimes a \otimes
1_{\bar \rho}) \circ \bar R_{\rho} \otimes 1_{ \rho \otimes \bar \rho}.$$
and analogously $ \mathcal{\hat F}: B \rightarrow A$ as $$ \mathcal{\hat F}(b):=  1_{\bar \rho \otimes \rho} \otimes 
R_{\rho}^* \circ (1_{\bar \rho} \otimes b \otimes
1_{\rho}) \circ R_{\rho} \otimes 1_{\bar \rho \otimes \rho}.$$
The maps $\mathcal{S}:= \mathcal{\hat F} \circ \mathcal{F}: \mathcal{A}
\rightarrow \mathcal{A} $ 
and $\mathcal{\hat S}:= \mathcal{ F} \circ \mathcal{\hat F}: \mathcal{B} \rightarrow
\mathcal{B}$ are  the antipodes, and are antimultiplicative (this is an easy
consequence of the conjugation relations).
\begin{prop}
\label{theratrace}
  The Fourier transform preserves the inner product given by the traces
  $Tr_A$ and $Tr_B$ in the following sense:
$$ \forall a, a' \in \mathcal{A}, \ \ Tr_A(a'^* \circ a) =
\theta_{\bar \rho}(Tr_B(\mathcal{F}(a')^* \circ \mathcal{F}(a))),$$
$$ \forall b, b' \in \mathcal{B}, \ \ Tr_B(b'^* \circ b) =
\theta_{ \rho}(Tr_A(\mathcal{\hat F}(b')^* \circ \mathcal{\hat F}(b))).$$  
\end{prop}
{\noindent{\it Proof.}} 
We give only a sketch of the proof and skip the tedious  exposition of all the equalities. 

We have $Tr_A(\hat{\mathcal{F}}(b')^* \circ \hat{\mathcal{F}}(b)) = R^*_{\rho} \circ
(1_{\bar \rho} \otimes X) \circ R_{\rho}$, where the expression for $X \in
\End(\rho)$ is

$$(\bar R^*_{\rho} \otimes 1_\rho) \circ (1_{\rho} \otimes R^*_{\rho} \otimes 1_{\bar \rho}
\otimes 1_{\rho}) \circ (1_{\rho} \otimes 1_{\bar \rho} \otimes b'^* \otimes
1_{\rho}) \circ $$
$$  (1_{\rho \otimes \bar \rho \otimes \rho} \otimes  (R_{\rho} \circ  R^*_{\rho})) \circ
(1_{\rho} \otimes 1_{\bar \rho} \otimes b \otimes 1_{\rho}) \circ (1_{\rho}
\otimes R_{\rho} \otimes 1_{\bar \rho} \otimes 1_{\rho}) \circ (\bar R_{\rho}
\otimes 1_\rho)$$ 
As $R_{\rho},\bar R_{\rho}$ is standard, we have $R^*_{\rho} \circ (1_{\bar
  \rho} \otimes X) \circ R_{\rho} = \theta_{\rho}^{-1}( \bar R^*_{\rho} \circ
  (X \otimes 1_{\bar \rho}) \circ \bar R_{\rho}).$
One checks that $\bar R^*_{\rho} \circ
  (X \otimes 1_{\bar \rho}) \circ \bar R_{\rho} = Tr_{\mathcal{B}}(b'^* \circ b)$, which
  proves the second statement. The first statement is proved analogously. \hfill \halmos

We also have the following proposition, which is a   consequence of
standardness. We omit the proof, cf.\ e.g.\  \cite{Mu}.
\begin{prop}
    The following relations hold:$$\mathcal{S} \circ \mathcal{S} = id_{\mathcal{A}} \ ;
\    \mathcal{\hat S} \circ \mathcal{\hat S}  = id_{\mathcal{B}}$$ where by $id_{\mathcal{A}}$ and
    $id_{\mathcal{B}}$ we indicate  the identity endomorphisms of $\mathcal{A}$ and $\mathcal{B}$ respectively.   \end{prop}

We can define ``convolution'' products on $\mathcal{A}$ and $\mathcal{B}$ the following way: $$ a,
a' \in A, \ \ a \star a':= \mathcal{F}^{-1}(\mathcal{F}(a)\mathcal{F}(a')) \
; \ b,b' \in B, \ \ b \star b':= \mathcal{\hat F}^{-1}(\mathcal{\hat
  F}(b)\mathcal{\hat F}(b').$$

We restrict our attention for a moment to the case $\End(\iota_A) \cong \mathbb{C},$ 
$\End(\iota_B) \cong \mathbb{C}.$
Then $\mathcal{A}$ and $\mathcal{B}$ are finite dimensional algebras and we
are able to define a bilinear pairing between $\mathcal{A}$ and $\mathcal{B}$, i.e.\ 
a non-degenerate linear form $ \langle \cdot,\cdot\rangle: \mathcal{A} \otimes_{\mathbb{C}} \mathcal{B} \rightarrow
\mathbb{C}$, by $ \langle a,b\rangle:= d_{\rho}^{-1} Tr_A(a \mathcal{F}^{-1}(b))$, thus establishing a
duality (as linear spaces) between $\mathcal{A}$ and $\mathcal{B}$.
This duality enables us to define coproducts $\Delta: \mathcal{A} \rightarrow \mathcal{A} \otimes
\mathcal{A}, \ \hat \Delta: \mathcal{B} \rightarrow \mathcal{B} \otimes
\mathcal{B}$ by 
\begin{equation}
\label{co1}
 \langle \Delta(a), x \otimes y\rangle:=
 \langle a,xy\rangle, \
\ a \in \mathcal{A}, x,y \in \mathcal{B},
\end{equation} 
\begin{equation}
 \label{co2} 
 \langle a \otimes b, \hat \Delta(x)\rangle:=  \langle ab,x\rangle, \ \ a,b \in \mathcal{A}, x \in \mathcal{B}.
\end{equation}
Coassociativity of $\Delta$ and $\hat \Delta$ are implied by associativity of
the multiplications $m$ and $\hat m$ of $\mathcal{A}$ and $\mathcal{B}$ respectively.
We can also define counits 
\begin{equation}
  \label{counit}
\varepsilon(a):=  \langle a,1\rangle, \ \ a \in \mathcal{A}; \ \ \hat
\varepsilon(b):=  \langle 1,b\rangle, \ \ b \in \mathcal{B}.
\end{equation}
These operations  endow the algebras $\mathcal{A}$ and $\mathcal{B}$ with the structure of Hopf
algebras. We quote the following lemma and the
following proposition from \cite{Mu}: 

\begin{lem} (cf.\ \cite{Mu}, Lemma 6.18) 
  \label{lemma.mi}
The maps 
$$\Phi_1: \mathcal{A} \otimes_{\mathbb{C}} \mathcal{B} \rightarrow \mathcal{D}, \ \ a \otimes b
\mapsto 1_{\rho} \otimes a \circ b \otimes 1_{\rho},$$
$$\Phi_2: \mathcal{A} \otimes_{\mathbb{C}} \mathcal{B} \rightarrow \mathcal{D}, \ \ a \otimes b
\mapsto b \otimes 1_{\rho} \circ 1_{ \rho} \otimes a$$
are bijections.
\end{lem}

\begin{prop} (cf.\ \cite{Mu}, Proposition 6.19)
\label{proposition.mi}
  Let $\varepsilon, \hat \varepsilon, \Delta, \hat \Delta$ be defined as
  above, Then
  \begin{itemize}
  \item $\varepsilon, \hat \varepsilon$ are multiplicative,
  \item $\Delta, \hat \Delta$ are multiplicative,
\item $S, \hat S$ are coinverses, i.e.\  $m(S \otimes id) \Delta = m(id \otimes
  S)\Delta = \eta \varepsilon$, etc.  (where $\eta$ is the unit map $: 
  \mathbb{C} \ni c \mapsto c 1_{\mathcal{A}} \in \mathcal{A}$). 
\item $\mathcal{A}$ and $\mathcal{B}$ are finite dimensional Hopf algebras in duality, and $\mathcal{C}$ is
  the Weyl algebra (in the sense of \cite{Ni}) of $\mathcal{A} $.
  \end{itemize}
\end{prop}

The fact that the 2-$C^*$-category $\mathcal{C}$ is of depth two means that there is  only one isomorphism class of irreducible $1$-arrows
connecting $A$ to $B$, the one determined by $\rho$. The same holds for
$\bar \rho$.
This implies that for any sub$1$-arrow $X$ of $\bar \rho \otimes \rho$ we have
$ \rho \otimes X \ \cong  \oplus_{i = 1, \dots ,n} \ \rho$, the direct sum of $n$
copies of $\rho$. As the dimension is additive, we have $d_{\rho \otimes X  } =
n\ d_{\rho}$. But the dimension is multiplicative as well, i.e.\  $d_{\rho \otimes X   } = d_{X} d_{\rho}$.   This implies that the dimension $d_X$ of any sub
$1$-arrow $X$ of $\bar \rho \otimes \rho$ is an integer.
Such elements $\{ X$ sub$1$-arrow of $\bar \rho \otimes \rho \}$  form a
rational tensor category. In fact, as $\bar \rho \otimes \rho \otimes \bar
\rho \otimes \rho \cong \oplus_{1, \dots, d_{\bar \rho \otimes \rho}} \ \bar
\rho \otimes \rho$, all isomorphism classes of $1$-arrows in  $\mathcal{C}$ connecting $A$ to
$A$ appear in this set.
 One can construct a
faithful tensor functor from this category into the category of finite Hilbert
spaces assigning to each $X$ the complex Hilbert space of dimension $d_X$.
The natural transformations of this functor have the structure of a Hopf
algebra, and one can show that this is exactly the Hopf algebra $\mathcal{A}$
introduced above (see for example \cite{Mu}, Proposition 6.20).

Now let's return to the general case. We fix an $\omega \in \Omega_A$. Then the
$1$-arrows $\theta^{* -1}_{\rho}(\omega) \xleftarrow{\rho} \omega, \ \omega
\xleftarrow{\bar \rho}  \theta^{* -1}_{\rho}(\omega)$ in the category
$\mathcal{C}^{\omega, A}$ satisfy all the conditions of the above
propositions.
In particular $\End(\rho \otimes \bar \rho)_{\theta^{*
    -1}_{\rho}(\omega)} = \mathcal{B}_{\theta^{* -1}_{\rho}(\omega)}$ and $\End(\bar \rho \otimes
   \rho)_{\omega} = \mathcal{A}_{\omega}$ are finite
   dimensional Hopf algebras in duality.

Depth two of $\theta^{* -1}_{\rho}(\omega) \xleftarrow{\rho} \omega$ in the
category $\mathcal{C}^{\omega, A}$ implies that the associated dimension
$d_{\mathcal{A} |_{\omega}}$ has integer values and its square coincides with
the dimension (as a vector spaces) of the algebra $\mathcal{A}_{|_{\omega}}$. The same conclusion applies
to $\mathcal{B}_{_{\theta^{* -1}_{\rho}(\omega)}}$.
As this function is continuous  with respect to $\omega$, we conclude that $d_{\mathcal{A}}$ is a constant function on $\Omega_A$.
Thus the dimensions (as vector spaces) of the fibre algebras $\mathcal{B}_{\theta^{* -1}_{\rho}(\omega)}$
and $\mathcal{A}_{\omega}$ are constant respect to $\omega$.
Lemma \ref{est.imm} tells us that for each $\omega \in \Omega$ we can find a
neighbourhood $U$ and an algebraic embedding of $U \times \mathcal{A}_{|_{\omega}}$  into
$\mathcal{A}_{|_{U}}$. This embedding is actually surjective, as the fibre
algebras of $\mathcal{A}_{|_{U}}$ have all the same finite dimension. Thus all fibre algebras $\mathcal{A}_{\omega}$ are
isomorphic. i.e.\  $\mathcal{A}$ is a locally trivial bundle. The same
conclusion  applies to
$\mathcal{B}$.

Thus we have the following: $\mathcal{A}$ and $\mathcal{B}$ are locally trivial $C^*$-algebra
bundles over $\Omega_A$ and $\Omega_B$, with fibres isomorphic to finite
dimensional algebras $\mathcal{A}^0$ and $\mathcal{B}^0$ 
respectively.

We can view $\mathcal{B}$ as a $C(\Omega_A)$-valued Hilbert module by means of the
isomorphism $\theta_{\rho}$. Thus for $f \in C(\Omega_A)$, we have $
\theta_{\rho}(f)b = b\theta_{\rho}(f) \in \mathcal{B}, \ \forall b \in
\mathcal{B}$. The $C(\Omega_A)$-valued inner product is given by
$\theta_{\rho}^{-1}( \langle \cdot,\cdot\rangle_{{\End}(\iota_B)})$ and we can form the tensor product $\mathcal{A} \otimes_{C(\Omega_A)} \mathcal{B}$ where, for
example,$ a \otimes_{C(\Omega_A)}  \theta_{\rho}(f) b = 
a f
\otimes_{C(\Omega_A)} b$, for any $ a
\in A,\ b \in B, \ f \in C(\Omega_A).$ 
In other words  the usual tensor product of the fibre bundles
relative to $\mathcal{A}$ and $\mathcal{B}$. 
Analogously, we can define a $C(\Omega_A)$-valued linear non-degenerate
form $$  \langle \cdot,\cdot\rangle: \mathcal{A} \otimes_{C(\Omega_A)} \mathcal{B} \rightarrow
C(\Omega_A), \ \  \langle a,b\rangle:= d^{-1}_{\rho}
Tr_{\mathcal{A}}(a\mathcal{F}^{-1}(b)).$$ This form is well defined, as one
checks from the definition of $\mathcal{F}^{-1}$  that $f \otimes
\mathcal{F}^{-1}(b) = \mathcal{F}^{-1}(\theta_\rho(f)
\otimes b), \ \forall f \in C(\Omega_A), \ \forall b \in \mathcal{B}$  holds.

If we would like to think of $\mathcal{A}, \mathcal{B}$ as continuous bundles of
Hopf algebras, the natural candidate for a continuous comultiplication
$\Delta$ would be
defined by  $$ a \in \mathcal{A}, \ x, y \in \mathcal{B}, \ \ \Delta: a
\mapsto \Delta(a) \in \mathcal{A} \otimes_{C(\Omega_A)} \mathcal{A},$$
 such that $$ 
   \langle \Delta(a),x \otimes_{C(\Omega_A)} y\rangle \ = \  \langle a,xy\rangle.$$ 
As mentioned above, for each point $\omega \in \Omega_A$ we can consider the
category $\mathcal{C}^{\omega,A}$. The evaluation of the bilinear form $ \langle \cdot,\cdot\rangle$
in $\omega$ gives the pairing between the Hopf algebras $\mathcal{A}_{|_{\omega}}$
and $\mathcal{B}_{|_{\theta^{* -1}_{\rho}(\omega)}}$: $$b_{\theta^{* -1}_{\rho}(\omega)} \in
\mathcal{B}_{\theta^{* -1}_{\rho}(\omega)}, \ a_{\omega} \in \mathcal{A}_{\omega}, \ \
 \langle a_{\omega}, b_{\theta^{* -1}_{\rho}(\omega)}\rangle:= \  \langle a,b\rangle_{|_{\omega}}$$ 
where $ a \in A$ is such that $a_{|_{\omega}} =
a_{\omega}$, i.e.  $a$ is a continuous section in the
fibre bundle associated to the (bi-)module $\mathcal{A}$ whose value
corresponding to the base point $\omega$ is the
element $a_{\omega}$ of the finite dimensional fibre algebra $\mathcal{A}^0$, and the same
way $ b \in B$
      such that $b_{|_{\theta^{* -1}(\omega)}} = b_{\theta^{* -1}(\omega)}$.
Thus the above defined $\Delta(a)$ indeed defines a section in the fibre bundle
associated to the bimodule $\mathcal{A} \otimes_{C(\Omega_A)} \mathcal{A}$.
What we have to prove is that this section is continuous, i.e.\  it belongs to
$\mathcal{A} \otimes_{C(\Omega_A)} \mathcal{A}.$ More precisely, we will show the following

\begin{lem}
Let $a \in \mathcal{A}$. Then $\forall \epsilon > 0$ there exist a finite set
of $a^k_1, a^k_2 \in \mathcal{A}$ such that $$| (\sum_k  \langle a^k_1,b\rangle \langle a^k_2,c\rangle -
 \langle a,bc\rangle)_{|_{\omega}} |  <  \epsilon,  \ \forall \omega  \in \Omega_A, \ \forall b,c
\in \mathcal{B}, \  \| b \| \leq 1, \| c \| \leq 1.$$  
\end{lem}

{\noindent{\it Proof.}}  For our convenience we will choose for  $\mathcal{A}$ and
$\mathcal{B}$ the norms given by the inner product, i.e.\  the $C(\Omega_A)$
(resp. $C(\Omega_B)$)-valued traces $Tr_{\mathcal{A}}$
(resp. $Tr_{\mathcal{B}}$). 
We choose locally trivial algebraic maps $\Phi_i: \mathcal{A}_{|_{U'_i}} \rightarrow U'_i
\times \mathcal{A}^0$, $\Psi_i: \mathcal{B}_{|_{V'_i}} \rightarrow V'_i \times \mathcal{B}^0$, where $ \{ U'_i \}$
and $\{ V'_i \}$ are open coverings of $\Omega_A$ and $\Omega_B$. Without loss of
generality, we suppose that $V'_i = \theta^{* -1}_{\rho}(U'_i)$.
For each $x \in \mathcal{A}^0$ we can consider the corresponding constant section in $U'_i
\times \mathcal{A}^0$, which we will indicate with the same symbol. Then
$\Phi_i (x)^{-1}$ will be an element in $\mathcal{A}_{|_{U'_i}}$.
The same way, for $y \in \mathcal{B}^0$, $\Psi^{-1}_i (y)$ will be an
element of $\mathcal{B}_{|_{V'_i}}.$

Then we have the following linear form on $\mathcal{A}^0
\otimes_{\mathbb{C}} \mathcal{B}^0$ defined as 
$$  \langle \Phi^{-1}_i(x),\Psi^{-1}_i(y)\rangle_{|_{\omega}}, \ \ \omega \in U'_i.$$
This linear form depends continuously on $\omega$. 
We define 
 $$\Delta^{\omega}_i: \mathcal{A}^{0} \rightarrow \mathcal{A}^{0} \otimes
\mathcal{A}^{0}, \ \ x \in \mathcal{A}^0, \ \Delta^{\omega}_i(x):= \sum_l x^l_1
\otimes x^l_2 $$ such that for any pair $y,z \in \mathcal{B}^0$ $$
(\sum_l  \langle \Phi^{-1}_i(x^l_1),\Psi^{-1}_i(y)\rangle  \langle \Phi^{-1}_i(x^l_2),\Psi^{-1}_i(z)\rangle)_{|_{\omega}}
=  \langle \Phi^{-1}_i(x),\Psi^{-1}_i(yz)\rangle_{|_{\omega}}.$$
As $\mathcal{A}^0$ is finite dimensional, all norms give equivalent
topologies. For our convenience we will choose the following as norm:
$$ \| x \|:= sup_{\omega \in U'_i} (Tr_{\mathcal{A}}
\Phi^{-1}_i(x^*x))^{\frac{1}{2}}_{|_{\omega}}, \ \ x \in \mathcal{A}^0, \
  \forall U'_i.$$
We may suppose that we have chosen the sets $U'_i$ such that 
\begin{equation}
\label{eq:eq1} 
\|   \Phi_i(a)_{|_{\omega}} - \Phi_i(a)_{|_{\omega'}} \| \  <  \epsilon, \ \ \omega,
\omega' \in U'_i 
\end{equation}
as $\Phi_i(a)$ is a continuous section of $ U'_i \times \mathcal{A}^0 $. 
We choose an $\omega_i$ in $U'_i$. Then, as the map $\Delta^{\omega}_i$ is continuous respect to $\omega$, we may also  suppose that we have
chosen a set  $W_i \ni \omega_i$ such that 
\begin{equation}
  \label{eq:eq2}
  \|
\Delta^{\omega}_i(\Phi_i(a)_{|_{\omega_i}}) -
\Delta^{\omega'}_i(\Phi_i(a)_{|_{\omega_i}}) \|  <  \epsilon, \ \omega,\omega' \in
W_i,
\end{equation}
where $\Delta^{\omega}_i(\Phi_i(a)_{|_{\omega_i}}), \
\Delta^{\omega'}_i(\Phi_i(a)_{|_{\omega_i}})$ are elements of $\mathcal{A}^0
\otimes \mathcal{A}^0$  with the norm specified just above.

Letting $\omega_i$ vary arbitrarily in $U'_i$, the $W_i$ form an open cover
of $\Omega_A$. We choose a finite refinement of the two open covers of
$\Omega_A$ given by the $\{ U'_i \} $ and $\{ W_i \}$  such that both
(\ref{eq:eq1}) and (\ref{eq:eq2}) hold. We will denote with little abuse of
notation these sets by $\{ U_i \}$ and by $\{ \Phi_i \}$ the relative local charts.

\begin{rem}
We note that for any $b \in \mathcal{B}, \ \| b \| = 
\| Tr_{\mathcal{B}}(b^*b)^{\frac{1}{2}} \|^{C(\Omega_B)}  <  1$
and $a \in \mathcal{A}, \ \| a \| = \| Tr_{\mathcal{A}}(a^*a)^{\frac{1}{2}} \|^{C(\Omega_A)}  <  1$ one has $|
 \langle a,b\rangle_{|_{\omega}} |  <  \| d_{\rho}^{} \|$ by the definition of the bilinear form
and Proposition \ref{theratrace}.
\end{rem}

Now let's consider the difference

$$|  \langle \Delta(a) - (\Phi^{-1}_i \otimes \Phi^{-1}_i) \circ
\Delta^{\omega_i}_i(\Phi_i(a)_{|_{\omega_i}}),b \otimes c\rangle_{|_{\omega}} | \ \ = $$$$|  \langle \Delta(a)   -  \langle  (\Phi^{-1}_i \otimes \Phi^{-1}_i) \circ
\Delta^{\omega}_i(\Phi_i(a)_{|_{\omega_i}}),b \otimes c\rangle_{|_{\omega}}  $$
$$ +  \langle  (\Phi^{-1}_i \otimes \Phi^{-1}_i) \circ
\Delta^{\omega}_i(\Phi_i(a)_{|_{\omega_i}}),b \otimes c\rangle_{|_{\omega}}   -  \langle  (\Phi^{-1}_i \otimes \Phi^{-1}_i) \circ
\Delta^{\omega_i}_i(\Phi_i(a)_{|_{\omega_i}}),b \otimes c\rangle_{|_{\omega}} |  
$$
$$\leq \ |  \langle  (\Delta(a),b \otimes c\rangle_{|_{\omega}}     -  \langle  (\Phi^{-1}_i \otimes \Phi^{-1}_i) \circ
\Delta^{\omega}_i(\Phi_i(a)_{|_{\omega_i}}),b \otimes c\rangle_{|_{\omega}}  |$$
$$+ |  \langle  (\Phi^{-1}_i \otimes \Phi^{-1}_i) \circ
\Delta^{\omega}_i(\Phi_i(a)_{|_{\omega_i}}),b \otimes c\rangle_{|_{\omega}}   -  \langle  (\Phi^{-1}_i \otimes \Phi^{-1}_i) \circ
\Delta^{\omega_i}_i(\Phi_i(a)_{|_{\omega_i}}),b \otimes c\rangle_{|_{\omega}} | $$
evaluated in $\omega \in U_i.$ Notice that $ \langle \Delta(a),b \otimes c\rangle_{|_{\omega}} =  \langle a,b
c\rangle_{|_{\omega}}$ (by definition of $\Delta$) $ =  \langle \Phi_i^{
  -1}(\Phi_i(a)_{|_{\omega}}),bc\rangle_{|_{\omega}}$. Also $ \langle (\Phi^{-1}_i \otimes
\Phi^{-1}_i) \circ \Delta^{\omega}_i(\Phi_i(a)_{|_{\omega_i}}),b\otimes c\rangle_{\omega} =
 \langle \Phi^{-1}_i (\Phi_i(a)_{|_{\omega_i}}),bc\rangle_{\omega}$. It follows from
(\ref{eq:eq1}) and the remark above that the first summand satisfies $|  \langle \Phi^{
  -1}_i(\Phi_i(a)_{|_{\omega}}),bc\rangle_{\omega} -  \langle \Phi^{ 
  -1}_i(\Phi_i(a)_{\omega_i}),bc\rangle_{\omega} |  <  \| d_{\rho}^{} \| \epsilon$.
Analogously, from continuity of $\Delta^{\omega}_i$ with respect to $\omega$ and  (\ref{eq:eq2})
the second summand is $ <  \| d_{\rho} \|^{2} \epsilon$. 

Thus we have proven that $$|  \langle \Delta(a) - (\Phi^{-1}_i \otimes \Phi^{-1}_i) \circ
\Delta^{\omega_i}_i(\Phi_i(a)_{|_{\omega_i}}),b \otimes c\rangle_{|_{\omega}} | \  <
(\| d_{\rho}^{-1} \| + \| d_{\rho}^{-1} \|^2) \epsilon, \ \forall \omega \in U_i. $$

Now take a partition of unity $ \{ f_i \} $ subordinate to the open covering of
$\Omega_A$ realized by the $U_i$.  Then  $  \sum_i  f_i \ (\Phi^{-1}_i \otimes \Phi^{-1}_i) \circ
\Delta^{\omega_i}_i(\Phi_i(a)_{|_{\omega_i}}) $ is a sum  $\sum_k
a^k_1 \otimes_{C(\Omega_A)}  a^k_2$ of elements of $\mathcal{A}
\otimes_{C(\Omega_A)} \mathcal{A}$
   which satisfies the claim of
the
lemma  (modulo the constant  $(\| d_{\rho}^{}\| + \| d_{\rho}^{}\|^2) )$. \hfill \halmos

An analogous result holds for $\mathcal{B}$. Thus we see that (\ref{co1}),
(\ref{co2}), (\ref{counit}) make sense even in the case 
$\End(\iota_A) \neq \mathbb{C}$. We can think of (\ref{counit}) as a
$\End(\iota_A)$-valued counit (resp. $\End(\iota_B)$-valued counit),
or as continuous
sections of counits for the fibre algebras. Continuity is obvious as in the
definition only continuous objects are involved. The only non  trivial part is
to check  continuity of the maps $\Delta: \ \mathcal{A} \rightarrow 
\mathcal{A}$ and $\hat \Delta: \ \mathcal{B} \otimes_{C(\Omega_B)} \mathcal{B}$,
which has been established by the previous lemma.
We have thus the following analogue of Proposition \ref{proposition.mi}:
\begin{prop}
 
  Let $\mathcal{A},\mathcal{B,}\varepsilon, \hat \varepsilon, \Delta, \hat \Delta$ be defined as
  above, then
  \begin{itemize}
\item $\mathcal{A}$ and $\mathcal{B}$ are locally trivial bundles of Hopf
  algebras, with fibres the finite dimensional algebra $\mathcal{A}^0$ and
  $\mathcal{B}^0$ respectively,  
\item $\varepsilon, \hat \varepsilon$ are multiplicative,
  \item $\Delta, \hat \Delta$ are multiplicative,
\item $S, \hat S$ are coinverses, i.e.\  $m(S \otimes id) \Delta = m(id \otimes
  S)\Delta = \eta \varepsilon$, etc. 
\item for every $\omega \in \Omega_A, \ \mathcal{A}_{|_{\omega}}$ and $\mathcal{B}_{|_{\theta^{* -1}_{\rho}(\omega)}}$ are finite dimensional Hopf algebras in duality, and $\mathcal{C}_{|_{\omega}}$ is
  the Weyl algebra of $\mathcal{A}_{|_{\omega}}$.
 \end{itemize}
\end{prop}

{\noindent{\it Proof.}}  All algebraic properties follow by applying Proposition
\ref{proposition.mi} to each fibre algebra. Continuity of the maps follows from the above
remarks and the above lemma. \hfill \halmos

\begin{rem}
 A  notion of $C(X)$ Hopf algebra with a continuous field of
 coproducts over the topological space $X$ was introduced in \cite{Bla96}. 
This is a more general definition dealing with fields (i.e.\  bundles) of
 infinite and  not necessarily unital $C^*$-algebras over a locally compact
 space $X$. It is fairly easy to see that our example fits this definition as well. 

\end{rem}

\section{Frobenius algebras and Q-systems}
\label{Frobenius}
In this section we would like to make some remarks concerning Frobenius
algebras and $Q$-systems. The notion of a $Q$-systems was introduced in the
context of 
von Neumann algebras in \cite{Lo94} (where, in some sense,  it plays the role of
a formalisation of a sub-factor)  and subsequently for general tensor
$C^*$-categories the notion of  abstract
$Q$-system was introduced in \cite{LR97}. 
The notion of Frobenius algebra in a tensor category is more general, but we will show that in the $C^*$-category  case the two notions coincide. We will
follow the exposition given in \cite{Mu} where not only the $C^*$-case is
treated but more general categories are considered and a general
correspondence between Frobenius algebras and pairs of conjugate $1-$arrows in $2-$categories is studied.    

\begin{defn}
  \label{frobeniusalgebra}
(cf.\ \cite{Mu}, Definition 3.1.)
Let $\mathcal{A}$ be a strict tensor category. A Frobenius algebra in
$\mathcal{A}$ is a quintuple $(\lambda,V,V',W,W')$, where $\lambda$ is an object in
$\mathcal{A}$ and $V: \iota \rightarrow \lambda, \ V': \lambda \rightarrow \iota, \ W: \lambda
\rightarrow \lambda^2, \ W': \lambda^2 \rightarrow \lambda$ are morphisms satisfying the
following conditions:

\begin{equation}
  W \otimes 1_{\lambda} \circ W \ = \ 1_{\lambda} \otimes W \circ W
\end{equation}
\begin{equation}
  W' \circ W' \otimes 1_\lambda \ = \ W' \circ 1_\lambda \otimes W'
\end{equation}

\begin{equation}
  V' \otimes 1_\lambda \circ W \ = 1_{\lambda} = \ 1_\lambda \otimes V' \circ W
\end{equation}
\begin{equation}
  W' \circ V \otimes 1_\lambda \ = \ 1_\lambda \ = \ W' \circ 1_\lambda \otimes V
\end{equation}
\begin{equation}
  W' \otimes 1_\lambda \circ 1_\lambda \otimes W \ = \ W \circ W' \ = \
  1_\lambda \otimes W'
  \circ W \otimes 1_\lambda.
\end{equation}
\end{defn}

\begin{defn}
  (cf.\ \cite{Mu}, Definition 3.3) Two Frobenius algebras $(\lambda,V,V',W,W'),$ $ \
  (\tilde \lambda, \tilde V, \tilde V ', \tilde W, \tilde W')$ in the strict tensor
  category $\mathcal{A}$ are isomorphic if there is an isomorphism $S: \lambda
  \rightarrow \tilde \lambda$ such that  $$S \circ V = \tilde V, \ \ V' = \tilde V'
  \circ S, \ \ S \otimes S \circ W = \tilde W \circ S, \ \ S \circ W' = \tilde
  W' \circ S \otimes S. $$
\end{defn}

We recall the definition of conjugation in the case of a (not necessarily $C^*$)
 2-category $\mathcal{C}$ (the term ``duality'' is often used instead of ``conjugation''):

\begin{defn}
  A $2$-category $\mathcal{C}$ is said to have left ( right) duals if for
  every 1-arrow $\rho:B \leftarrow A \in \mathcal{C}$ there is $\bar{ } \rho
 :A \leftarrow B$
  ($\rho \bar{ }: A \leftarrow B$ )  together
  with 2-arrows $e_{\rho} \in \Hom(\iota, \rho \otimes \bar{ } \rho), \
  d_{\rho} \in \Hom ( \bar{}
  \rho \otimes \rho,  \iota) \ \ (\varepsilon_\rho \in \Hom(\iota, \rho \bar{}
  \otimes \rho), \ \eta_\rho \in \Hom (\rho \otimes \rho \bar{}, \iota))$
satisfying: $$1_\rho \otimes d_\rho \circ e_\rho \otimes 1_\rho = 1_\rho \ , \ \
  d_\rho \otimes 1_{\bar{} \rho} \circ 1_{\bar{} \rho} \otimes
  e_{\rho}=1_{\bar{} \rho}$$
$$\eta_\rho \otimes 1_{\rho \bar{}} \circ 1_{\rho \bar{}} \otimes
\varepsilon_{\rho}=1_{\rho \bar{}} \ , \ \ 1_\rho \otimes \eta_\rho \circ \varepsilon_\rho \otimes 1_\rho = 1_\rho$$

If $\bar{} \rho =
 \rho \bar{}$, $\bar{} \rho$ is said to be a    two-sided dual, and we
 indicate it by $\bar \rho$. 

\end{defn}
 We will assume in the sequel duals to be two-sided.
Duals
are automatically two sided in a $*$-category. It is easy to see that the
above definition of duality (i.e.\  conjugation) reduces to the one already
introduced for the $C^*$-case.

\begin{lem}
The object $\lambda$ of a Frobenius algebra is self-conjugate.
  
\end{lem}
{\noindent{\it Proof.}}  Set $e_\lambda:= W \circ V$ , $d_\lambda:= V' \circ W'$. It is easy
to see that they satisfy the claimed relations.\hfill \halmos

The following is an important example:
\begin{lem}
  (cf.\ \cite{Mu}, Lemma 3.4) Let $ \rho: B \leftarrow A $ be a $1$-arrow in a
  $2$-category $\mathcal{E}$ and let $\bar \rho: A \leftarrow B$ be a two sided
  dual with duality $2$-morphisms $d_{\rho}, e_{\rho}, \varepsilon_{\rho}, \eta_{\rho}$. Positing
  $\lambda = \bar \rho \otimes  \rho: A \leftarrow A$ there are $ V, V', W,
  W'$ such that $(\lambda,V,V',W,W')$ is a Frobenius algebra in the tensor category
  $\mathcal{A} = HOM_{\mathcal{E}}(A,A)$. 
\end{lem}
{\noindent{\it Proof.}}  It suffices to choose $$ V:= \varepsilon_\rho, \ V':= d_\rho, \ W:= 1_{\bar
  \rho} \otimes e_\rho \otimes 1_\rho, \ W':= 1_{\bar \rho} \otimes \eta_\rho
  \otimes 1_\rho.$$ It is not difficult to check that the Frobenius algebra
  relations hold. \hfill \halmos 

 In the sequel we will prove explicitly a similar result in
  the $C^*$-case.  
The following propositions show to which extent a generic Frobenius algebra
can be realized as a couple of conjugate 1-arrows in a 2-category as in the
example above.

\begin{defn}
  An almost-2-category is defined as a 2-category except that we do not
  require the existence of a unit 1-arrow $\iota_{\sigma}$ for every object $\sigma$.
\end{defn}

\begin{prop}
\label{3.8}
  (cf.\ \cite{Mu}, Proposition 3.8) Let $\mathcal{A}$ be a strict tensor
  category and $\lambda = (\lambda,V,V',W,W')$ a Frobenius algebra in
  $\mathcal{A}$. Then there is  an almost-2-category $\mathcal{E}_0$
  satisfying:
  \begin{itemize}
  \item Obj $\mathcal{E}_0 = \ \{ A,B \}$.
  \item There is an isomorphism $I: \mathcal{A} \rightarrow
  HOM_{\mathcal{E}_0}(A,A)$ of tensor categories.
  \item There are 1-arrows $\rho: B \leftarrow A$ and $\bar \rho: A
  \leftarrow B$ such that $ \bar \rho \otimes \rho = I(\lambda)$.
  \end{itemize}

If $\mathcal{A}$ is $\mathbb{K}$-linear then so is $\mathcal{E}_0$. Isomorphic
Frobenius algebras give rise to isomorphic almost-2-categories.
\end{prop}

\begin{thm}
 \label{conj2cat} (cf.\ \cite{Mu}, Theorem 3.11) Let $\mathcal{A}$ be a strict tensor category
  and $\lambda = (\lambda,V,V',W,W')$ a Frobenius algebra in $\mathcal{A}$. Assume that
  one of the following conditions is satisfied:
  \begin{itemize}
  \item $W' \circ W = 1_\lambda.$
  \item $A$ is $\End(\iota)$-linear and $$W' \circ W = z_1 \otimes 1_\lambda,$$
  where $z_1$ is an invertible element of the commutative monoid
  $\End(\iota)$.
\end{itemize}
Then the completion $\mathcal{E} =  \mathcal{E}^P$ of the $\mathcal{E}_0$
defined in Proposition (\ref{3.8}) is  a bicategory such that
\begin{itemize}
\item $Obj \mathcal{E} = \{ A,B \}$.
\item There is a fully faithful tensor functor $I: \mathcal{A} \rightarrow
  HOM_{\mathcal{E}}(A,A)$ such that for every $Y \in HOM_{\mathcal{E}}(A,A)$
  there is $X \in \mathcal{A}$ such that $Y$ is a retract (i.e.\  sub-$1$-arrow) of $I(X)$. 
\item There are 1-arrows $\rho: B \leftarrow A$ and $ \bar \rho: A \leftarrow B$
  and     2-arrows $$ e_\rho:\iota_B \rightarrow  \rho \otimes
  \bar \rho, \
  \varepsilon_\rho: \iota_A \rightarrow  \bar \rho \otimes  \rho, \ d_\rho:
  \bar \rho \otimes \rho \rightarrow \iota_A, \
  \eta_\rho:  \rho \otimes  \bar \rho \rightarrow \iota_B$$ satisfying the
  conjugation (i.e.\  duality)  relations.
\item We have the identity $$I(\lambda,V,V',W,W') = (\bar \rho \otimes \rho , e_\rho,\eta_\rho,1_\rho
  \otimes \varepsilon_\rho \otimes 1_{\bar \rho}, 1_\rho \otimes d_\rho
  \otimes 1_{\bar \rho}) $$
of Frobenius algebras in $HOM_{\mathcal{E}}(A,A)$.
\item If $\mathcal{A}$ is a preadditive category, then $\mathcal{E}$ is a
  preadditive 2-category.
\item If $\mathcal{A}$ has direct sums then $\mathcal{E}$ has direct sums of 1-arrows.
\end{itemize}
Isomorphic Frobenius algebras $\lambda, \tilde \lambda$ give rise to isomorphic
bicategories $\mathcal{E}, \tilde{\mathcal{E}}$.

\end{thm}

\begin{rem}
As we have seen, a generic Frobenius algebra can be realized as the product of
a couple of
1-arrows in a bicategory. In order for these 1-arrows to be conjugate, in Theorem
\ref{conj2cat} additional hypotheses were required.     
With further requirements one can prove the universality of this construction.
\end{rem}

\begin{defn}
  (cf.\ \cite{Mu}, Definition 3.13) Let $\mathcal{A}$ be an
  $\End(\iota)$-linear category. A Frobenius algebra $(\lambda,V,V',W,W')$ in
  $\mathcal{A}$ is ``strongly separable'' iff $$W' \circ W = z_1 \otimes 1_\rho,$$
  $$V' \circ V = z_2,$$ where $z_1, z_2 \in \End(\iota)$
  are invertible. $(\lambda, V, V', W, W')$ is said to be normalised if $z_1 = z_2$.
\end{defn}

\begin{thm}
  (cf.\ \cite{Mu}, Theorem 3.17) Let $\mathcal{A}$ be $\End(\iota)$-linear
  and $(\lambda,V,V',W,W')$ a strongly separable Frobenius algebra in
  $\mathcal{A}$. Let $\mathcal{E}$ be as constructed in Theorem \ref{conj2cat}
and
  $\tilde{\mathcal{E}}$ be any bicategory such that:
  \begin{itemize}
  \item Obj $\tilde{\mathcal{E}} = \{ A,B \}.$
  \item Idempotent 2-arrows in $\tilde{\mathcal{E}}$ split.
  \item There is a fully faithful tensor functor $\tilde I: \mathcal{A}
  \rightarrow HOM_{\tilde{\mathcal{E}}}(A,A)$ such that every object of
  $HOM_{\tilde{\mathcal{E}}}(A,A)$ is a retract of $\tilde I(X)$ for some $X
  \in \mathcal{A}$.
  \item There are mutually two-sided dual 1-arrows $\tilde \rho: B \leftarrow
  A, \ \tilde{\bar \rho}: A \leftarrow B $ and an isomorphism $\tilde S: I(\lambda)
  \rightarrow \tilde {\bar \rho} \otimes \tilde{\rho}$ between the Frobenius algebras
  $I(\lambda,V,V',W,W')$ and $ (\tilde {\bar \rho} \otimes \tilde{  \rho}, \tilde{e}_{\tilde \rho},
  \dots)$ in $HOM_{\tilde{\mathcal{E}}}(A,A)$.
  \end{itemize}

Then there is an equivalence $E: \mathcal{E} \rightarrow \tilde{\mathcal{E}}
$ of bicategories such that there is a tensor isomorphism between the tensor
functors $\tilde I$ and $(E |_ {HOM_{\tilde{\mathcal{E}}}(A,A)}) \circ I.$
\end{thm}

We recall the notion of a $Q$-system:

\begin{defn} \label{defQ}
Let $\mathcal{A}$ be a tensor $*-$category. A $Q$-system in
$\mathcal{A}$ is a triple $(\lambda, T , S )$ where $\lambda $ is an object
in $\mathcal{A}$ and $T \in \Hom(\iota , \lambda)$ , $S \in \Hom(\lambda, \lambda^2)$
are arrows satisfying the following relations:

\begin{equation}
  \label{eq:Q1}
 T^* \otimes 1_{\lambda} \circ S = 1_{\lambda} = 1_{\lambda} \otimes T^*
  \circ S
\end{equation}

\begin{equation}
  \label{eq:Q2}
 S^* \circ S = 1_\lambda 
\end{equation}
 
\begin{equation}
  \label{eq:Q3}
   S \otimes 1_{\lambda} \circ S = 1_{\lambda} \otimes S \circ S  
\end{equation}

\begin{equation}
\label{eq:Q4}
  S^* \otimes 1_{\lambda} \circ 1_{\lambda} \otimes S = S \circ S =
  1_{\lambda} \otimes S^* \circ S \otimes 1_{\lambda}.  
\end{equation}

\end{defn}

\begin{rem}
  It follows from the definition that $\lambda$ is self conjugate and
  $S_l(\lambda) = S_r(\lambda)$.
\end{rem}

\begin{rem}
The relations above are  the same as in the original definition of $Q$-system (cf.\ \cite{LR97}) given for the case of a tensor $C^*$-category with simple unit. $\End(\iota)$-linearity was not explicitly assumed, as it holds trivially when $\End(\iota) \cong \mathbb{C}$. 
Here we have only generalised the context. In   case the category $\mathcal{A}$ is $\End(\iota)$-linear, we will talk about an $\End(\iota)$-linear $Q$-system.  
 
  In \cite{Mu}  a $Q$-system was defined as a strongly separable
  Frobenius algebra $(\lambda, T, T^*, S, S^*)$  in a tensor $*$-category (thus assuming $\End(\iota)$-linearity). 
In the $C^*$ case this latter  definition is almost  equivalent with our definition of $\End(\iota)$-linear $Q$-system: given a strongly separable
  Frobenius algebra $(\lambda,T,T^*,S,S^*)$ it is sufficient to renormalise
  $T$ and $S$ by the invertible $z_1$ in order to turn $S$ into an isometry.
On the other hand, given a $Q$-system as defined above, $z_1:= T^* \circ T$ and
  $z_2:= T^* \circ S^* \circ S \circ T$ are positive and invertible on
  $S_l(\lambda)$, as the following lemma shows:
\end{rem}

\begin{lem}
 Let $(\lambda,T,S)$ be a $Q$-system in a tensor $C^*$-category. Then $T^* \circ S^*
  \circ S \circ T$ and $T^* \circ T$ are positive  elements of
  $\End(\iota)$ invertible on $S_l(\lambda) = S_r(\lambda).$
\end{lem}
$Proof. $ $S \circ T$, $T^* \circ S^*$ satisfy the conjugation relations for $\lambda$ and this
implies that $T^* \circ S^* \circ S \circ T$ is invertible on $S_l(\lambda)$
by Corollary \ref{closets}. The inequality $$
\| S^* \circ S \| \ \ T^* \circ T \ \geq T^* \circ S^* \circ S \circ T$$ implies that
$T^* \circ T$ is invertible on $S_l(\lambda)$ as well. \hfill \halmos

Thus, bearing in mind the observations at the end of
section \ref{Preliminaries}, we see that it is always possible to
embed an $\End(\iota)$-linear Q-system $(\lambda,T,S)$ (in fact, the whole tensor $C^*$-category
generated by its tensor powers) into a tensor $C^*$-category such that
$S_l(\lambda) = \Omega$, the topological space associated to
the centre of the new tensor category. In this category $T^* \circ S^*
\circ S \circ T$ and $T^* \circ T$ are invertible elements in
$\End(\iota)$, i.e.\  $(\lambda,T,S)$ is a strongly separable Frobenius
algebra.

\begin{rem}

We will in the sequel indicate by $(T^* \circ
T)^{-1}_{| _{S_l(\lambda)}}$ the positive element in $\End(\iota)
\otimes E_{| _{S_l(\lambda)}}$ such that $(T^* \circ
T)^{-1}_{| _{S_l(\lambda)}} \circ T^* \circ T = E_{| _{S_l(\lambda)}}.$ Thus, in particular, $((T^*
\circ T)^{-1}_{| _{S_l(\lambda)}} \circ T^* \circ T) \otimes 1_\lambda = 1_\lambda.$  

\end{rem}

 We will see now that the  four relations in Definition \ref{defQ} are not independent. Having assumed (\ref{eq:Q1}) to
 hold, we can choose any pair of  (\ref{eq:Q2}), (\ref{eq:Q3}),
(\ref{eq:Q4}) and the  remaining relation will follow (up to
 isomorphism). That (\ref{eq:Q1}, \ref{eq:Q2}, \ref{eq:Q3}) imply
(\ref{eq:Q4}) was proven in  \cite{LR97}. That (\ref{eq:Q1}, \ref{eq:Q2},
 \ref{eq:Q4}) imply (\ref{eq:Q3}) was proven in \cite{IK}.  
It is not necessary to suppose $\End(\iota) \cong \mathbb{C}$ and
 $\End(\iota)$-linearity plays no role in the proof.

 \begin{prop}
\label{qequivalenza}
  Let $\lambda, S, T$ be as in Definition \ref{defQ} in a (not necessarily
  $\End(\iota)$-linear) tensor $C^*$-category $\mathcal{A}$. 
Assume (\ref{eq:Q1}) to hold. Then we have the following implications
\begin{itemize}
\item $ (\ref{eq:Q2})+(\ref{eq:Q3}) \ \Rightarrow (\ref{eq:Q4})$
\item $ (\ref{eq:Q2})+(\ref{eq:Q4}) \ \Rightarrow (\ref{eq:Q3})$
\item $ (\ref{eq:Q3})+(\ref{eq:Q4}) \ \Rightarrow (\ref{eq:Q2})$ for $S', \
  T'$ isomorphic to $S, \ T$.
\end{itemize}
 \end{prop}
{\noindent{\it Proof.}} 
$ (\ref{eq:Q2})+(\ref{eq:Q3}) \ \Rightarrow (\ref{eq:Q4})$.
Consider the inequality $$(S^* \otimes 1_\lambda) \circ (1_\lambda \otimes S)
\circ (S \circ S^*) \circ (1_\lambda \otimes S^*) \circ (S \otimes 1_\lambda)
\leq $$
$$\  \| S \circ S^* \| (S^* \otimes 1_\lambda) \circ (1_\lambda \otimes S) \circ (1_\lambda \otimes S^*) \circ (S \otimes 1_\lambda) = $$
$$ (S^* \otimes 1_\lambda) \circ (1_\lambda \otimes S) \circ (1_\lambda
\otimes S^*) \circ (S \otimes 1_\lambda)$$
(as $\| S \circ S^* \| = 1$). 
The difference $X:= (S^* \otimes 1_\lambda) \circ (1_\lambda \otimes S)
\circ (1_\lambda \otimes S^*) \circ (S \otimes 1_\lambda) - (S^* \otimes
1_\lambda) \circ (1_\lambda \otimes S) \circ (S \circ S^*) \circ (1_\lambda
\otimes S^*) \circ (S \otimes 1_\lambda)$ is a positive element in
$(\lambda,\lambda)$. Checking that ${ \langle 1_\lambda,
X  \rangle }^{(\lambda,\lambda)}_{{\End}(\iota)} = 0$ one concludes that the inequality is actually an
  equality, as the above product is non-degenerate. We put $X':= (1_\lambda \otimes S^*) \circ (S \otimes
  1_\lambda) - S \circ S^*$. Using the preceding result, one checks that $X'^* \circ X' = 0$, which is a
  restatement of (\ref{eq:Q4}).

(\ref{eq:Q2})+(\ref{eq:Q4}) $\ \Rightarrow (\ref{eq:Q3})$
One has: $$(1_\lambda \otimes S) \circ (1_\lambda \otimes S^*) \circ (S \otimes
1_\lambda) = (1_\lambda \otimes S) \circ (S \circ S^*).$$
On the other hand: $$ (1_\lambda \otimes S) \circ (1_\lambda \otimes S^*) \circ
(S \otimes 1_\lambda) = (1_\lambda \otimes (S \circ S^*)) \circ (S \otimes
1_\lambda) =$$ 
$$ (1_\lambda \otimes ((S^* \otimes 1_\lambda) \circ (1_\lambda \otimes S)))
\circ (S \otimes 1_\lambda) = \ ((1_\lambda \otimes S^* \otimes 1_\lambda) \circ (S \otimes S) =$$
$$(1_\lambda \otimes S^* \otimes 1_\lambda) \circ (S \otimes 1_\lambda \otimes
1_\lambda) \circ (1_\lambda \otimes S) = ((S \circ S^*) \otimes 1_\lambda)
\circ (1_\lambda \otimes S) = $$
$$ (S \otimes  1_\lambda) \circ (S^* \otimes 1_\lambda) \circ (1_\lambda \otimes S) =
\ (S \otimes 1_\lambda) \circ (S \circ S^*).$$

Thus $(1_\lambda \otimes S) \circ (S \circ S^*) = (S \otimes 1_\lambda) \circ
(S \circ S^*)$, which is equivalent to (\ref{eq:Q3}), as $S$ is an isometry.

(\ref{eq:Q3})+(\ref{eq:Q4}) $\ \Rightarrow (\ref{eq:Q2})$. Consider the
positive element $H:= S^* \circ S \ \in (\lambda,\lambda)$. We show that $ S
\circ H = \ H \otimes 1_\lambda \circ S = \ 1_\lambda \otimes H \circ S$
holds:
$$H \otimes 1_\lambda \circ S = \ (S^* \circ S) \otimes 1_\lambda \circ S =
S^* \otimes 1_\lambda \circ (S \otimes 1_\lambda \circ S) = \ S^* \otimes
1_\lambda \circ (1_\lambda \otimes S \circ S) = $$
$$  (S^* \otimes
1_\lambda \circ 1_\lambda \otimes S) \circ S = (S \circ S^*) \circ S = S \circ
H;$$
the same way
$$H \otimes 1_\lambda \circ S = \ (S^* \circ S) \otimes 1_\lambda \circ S =
S^* \otimes 1_\lambda \circ (S \otimes 1_\lambda \circ S) = \ S^* \otimes
1_\lambda \circ (1_\lambda \otimes S \circ S) = $$
$$  (S^* \otimes
1_\lambda \circ 1_\lambda \otimes S) \circ S = (1_\lambda \otimes
S^* \circ S \otimes 1_\lambda) \circ S =  $$
$$  1_\lambda \otimes S^* \circ (S \otimes 1_\lambda \circ S) =  1_\lambda
\otimes S^* \circ 1_\lambda \otimes S \circ S = 1_\lambda \otimes H \circ S. $$

Now we show that $H$ is invertible. We have$$ 1_\lambda = S^* \circ T \otimes
1_\lambda \circ T^* \otimes 1_\lambda \circ S \leq  S^* \circ ((T^* \circ T) \otimes
1_\lambda \otimes 1_\lambda) \circ S = (T^* \circ T) \otimes H$$
 where we have used (\ref{eq:Q1}) and the inequality $T \circ T^* \leq (T^* \circ
 T) \otimes 1_\lambda$ which follows by the considerations at the end of the
 proof of Lemma \ref{cupunit}. As $T^* \circ T$ is invertible on $S_l(\lambda)$, we can write $$(T^*
 \circ T)^{-1}_{| _{S_l(\lambda)}} \otimes 1_\lambda \leq H.$$ Thus $H$ is positive and greater or
 equal to a positive invertible element, and this implies that $H$ is
 invertible as well. 

We can define $$S':= H^{- \frac{1}{2}} \otimes 1_\lambda \circ S = \ 1_\lambda
\otimes H^{- \frac{1}{2}} \circ S = \ S \circ H^{- \frac{1}{2}} = \ H^{-
  \frac{1}{2}} \otimes H^{- \frac{1}{2}} \circ S \circ H^{\frac{1}{2}}$$ and
$T':= H^{\frac{1}{2}} \circ T$.
Then $(\lambda,S',T')$ is a $Q$-system isomorphic to $(\lambda,S,T)$ and $S'$ is an isometry.\hfill \halmos
\begin{cor}
\label{frbstrspb}
  Each Frobenius algebra $\lambda$ in a tensor $C^*$-category $\mathcal{T}$ is
  equivalent to a $Q$-system. 
\end{cor}

\begin{cor}

  Each Frobenius algebra $\lambda$ in an $\End(\iota)$-linear tensor $C^*$-category $\mathcal{T}$
  can be embedded in a new tensor $C^*$-category $\mathcal{T'}$ such that the
  embedding is equivalent to a strongly separable Frobenius algebra. 
\end{cor}

\begin{prop}
\label{thmmichaelplus}
  Let $\mathcal{A}$ be a tensor $C^*$-category. The conclusions of Theorem
  \ref{conj2cat} are valid without assuming anyone of the two conditions
  stated in the hypothesis, i.e.\ :

Let $\lambda = (\lambda,T,T^*,S,S^*)$ be a Frobenius algebra in $\mathcal{A}$. 
Then the completion $\mathcal{E} =  \mathcal{E}^P$ of the $\mathcal{E}_0$
defined in Proposition \ref{3.8} is  a bicategory such that
\begin{itemize}
\item $Obj \mathcal{E} = \{ A,B \}$.
\item There is a fully faithful tensor functor $I: \mathcal{A} \rightarrow
  HOM_{\mathcal{E}}(A,A)$ such that for every $Y \in HOM_{\mathcal{E}}(A,A)$
  there is $X \in \mathcal{A}$ such that $Y$ is a retract (i.e.\  sub-$1$-arrow) of $I(X)$. 
\item There are 1-arrows $\rho: B \leftarrow A$ and $ \bar \rho: A \leftarrow B$
  and     2-arrows $$ e_\rho:\iota_B \rightarrow  \rho \otimes
  \bar \rho, \
  \varepsilon_\rho: \iota_A \rightarrow  \bar \rho \otimes  \rho, \ d_\rho:
  \bar \rho \otimes \rho \rightarrow \iota_A, \
  \eta_\rho:  \rho \otimes  \bar \rho \rightarrow \iota_B$$ satisfying the
  conjugation (i.e.\  duality)  relations.
\item We have the identity $$I(\lambda,T,T^*,S,T^*) = (\bar \rho \otimes \rho , e_\rho,\eta_\rho,1_\rho
  \otimes \varepsilon_\rho \otimes 1_{\bar \rho}, 1_\rho \otimes d_\rho
  \otimes 1_{\bar \rho}) $$
of Frobenius algebras in $HOM_{\mathcal{E}}(A,A)$.
\item $\mathcal{E}$ is a
  preadditive 2-category.
\item $\mathcal{E}$ has direct sums of 1-arrows.
\end{itemize}
Isomorphic Frobenius algebras $\lambda, \tilde \lambda$ give rise to isomorphic
bicategories $\mathcal{E}, \tilde{\mathcal{E}}$.

\end{prop}

As a last step in the sequence of propositions dealing with the relationship
between $Q$-systems and 2-categories, we quote the following, which shows that
starting from a $Q$-system one can extend the $*$ operation and, when starting
with a $Q$-system in a tensor $C^*$-category, actually
recover a $2-C^*-$category.
 
\begin{prop}
  (cf.\ \cite{Mu}, Proposition 5.5) Let $\mathcal{A}$ be a tensor *-category
  and  $\lambda$ a $Q$-system in $\mathcal{A}$. Then $\mathcal{E}_0$ has
  a positive *-operation which extends the given one on $\mathcal{A}$. Let
  $\mathcal{E}_*$ be the full sub-bicategory of $\mathcal{E}$ whose 1-arrows
  are $(X,P)$, where $X$ is an object in $\mathcal{A}$ and where $P = P \circ P = P^*$. Then $\mathcal{E}_*$ is equivalent
  to $\mathcal{E}$ with positive involution *. 
\end{prop}

The following proposition answers positively a question left open in
$\cite{LR97}$, namely, whether a couple of conjugate elements $\rho, \bar \rho$
in a tensor $C^*$-category does give rise   to a Frobenius algebra $(\bar \rho
\otimes
\rho, S,T)$ such that $S$ is an isometry also in the case $\End(\iota) \ncong \mathbb{C}$.

\begin{prop}
  Let $ B \xleftarrow {\rho} A$ be a 1-arrow in a 2-$C^*$-category $\mathcal{C}$and let
  $A \xleftarrow{\bar \rho} B$ be a conjugate for $\rho$. Let $\lambda:= \bar
  \rho \otimes \rho$. Then there exist $S, \ T$ such that $(\lambda, S, T)$
  satisfy the defining relations for a 
  $Q$-system (\ref{eq:Q1}-\ref{eq:Q4}) in the tensor $C^*$-category $\mathcal{HOM}(A,A)$ generated by 1-arrows
  connecting the object $A$ to itself.
\end{prop}

{\noindent{\it Proof.}}  As $\rho$ and $\bar \rho$ are conjugate, there exist $R_\rho, \bar R_\rho$
satisfying the conjugation equations. Defining $T':= R_\rho,$ $ S':= (1_{\bar \rho}
\otimes \bar R_\rho \otimes 1_{\rho}) \circ R_\rho$ we see that $T' \in \Hom(\iota, \lambda)$
and $S' \in (\lambda, \lambda \otimes \lambda)$. $(\lambda,S',T')$ satisfy the
relations (\ref{eq:Q1}), (\ref{eq:Q3}), (\ref{eq:Q4}) (we leave to the reader
the easy proof, which relies on the conjugation relations for $\rho, \bar \rho$).
Defining $S:= S' \circ (S^* \circ S)^{-\frac{1}{2}}, \ T:= (S^* \circ
S)^{\frac{1}{2}} \circ T'$, Proposition \ref{qequivalenza} tells us that
$(\lambda, S, T)$ satisfy all  of (\ref{eq:Q1}-\ref{eq:Q4}). \hfill \halmos

\begin{rem}
  A different choice of solutions $R,\bar R$ or of the conjugate $\bar \rho$
  changes the Frobenius algebra only up to isomorphism. 
\end{rem}

\begin{ex}
As $R_\rho^* \circ R_\rho$ and $\bar R_\rho^* \circ \bar
  R_\rho$ are invertible on $S_l(\rho)$ and $S_r(\rho)$ respectively, the elements $$ E_\rho:= R_\rho \circ (R^*_\rho \circ 
  R_\rho)^{-1}_{| _{S_l(\rho)}} 
\circ R^*_\rho, \ \bar E_\rho:=  \bar R_\rho \circ (\bar
  R^*_\rho \circ \bar R_\rho)^{-1}_{| _{S_r(\rho)}} \circ \bar R^*_\rho.$$ are easily seen to
  be projections (where we use the same notation as before, i.e.\  $(R^*_\rho \circ
   R_\rho)^{-1}_{| _{S_l(\rho)}}$ is the element in $ \End(\iota_A)
  \otimes E_{S_l(\rho)}$ s.t. $(R^*_\rho \circ \ R_\rho)^{-1}_{|
  _{S_l(\rho)}} \otimes (R^*_\rho \circ R_\rho) = E_{S_l(\rho)}$ and 
  analogously for $(\bar
  R^*_\rho \circ \bar R_\rho)^{-1}_{| _{S_r(\rho)}}$.)
  If we suppose $\rho$ to be centrally balanced and Assumption \ref{balancedass} to
  hold (thus we have $\End(\iota)$-linearity of the category generated by $ \bar \rho \otimes \rho $)  we can
  also suppose to have chosen $R_\rho$ and $\bar R_\rho$ such that $
  \bar R_\rho^* \circ \bar R_\rho = \theta_\rho(R_\rho^* \circ R_\rho)$ (it
  suffices to renormalise $R_\rho$ and $\bar R_\rho$ by tensoring with
  $(R_\rho^*\circ R_\rho)^{- \frac{1}{4}}_{| _{S_l(\rho)}}  \otimes \theta_\rho^{-1}((\bar
  R_\rho^* \circ \bar R_\rho^*)^{\frac{1}{4}})$ and $(\bar R_\rho^* \circ \bar
  R_\rho )^{- \frac{1}{4}}_{| _{S_r(\rho)}}  \otimes \theta_\rho(R_\rho^* \circ R_\rho)^{\frac{1}{4}})$.

Then the following relations hold (cf.\ \cite{LR97}): 
$$ 1_\rho \otimes (R_\rho \circ R_\rho^*) \circ (\bar R_\rho \circ \bar R_\rho^*) \otimes 1_\rho
\circ 1_\rho \otimes (R_\rho^* \circ R_\rho) = 1_\rho \otimes (R_\rho^* \circ R_\rho) =
\theta_\rho(R_\rho^* \circ R_\rho) \otimes 1_\rho . $$

This implies $$ 1_\rho \otimes E_\rho \circ \bar E_\rho \otimes  1_\rho
\circ 1_\rho \otimes E_\rho = \theta_\rho(R_\rho^* \circ R_\rho)^{-2}_{| _{S_r(\rho)}}  \otimes 1_\rho \otimes E_\rho = (\bar R_\rho^* \circ \bar R_\rho)^{-2}_{| _{S_r(\rho)}}  \otimes 1_\rho \otimes E_\rho.$$

Analogously one obtains 
$$ 1_{\bar \rho} \otimes \bar E_\rho \circ E_\rho \otimes 1_\rho \circ 1_{\bar
  \rho} \otimes \bar E_\rho = (R_\rho^* \circ R_\rho)^{-2}_{| _{S_l(\rho)}}  \otimes 1_{\bar \rho} \otimes
  \bar E_\rho.$$ 

These are the Jones relations which in the case $\End(\iota) = \mathbb{C}$
lead to a representation of the Temperley-Lieb algebra related to the
parameter $R_\rho^* \circ R_\rho$. The difference here is that $R_\rho^* \circ R_\rho$ is,
in general, a positive  function in the  $C^*$-algebra $\End(\iota) \cong
C(\Omega)$. Obviously if the function $R_\rho^* \circ R_\rho$ takes values in the discrete part of
the spectrum of the Jones index, i.e.\  $ (R_\rho^* \circ R_\rho)^2(\omega) \in \{ 4 \cos^2
\pi / k, \ k \in \mathbb{N}, \ k \geq 3 \}$, then it is a locally constant function,
as it has to be continuous on the connected subspaces of $\Omega$.

\end{ex}

\begin{lem}
  Let $\mathcal{A}$ be a tensor $C^*$-category and $(\lambda,S,S^*,T,T^*)$ a Frobenius
  algebra in $\mathcal{A}$. Suppose that Assumption \ref{balancedass} holds 
  in $\mathcal{A}$. Then the tensor $C^*$-category generated by $\lambda$ (i.e.\  its
  tensor powers and their sub-objects) is $\End(\iota)$-linear.
\end{lem}
{\noindent{\it Proof.}} 
It suffices to prove $\End(\iota)$-linearity only for $\lambda$, i.e.\  that
for any $z \in \End(\iota)$, $ z \otimes 1_\lambda = 1_\lambda \otimes z $. The
same relation for powers of $\lambda$ and sub-objects follows immediately.

Let's suppose for the moment  $\lambda$ to be  centrally balanced. By
Assumption \ref{balancedass} there exists a $w \in \End(\iota)$ such that $z
\otimes 1_\lambda = 1_\lambda \otimes w$. We suppose, without loss of
generality, $w,z \in E_{S_l(\lambda)} \otimes \End(\iota).$ Thus we have:   
$$(1_\lambda \otimes T^*) \circ (1_\lambda \otimes z \otimes 1_\lambda) \circ
S = (1_\lambda \otimes T^*) \circ (1_\lambda \otimes 1_\lambda \otimes w) \circ
S.$$  

For the left hand side the following hold: 
$$(1_\lambda \otimes T^*) \circ (1_\lambda \otimes z \otimes 1_\lambda)
\circ S = (1_\lambda \otimes z \otimes T^*) \circ S = (1_\lambda \otimes (z
\circ T^*) \circ S = $$
$$(1_\lambda \otimes z) \circ (1_\lambda \otimes T^*) \circ
S = 1_\lambda \otimes z.$$ 

For the right hand side we have:
$$ (1_\lambda \otimes T^*) \circ (1_\lambda \otimes  1_\lambda \otimes
w) \circ S = ((1_\lambda \otimes T^*) \circ S) \otimes w = 1_\lambda \otimes
w.$$

Thus $1_\lambda \otimes z = 1_\lambda \otimes w$, i.e.\  $z = w$.

Now consider $\lambda = \oplus_i \lambda_i$, where each
$\lambda_i$ is centrally balanced. Then for any $z$ in $\End(\iota)$ we have $ z
\otimes 1_\lambda = z \otimes ( \oplus_i 1_{\lambda_i} ) = \oplus_i (1_{\lambda_i}
\otimes \theta_{\lambda_i}(z)).$  But we have just seen that each $\theta_{\lambda_i}$
is trivial, i.e.\  $\theta_{\lambda_i}(z) = z $ for any $z$ in $ \Hom(\iota,\iota)$.

Thus $z \otimes 1_\lambda = 1_\lambda \otimes z$. \hfill \halmos

\begin{cor}
  Let $\mathcal{A}$ be a tensor $C^*$-category for which Assumption
  \ref{balancedass} holds. Then each Frobenius algebra $\lambda \in \mathcal{A}$
  is an $\End(\iota)$-linear $Q$-system in the tensor $C^*$-category generated by the tensor powers of
  $\lambda$ itself.
\end{cor}

\section{Conclusions}
In all of the present work a fundamental role has been played by the
$C^*$-property of the norm and the conjugation relations. A bundle
structure of the spaces of 2-arrows appears, whereas in preceding works it was
(implicitly or explicitly) a starting hypothesis. The fact that the bundle structure
is preserved by the $\circ$ composition (Proposition \ref{bundleb}) is a direct consequence of the
conjugation equations. In order to give a reasonable description of the
behaviour of this structure under the $\otimes$ composition, we have
introduced an additional hypothesis, i.e.\  Assumption \ref{balancedass}. This has
enabled us to describe our initial 2-$C^*$-category as a collection of
2-$C^*$-categories with simple units indexed by the elements of a compact topological
space (the categories $C^{\omega_0 \mathcal{A}}$ at the end of Section
\ref{bundles}), a structure resembling that of fibre bundles which assumes a
particularly nice form in the case of tensor categories (Proposition
\ref{fibre of tensor categories}). This result suggests that tensor $C^*$-categories with non-simple units might play an  important role for a duality
theory of compact groupoids, as the simple unit case has played for compact groups (\cite{DR89}).       
Properties analogous to the case of simple units have been proven (as, for
example, the finite dimensionality of the fibres), as well continuity of the
sections in this ``fibre picture''.
Assumption \ref{balancedass} seems general enough to handle many interesting
cases. It would be interesting to study to which extent one can consider this
hypothesis as valid, or give some counterexamples.

A second open question is that of the existence of standard solutions. We gave
some partial answers, i.e.\  a ``weak'' positive answer 
(Proposition \ref{weakstandardsolution}) in the general case and a ``global'' positive
answer in the case of locally trivial bundles (Proposition \ref{locallytrivialstandard}).
We don't know the answer for the general case (nor do we have any counterexamples).  

The bundle approach seems to have proven itself fruitful giving an example
of bundles of Hopf algebras in
Section \ref{hopfbundles}.

Finally, the remarks in Section \ref{Frobenius} show that, once more, the
conjugation relations together with the $C^*$-property of the norm may hide more structure than what appears at first glance.

\section*{Acknowledgements}
I would like to thank J. E. Roberts for suggestions, comments, observations
and encouragement: his unpublished handwritten notes were fundamental for the initial
settings of this work. I also would like to thank E. Vasselli for
fruitful discussions and M. M\"uger for comments on previous versions of this
manuscript,  as well as R. Schrader for kind hospitality at the Department
of Physics - Freie
Universit\"at Berlin, where part of this work was written.

\end{document}